\documentclass{amsart}

\usepackage{amsmath, amssymb, amsthm, mathtools, tikz-cd, xcolor}
\usepackage[hidelinks]{hyperref}
\usepackage{environ}
\usepackage{comment}
\usepackage{enumitem}
\theoremstyle{plain}
\newtheorem{thm}{Theorem}
\newtheorem{prop}[thm]{Proposition}

\newtheorem{cor}[thm]{Corollary}
\newtheorem{lem}[thm]{Lemma}

\theoremstyle{definition}
\newtheorem{defi}[thm]{Definition}
\newtheorem{ex}[thm]{Example}
\newtheorem{rmk}[thm]{Remark}

\newcommand{\cF}{\mathcal{F}}
\newcommand{\cG}{\mathcal{G}}
\newcommand{\cK}{\mathcal{K}}
\newcommand{\cO}{\mathcal{O}}

\newcommand{\cZ}{\mathcal{Z}}
\newcommand{\cC}{\mathcal{C}}
\newcommand{\cU}{\mathcal{U}}

\newcommand{\cT}{\mathcal{T}}
\newcommand{\cM}{\mathcal{M}}
\newcommand{\cQ}{\mathcal{Q}}
\newcommand{\bC}{\mathbb{C}}
\newcommand{\bT}{\mathbb{T}}
\newcommand{\bN}{\mathbb{N}}

\newcommand{\into}{\hookrightarrow}

\newcommand{\restrict}[2]{\left.{#1}\right|_{#2}}
\newcommand{\norm}[1]{\left\lVert {#1}\right\rVert}
\newcommand{\abs}[1]{\left| {#1}\right|}
\newcommand{\ndim}[1]{\operatorname{dim}_\mathrm{nuc}({#1})}
\newcommand{\nndim}[1]{\operatorname{dim}_\mathrm{nuc}^\omega({#1})}
\newcommand{\fndim}[1]{\mathrm{f}\mathrm{dim}_\mathrm{nuc}({#1})}
\NewEnviron{equa}{ \begin{equation}\begin{split} \BODY \end{split} \end{equation} }
\usepackage[dvipsnames]{xcolor}

\title{The nuclear dimension of extensions and graph $C^*$-algebras}

\author{Samuel Evington}
\address{Samuel Evington, Mathematical Institute, University of M\"unster, Einsteinstr.\ 62, 48149 M\"unster, Germany}
\email{evington@uni-muenster.de}
\curraddr{School of Maths and Stats, University College Dublin, D04 V1W8, Ireland}

\author{Gregory Faurot}
\address{Gregory Faurot, Department of Mathematics, The Ohio State University, Columbus, OH, USA 43210}
\email{faurot.3@osu.edu}

\author{Christopher Schafhauser}
\address{Christopher Schafhauser, Department of Mathematics, University of Nebraska--Lincoln, Lincoln, NE, USA  68588}
\email{cschafhauser2@unl.edu}

\thanks{SE was partially funded by ERC Advanced Grant 834267 (AMAREC), the German Research Foundation (DFG) under CRC 1442 (Project-ID 427320536), and Germany’s Excellence Strategy EXC 2044/2 (Project-ID 390685587, Mathematics Münster: Dynamics–Geometry–Structure).  GF was partially supported by NSF grant DMS-2154389.  CS was partially supported by NSF grant DMS-2400178.}

\date{\today}


\begin{document}
	
	\begin{abstract}
		We prove that the nuclear dimension of any graph $C^*$-algebra is at most two, and we compute the exact value of the nuclear dimension for a large class of graph $C^*$-algebras. In particular, this class covers all graphs with Condition (K).  Our methods --- based on Szab{\'o}'s recent $KK$-uniqueness theorem --- provide a unified approach to computing the nuclear dimension of various extensions. This includes the nuclear dimension of the Toeplitz algebra (first computed by Brake--Winter) and that of the Cuntz--Toeplitz algebras (first computed by Easo \emph{et al}).
	\end{abstract}
	
	\maketitle
	
	\renewcommand{\thethm}{\Alph{thm}}
	\numberwithin{equation}{section} 
	
	\section{Introduction}\label{sec:intro} 
	Nuclear dimension for $C^*$-algebras was introduced by Winter and Zacharias (\cite{WZ10}). Defined by ``coloring" the finite-dimensional approximations from the Choi--Effros--Kirchberg characterization of nuclearity (\cite{CE76b,Kir77}), 
	nuclear dimension can be viewed as a noncommutative analogue of covering dimension for topological spaces in the sense that $\ndim{C(X)}=\dim(X)$ holds for every compact metrizable space $X$. Indeed, in one direction, if $\dim(X) \leq n$, then, by Ostrand's theorem (\cite{Os65}), every finite open cover of $X$ has a refinement which can be $(n+1)$-colored so that open sets of the same color are disjoint. The finite-dimensional approximations witnessing $\ndim{C(X)} \leq n$ can then be defined using a partition of unity subordinate to the $(n+1)$-colored open cover. 
	
	Finite nuclear dimension has played a major role in the Elliott classification program for simple separable nuclear $C^*$-algebras (\cite{Ell76, Ell95, Elliott-Toms}): 
	counterexamples of Villadsen (\cite{Villadsen1, Villadsen2}), R{\o}rdam (\cite{Rordam}), and Toms (\cite{Toms}) showed that finite nuclear dimension is a necessary condition for a general $K$-theoretic classification theorem, 
	and the state-of-the-art classification theorems include finite nuclear dimension (or an equivalent regularity condition) as a hypothesis (\cite{GLN20A, GLN20B, EGLN23, TWW17, CGSTW23}). 
	
	The study of the nuclear dimension of simple $C^*$-algebras has a long history (\cite{WZ10, MS14, Enders, RSS15, BBSTWW19}), with the definitive result obtained in \cite{W12, CETWW21, CE20}: a simple $C^*$-algebra $A$ has $\ndim{A} \in \{0,1,\infty\}$.  In the separable non-elementary setting, finite nuclear dimension is equivalent to $\mathcal Z$-stability (in the sense of \cite{JiangSu}) and nuclear dimension zero characterizes approximately finite-dimensional (AF) $C^*$-algebras.
	
	The commutative case shows that the nuclear dimension trichotomy that occurs in the simple setting does not generalize to all $C^*$-algebras, yet there is evidence that this extends well beyond the simple case; see \cite[Section~23]{STW26} for a recent survey.  For example, \cite[Problem~LXXXI]{STW26} asks if the trichotomy holds for all $C^*$-algebras for which no quotient has a commutative hereditary algebra $C_0(X)$ with $\dim(X) \geq 2$.  A definitive answer in the purely infinite case, predating \cite{STW26}, is given in the work of Bosa, Gabe, Sims, and White in \cite{BGSW22}, showing separable nuclear $\mathcal O_\infty$-stable $C^*$-algebras have nuclear dimension at most one.
	
	Since their introduction by Kumjian, Pask, and Raeburn (\cite{KPR98}), the study of graph $C^*$-algebras has grown into a significant field of operator algebras in its own right (\cite{R05}). 
	In particular, these $C^*$-algebras have proven to be a useful testing ground for $C^*$-algebraic conjectures as their internal structure can often be understood in terms of the underlying graph.  
	For example, graph algebra results played a prominent role in computing the nuclear dimension of Kirchberg algebras. A series of works from various authors (\cite{WZ10, MS14, Enders}) proved that the nuclear dimension of Kirchberg algebras was at most three, with an upper bound of one under an assumption of torsion-free $K_1$-group. Ruiz, Sims, and S{\o}rensen (\cite{RSS15}) proved all UCT Kirchberg algebra have nuclear dimension at most one by constructing them as a direct limits of $2$-graph C$^*$-algebras. This played a major role in motivating the definitive result of \cite{BBSTWW19}, where the UCT assumption was able to be dropped using other methods. 
	
	Since the only elementary subquotients that can occur in a graph $C^*$-algebra have one-dimensional spectrum, it has been conjectured that all graph $C^*$-algebras have nuclear dimension at most one (with zero occurring if and only if the graph has no cycles; equivalently, the corresponding $C^*$-algebra is AF); this is recorded as \cite[Problem~LXXXIII]{STW26}. 
	
	The first main result of this paper is a partial solution to this problem.  Note that, for graph $C^*$-algebras, we follow the edge convention of Raeburn's monograph (\cite{R05}). 
	
	\begin{thm}[Theorem~\ref{thm:graphs}]\label{intro-thm:graphs-with-dim-one} 
		If $E$ is a directed graph such that there is no path from an isolated cycle to another cycle (e.g.\ if $E$ has Condition~(K)) then 
		\begin{equation}
			\ndim{C^*(E))} = \begin{cases} 0 & $E$ \text{ has no cycles}; \\ 1 & $E$ \text{  has a cycle}. \end{cases}
		\end{equation}
	\end{thm}
	
	Loosely speaking, in a graph $E$, \emph{return paths} are paths from a vertex $v$ to itself that do not otherwise visit $v$; \emph{cycles} are return paths where $v$ is the only vertex visited more than once; and a cycle in a graph is \emph{isolated} if there is no other return path based on this cycle. With this terminology, a graph has \emph{Condition~(K)} if and only if it has no isolated cycles (see Definition~\ref{def:isolatedcycles} for the formal definitions). 
	Theorem~\ref{intro-thm:graphs-with-dim-one} improves on a previous result of the second- and third-named authors, where it was shown that all $C^*$-algebras of graphs with Condition~(K) have nuclear dimension at most two (\cite{FC24}).
	
	Our methods also show that all graph $C^*$-algebras have finite nuclear dimension, and in fact, we obtain a common, near-optimal upper bound.  A slightly weaker upper bound was independently obtained by J.\ Wu and V.\ Wu using different methods.
	
	\begin{thm}[{Theorem~\ref{thm:graphs}}]\label{intro-thm:dim-at-most-two}
		If $E$ is a directed graph, then $\ndim{C^*(E)}\leq 2$.
	\end{thm}
	
	Previously, Ruiz, Sims, and Tomforde proved that all purely infinite graph $C^*$-algebras have nuclear dimension at most two, and produced bounds of one or two when the graph algebra is an extension of an AF algebra by a purely infinite ideal (\cite{RST15}). Moreover, the second- and third-named authors proved that all graphs with Condition~(K) have graph $C^*$-algebras whose nuclear dimension is at most two (\cite{FC24}), and reduced this bound to one in certain cases using a result of the first-named author (\cite{E22}). 
	
	Our strategy for proving Theorem~\ref{intro-thm:graphs-with-dim-one} and Theorem~\ref{intro-thm:dim-at-most-two} centers around understanding how nuclear dimension behaves under extensions. 
	Indeed, every graph $C^*$-algebra can be constructed from $\mathbb C$, $C(\mathbb T)$, and Kirchberg algebras using the operations of stable isomorphism, extension, and direct limit. 
	Furthermore, all of these building blocks have nuclear dimension at most one, stable isomorphism preserves nuclear dimension, and nuclear dimension is lower semicontinuous under inductive limits (\cite{WZ10, RSS15, BBSTWW19}). 
	
	Given an extension of $C^*$-algebras 
	\begin{equation}\label{eq:intro-extension}
		\begin{tikzcd}
			0 \arrow{r} & I \arrow{r} & A \arrow{r} & B \arrow{r} & 0,
		\end{tikzcd}
	\end{equation}
	Winter and Zacharias proved the following bounds on the nuclear dimension of $A$ in terms of the nuclear dimensions of $I$ and $B$ (\cite{WZ10}):
	\begin{equation}
		\max\{\ndim{I}, \ndim{B} \} \leq \ndim{A} \leq \ndim{I} + \ndim{B} + 1,
	\end{equation}
	and these are still the best bounds that are known to hold generally.
	
	When $A$ is separable and abelian, one can do much better. Translating from nuclear dimension to covering dimension of spaces, the problem changes to computing the dimension of a metric space from the dimension of an open set and its complement. A classical result in topology (see, for example, \cite[Corollary 5.8]{Pears75book}) then proves that
	\begin{equation}\label{eq:optimal}
		\ndim{A} = \max\{ \ndim{I}, \ndim{B} \},
	\end{equation}
	i.e.\ the lower bound in the Winter--Zacharias estimate is attained.   
	It is an open question whether \eqref{eq:optimal} holds generally; see \cite[Problem~LXXXII]{STW26}, for example.
	
	It is straightforward to show that \eqref{eq:optimal} also holds when the extension is a direct sum, and this result generalizes to ``approximate'' direct sums in the sense of quasidiagonal extensions; see \cite[Proposition 6.3]{RST15}. 
	The exact value of the nuclear dimension of several further interesting examples of extensions has been computed in recent years, and in all these cases \eqref{eq:optimal} has been verified.
	
	Brake and Winter showed that the Toeplitz algebra $\cT$ -- an extension of $C(\bT)$ by the compact operators $\mathcal K$ -- has nuclear dimension one (\cite{BW19}); 
	Gardner and Tikuisis subsequently generalized this result showing that \eqref{eq:optimal} holds for unital essential extensions of $C(X)$ by $\mathcal K$, where $X$ is a compact metric space (\cite{GT22}); 
	a large team working at the University of Glasgow (including the first-named author) showed that the Cuntz--Toeplitz algebras $\cT_n$ -- an extension of the Cuntz algebra $\cO_n$ by  $\mathcal K$ -- has nuclear dimension one (\cite{15A20}) by combining the Brake--Winter technique with classification machinery for $\cO_\infty$-stable maps (\cite{Ki95,Ph00,Gabe2}); 
	this was extended to unital full extensions of unital separable nuclear $\mathcal O_\infty$-stable $C^*$-algebras $B$ by stable approximately finite-dimensional algebras $I$ by the first-named author (\cite{E22}).  Recently,  the first-named author, Ng, Sims, and White proved that \eqref{eq:optimal} holds for all essential extensions of unital separable commutative $C^*$-algebras $B$ by stable Kirchberg algebras $I$ (\cite{ENSW25}).
	
	The new techniques introduced in this paper provide a common framework for these results on the nuclear dimension of extensions.
	We consider full extensions of the form \eqref{eq:intro-extension}, in the case where $I$ is stable and $B$ satisfies a technical condition (discussed further below). 
	We then obtain the near-optimal upper bound
	\begin{equation}\label{eq:near-optimal-upper-bound}
		\ndim{A} \leq \max \{ \ndim{I} + 1, \ndim{B} \}.
	\end{equation}
	Furthermore, when $I$ is a stable Kirchberg algebra, we can eliminate the +1 in \eqref{eq:near-optimal-upper-bound} using classification theory in a similar manner to  \cite{ENSW25}. In this case, we therefore obtain \eqref{eq:optimal}.
	
	For clarity of exposition, we shall introduce a variation of nuclear dimension, which we shall call \emph{faithful nuclear dimension} and denote by $\fndim{\,\cdot\,}$. It is defined by an additional injectivity condition on one of the ``colors'' in the finite-dimensional approximation (see Definition \ref{def:fndim} below).  By definition, $\ndim{B} \leq \fndim{B}$ always holds, but the inequality can be strict in general (see Examples~\ref{ex:nndim:0} and \ref{ex:nndim:finitedimensional}). 
	Nevertheless, $\fndim{B}= 1$ holds if $B$ is a separable simple nuclear $\mathcal Z$-stable $C^*$-algebra (Example~\ref{ex:nndim:simplezstab}) or a separable nuclear $\cO_\infty$-stable $C^*$-algebra (Example~\ref{ex:nndim:purelyinfinite}). Further, if $X$ is a connected locally connected compact metrizable space, then $\fndim{C(X)} = \dim(X)$ (Theorem~\ref{thm:nndim:abelian}).
	
	We can now give a more formal statement of our main technical result on the nuclear dimension of extensions.
	
	\begin{thm}[{Theorems~\ref{body-thm:extension-result-1}~and~\ref{body-thm:extension-result-2}}]\label{intro-thm:extension-result}
		Suppose that $0 \to I \to A \to B \to 0$
		is a full extension of separable $C^*$-algebras and $I$ is stable. Then
		\begin{equation}\label{eq:extension-result}
			\ndim{A}\leq \max\{\ndim{I}+1, \fndim{B}\}.
		\end{equation}
		Further, if $I$ is a stable Kirchberg algebra, then
		\begin{equation}
			\ndim{A} \leq \max\{ \ndim{I}, \fndim{B}\}.
		\end{equation}
	\end{thm}
	
	Theorem~\ref{intro-thm:extension-result} recovers the results of \cite{BW19, 15A20, E22} because replacing $\fndim{B}$ with $\ndim{B}$ has no effect on the value of the right-hand side in these special cases. 
	Moreover, Theorem~\ref{intro-thm:extension-result}  recovers the results of \cite{GT22,ENSW25} under mild restrictions on the allowed spectrum of the commutative $C^*$-algebra (see Theorem~\ref{thm:nndim:abelian}). 
	A strengthened version of Theorem~\ref{intro-thm:extension-result} (Theorem~\ref{thm:superniceapproximations}) is required to facilitate an induction through the ideal lattice of the $C^*$-algebra associated to a finite graph. This, together with a combinatorial reduction argument to the finite case, leads to the proof of Theorems \ref{intro-thm:graphs-with-dim-one} and \ref{intro-thm:dim-at-most-two}
	(cf.\ \cite[Section~1]{RS04}, \cite[Proposition~3.1]{JP02}, and \cite[Section~2]{FC24}).

	\subsection*{Discussion of methods}
	
	The proof of Theorem~\ref{intro-thm:extension-result} (and its variants that are needed to obtain Theorems~\ref{intro-thm:graphs-with-dim-one} and \ref{intro-thm:dim-at-most-two}) is rather technical, so we take some space here to outline the general ideas involved. In this discussion, we assume the reader has some familiarity with the theory of completely positive maps, order zero maps, and the definition of nuclear dimension (see for example \cite{WZ09} and \cite{WZ10}).
	
	Consider an extension
	\begin{equation}
		\begin{tikzcd}
			0 \arrow{r} & I \arrow{r} & A \arrow[r, "q"] & B \arrow{r} & 0. 
		\end{tikzcd}	
	\end{equation}
	The basic setup is similar to that of earlier results on nuclear dimension of extensions from \cite{BW19,GT22,15A20,E22,ENSW25}, although we keep more of the proof at the level of ultrapowers, which we believe helps remove some notational clutter and highlights the use of classification techniques.
	
	Suppose $\fndim{B} < \infty$. The faithful nuclear dimension approximations for $B$ induce completely positive contractive order zero (cpc$_\perp$) maps $\psi^0, \psi^1 \colon B \rightarrow B_\omega$ with $\psi^0(b) + \psi^1(b) = b$ for all $b \in B$. These approximations are also constructed such that $\psi^0$ and $\psi^1$ lift to sequences of completely positive (cp) maps $B \rightarrow B$ of nuclear dimension 0 and (at most) $\fndim{B} - 1$, respectively. For brevity, we denote this by $\nndim{\psi^0} = 0$ and $\nndim{\psi^1} \leq \fndim{B}-1$.  The difference between the definition of $\fndim{B} \leq n$ and that of $\ndim{B} \leq n$ is precisely the following: we may further arrange for $\psi^0$ to be \emph{totally faithful} in the sense that the induced $^*$-homomorphism $\Psi^0 \colon C_0(0, 1] \otimes B \rightarrow B_\omega$ from \cite{WZ09} is faithful.
	
	Let $\mu \colon B \rightarrow I$ be a completely positive contractive (cpc) splitting of $q$.  A now fairly standard use of quasicentral approximate units, in the form of Kirchberg and R{\o}rdam's ``$\sigma$-ideals'' (\cite{KR14}), implies that for any separable set $S \subseteq A_\omega$, there are positive contractions $f, g \in I_\omega \cap S'$ with $gf = f$ and $fx = x$ for all $x \in I_\omega \cap S$ (see Proposition~\ref{prop:qcapproxunit}).
	After applying this to a sufficiently large separable set $S$, define $\phi^0, \phi^1 \colon A \rightarrow A_\omega$ by
	\begin{align}
		\phi^0(a) &= (g-f) a + (1 - g) \mu_\omega(\psi^0(q(a))) \\
		\shortintertext{and}
		\phi^1(a) &= fa + (1 - g) \mu_\omega(\psi^1(q(a)))
	\end{align}
	It is easy to see $\phi^0$ and $\phi^1$ are cpc$_\perp$ and that $\phi^0(a) + \phi^1(a) = a$. Since $(1-g)\mu_\omega(\,\cdot\,)$ is cpc$_\perp$, and $q$ is a $^*$-homomorphism, the second terms of $\phi^i$ have nuclear dimension at most $\ndim{\psi^i}$; with a careful choice of $g$ (Lemma~\ref{lem:constructing-g}), we obtain the same bound in terms of $\nndim{\,\cdot\,}$. As $f$ and $(1-g)$ are orthogonal, we have $\nndim{\phi^1} \leq  \max\{\ndim{I}, \fndim{B} - 1\}$.  It thus suffices to show $\nndim{\phi^0} = 0$. 
	
	Let $\tilde\phi^0 \colon A \rightarrow A_\omega$ denote the second term in the definition of $\phi^0$, i.e.\
	\begin{equation}
		\tilde\phi^0(a) = (1 -g) \mu_\omega (\psi^0(q(a))).
	\end{equation}
	Again, if $g$ is chosen with care (Lemma~\ref{lem:constructing-g}), we have $\nndim{\tilde\phi^0}  = 0$.  In \cite{15A20} and \cite{E22}, the quotient $B$ is assumed to be $\mathcal O_\infty$ stable, and then Gabe's generalization of the Kirchberg--Phillips theorem (\cite{Gabe2}) is used to show $\phi^0$ is unitarily equivalent to $\tilde\phi^0$, which implies $\nndim{\phi^0} = 0$. Therefore,
	\begin{equation}
		\begin{split}
			\ndim{A} \leq \nndim{\phi^0} + \nndim{\phi^1} +1\\
			\leq \max\{\ndim{I} + 1, \fndim{B}\}.
		\end{split}
	\end{equation}
	In our case, the underlying classification theorem is Szab\'o's recent $KK$-uniqueness theorem, giving a strong classification theorem for maps into the multiplier algebra of (separable subalgebras of) $I_\omega$.  In the spirit of the lifting techniques in \cite{CGSTW23, CGSTW26}, we obtain the following uniqueness theorem for cpc$_\perp$ lifts of $\psi^0$, which gives the desired unitary equivalence between $\phi^0$ and $\tilde\phi^0$. 
	
	\begin{thm}[{cf.\ Section~\ref{sec:proof-of-thm-D}}]\label{thm:D}
		Let $0 \rightarrow I \rightarrow A \xrightarrow q B \rightarrow 0$ be a full extension of separable $C^*$-algebras such that $I$ is stable and has finite nuclear dimension.  Suppose $D$ is a separable nuclear $C^*$-algebra and $\psi \colon D \rightarrow B_\omega$ is a cpc$_\perp$ map inducing a faithful $^*$-homomorphism $\Psi \colon C_0(0, 1] \otimes D \rightarrow B_\omega$.  If $\tilde\phi, \phi \colon D \rightarrow A_\omega$ are completely positive, contractive, order zero maps with $q_\omega \circ \tilde \phi = q_\omega \circ \phi = \psi$, then $\tilde\phi$ and $\phi$ are unitarily equivalent.
	\end{thm}
	
	When $I$ is a stable Kirchberg algebra, the stronger bound in Theorem~\ref{intro-thm:extension-result} is obtained in a similar way.  In addition to $f$ and $g$, we may arrange for an element $e \in I_\omega \cap A'$ with $e x = x$ for all $x \in I$ and $fe = e$.  The proof from \cite{BGSW22} that Kirchberg algebras have nuclear dimension one gives cpc$_\perp$ maps $\theta^0, \theta^1 \colon  A \rightarrow I_\omega$ with nuclear dimension zero satisfying $\theta^0(a) + \theta^1(a) = ea$ for all $a \in A$. Then $\phi^0, \phi^1 \colon A \rightarrow A_\omega$ are defined by 
	\begin{align}
		\label{eq:intro:approximation1}\phi^0(a) = \theta^0(a) + (g - f)a + (1 - g) \mu_\omega(\psi^0(q(a))) \\
		\shortintertext{and}
		\label{eq:intro:approximation2}\phi^1(a) = \theta^1(a) + (f - e)a + (1 - g) \mu_\omega(\psi^1(q(a))).
	\end{align}
	Theorem~\ref{thm:D} is used to absorb the middle term of $\phi^0$ into the third term, and the Kirchberg--Phillips theorem (in the form of \cite{Gabe2}) is used to absorb the middle term of $\phi^1$ into the first term.  Using the orthogonality of the first and third terms in both $\phi^0$ and $\phi^1$ then gives $\ndim{A} \leq \max\{\ndim{I},\fndim{B}\} = \fndim{B}$ (recalling that $\fndim{B}\geq 1$ by Example~\ref{ex:nndim:0}).  
	
	Returning to graph $C^*$-algebras and the proof of Theorem~\ref{intro-thm:dim-at-most-two}, combinatorial methods first reduce the problem to the case of finite graphs.  If $E$ is a finite graph, then $C^*(E)$ has finitely many gauge-invariant ideal, allowing us to prove Theorem~\ref{intro-thm:dim-at-most-two} by inducting on the number of ideals.  In the gauge-simple case, $C^*(E)$ is one of $M_n$, $C(\mathbb T)\otimes M_n$, or a unital Kirchberg algebras, all of which have nuclear dimension at most 1.  In general, for a minimal gauge-invariant ideal $I$ of $C^*(E)$, there is an extension
	\begin{equation}
		\begin{tikzcd}
			0 \arrow{r} & I \arrow{r} & C^*(E) \arrow{r} & C^*(F) \arrow{r} & 0,
		\end{tikzcd}
	\end{equation}
	where $F$ is a suitable subgraph of $E$.  Further, $I$ is either a direct summand of $C^*(E)$ or is stable.  The nuclear dimension of a direct sum is well-known to be the maximum of the nuclear dimension of the summands, so this reduces the result to the case that $I$ is stable.  Theorem~\ref{intro-thm:extension-result} does not apply directly since the extension need not be full, but the induced extension
	\begin{equation}\label{eq:intro-quotient-extension}
		\begin{tikzcd}
			0 \arrow{r} & I \arrow{r} & C^*(E)/I^\perp \arrow{r} & C^*(F)/I^\perp \arrow{r} & 0,
		\end{tikzcd}
	\end{equation}	
	is full.  In the case $I = \mathcal K$ or $I$ is a stable Kirchberg algebra, this follows from the simplicity of the corona algebra $\mathcal Q(I)$.  In the remaining case, when $I = \mathcal C(\mathbb T) \otimes \mathcal K$, fullness is obtained via a more careful analysis of the graph and an absorption theorem of Pimsner, Popa, and Voiculescu (\cite{PPV79}).
	
	A refinement of the proof of Theorem~\ref{intro-thm:extension-result} gives a version of this theorem for extensions which are full modulo the annihilator in this sense.  This generalization requires a somewhat delicate construction since we need to know nuclear dimension approximations have one color $C^*(F) \rightarrow C^*(F)_\omega$ which is totally faithful modulo $I^\perp$ in the sense that it induces a cpc$_\perp$ map $C^*(F)/I^\perp \to C^*(F)_\omega/(I^\perp)_\omega$ with the induced $^*$-homomorphism on the cone being faithful. In order to facilitate the induction, we assume one color of the nuclear dimension approximations for $C^*(F)$ is ``residually totally faithful'' in the sense that this property holds modulo every gauge-invariant ideal of $C^*(F)$, and we show $C^*(E)$ has the same property. The approximations are still of the form in \eqref{eq:intro:approximation1} and \eqref{eq:intro:approximation2}.
	
	In the case $I$ is $C(\mathbb T)\otimes \cK$ or a stable Kirchberg algebra, the extension result is similar to the proof of Theorem~\ref{intro-thm:extension-result} discussed above.  When $I = \mathcal K$, we are unable to maintain the faithfulness conditions of the approximations for $C^*(E)$, but a combinatorial argument at the level of the graph allows us to avoid this issue by first showing the ``sourceless'' part of the graph has faithful nuclear dimension at most 2. Then inducting through the sources, using that the annihilator of each source projection contains the other source projections, allows for the requisite faithfulness condition on \eqref{eq:intro-quotient-extension} to be satisfied at each stage.
	
	At each step of the induction we have
	\begin{equation}
		\ndim{C^*(E)} \leq \max \{\ndim{I} +1, \ndim{C^*(F)}\}.
	\end{equation}
	Since $\ndim{I} \leq 1$ and the inductive hypotheses gives $\ndim{C^*(F)} \leq 2$, we have $\ndim{C^*(E)} \leq 2$, proving Theorem~\ref{intro-thm:dim-at-most-two}. The proof of Theorem~\ref{intro-thm:graphs-with-dim-one} is essentially identical, using the combinatorial condition on the graph to ensure that we never have $I = C(\mathbb T) \otimes \mathcal K$ at any stage in the induction, and using the stronger version of Theorem~\ref{intro-thm:extension-result} in the case $I$ is a stable Kirchberg algebra in order to avoid increasing the nuclear dimension by one.
	
	The following two graphs, which were highlighted at the end of \cite{ENSW25}, remain the two smallest examples where we cannot compute the nuclear dimension of their graph $C^*$-algebras.
	\begin{equation}
		\begin{tikzcd}
			\bullet \arrow[out=230, in=130, distance=3.5em] \arrow[r] & \bullet \arrow[out=310, in=50, distance=3.5em]
			& & & & & 
			\bullet \arrow[out=230, in=130, distance=3.5em] \arrow[r] &  \bullet \arrow[out=230, in=310, distance=3.5em] \arrow[out=130, in=50, distance=3.5em]
		\end{tikzcd}
	\end{equation}
	As suggested in \cite{ENSW25}, computing the nuclear dimension of the $C^*$-algebra associated to either of these graphs would likely lead to a full solution to the nuclear dimension of graph $C^*$-algebras and would be a major step towards the general problem of computing nuclear dimension of extensions.  In both cases, the difficulty lies in removing the ``$+1$'' in \eqref{eq:extension-result} in the case when $I = C(\mathbb T) \otimes \mathcal K$.
	
	\subsection*{Acknowledgments}
	We thank Stuart White for several helpful conversations over the course of this project.  We are also grateful to Jos\'e Carri\'on, Jamie Gabe, Aaron Tikuisis, and Stuart White for allowing us to include some results on purely large extensions from \cite{CGSTW26}, used in Section~\ref{sec:proof-of-thm-D}, and to Jianchao Wu and Victor Wu for informing us about their forthcoming result on graph $C^*$-algebras having finite nuclear dimension, mentioned above Theorem~\ref{intro-thm:dim-at-most-two}.
	
	The authors would like to thank the Isaac Newton Institute for Mathematical Sciences, Cambridge, for support and hospitality during the program Topological groupoids and their $C^*$-algebras, where work on this paper was undertaken. This work was supported by EPSRC grant EP/V521929/1.  CS acknowledges the support of the Institut Henri Poincar\'e (UAR 839 CNRS-Sorbonne Universit\'e), and LabEx CARMIN (ANR-10-LABX-59-01).
	
	\subsection*{AI Statement}  No AI tools were used to obtain the results of this paper or in the preparation of the manuscript.
	
	\setcounter{thm}{0}
	\numberwithin{thm}{section}
	
	\section{Faithful nuclear dimension}
	
	As is standard, we abbreviate completely positive and completely positive contractive maps as cp and cpc maps, respectively.  For $C^*$-algebras $A$ and $B$, a cpc map has \emph{order zero} (together written cpc$_\perp$) if $\phi(a)\phi(b)=0$ whenever $a, b \in A_+$ with $ab=0$. Of particular importance for us is the following structure theorem for order zero maps.
	
	\begin{thm}[{\cite[Corollary~4.1]{WZ09}}]\label{thm:orderzerodilations}
		Let $\phi\colon A \to B$ be a cpc$_\perp$ map between $C^*$-algebras $A$ and $B$. Then $z \otimes a \mapsto \phi(a)$ uniquely extends to a $^*$-homomorphism $\Phi\colon C_0(0,1] \otimes A \to B$, where $z$ is the canonical inclusion $(0,1] \into \bC$.
	\end{thm}
	
	\begin{rmk}\label{remark:multiplicationoncone}
		The map $\Phi$ in Theorem~\ref{thm:orderzerodilations} may be described more explicitly in many of the cases we will be interested in, including $\cO_\infty$-stable $C^*$-algebras (Example~\ref{ex:nndim:purelyinfinite}) and abelian $C^*$-algebras (Theorem~\ref{thm:nndim:abelian}). If $\phi\colon A \to A_\omega$ is given by $\phi(a)=ha$ for a positive contraction $h \in A_\omega \cap A'$, then $\Phi(\xi \otimes a)=\xi(h) a$. When $\xi(z)=z^n$ and $a\geq 0$, this holds because
		\begin{equation}
			\Phi(z^n \otimes a)=\Phi(z\otimes a^{1/n})^n=(ha^{1/n})^n=h^na,
		\end{equation}
		and the result follows from linearity and continuity of $\Phi$.
	\end{rmk}
	
	We will typically use lower-case Greek letters, such as $\phi$, $\psi$, and $\theta$, for cpc$_\perp$ maps and the corresponding upper-case letters, such as $\Phi$, $\Psi$, and $\Theta$, for the induced $^*$-homomorphisms as in the previous theorem. Our theorems will require certain order zero maps to be faithful after dilating to the cone, resulting in the following definition.
	
	\begin{defi}
		A cpc$_\perp$ map $\phi$ is said to be \emph{totally faithful} if the corresponding $^*$-homomorphism $\Phi$ is faithful.
	\end{defi}
	
	We will need nuclear dimension for both $^*$-homomorphisms and for cp maps.  We recall the definition here.
	
	\begin{defi}[{\cite[Definition~2.1]{WZ10}}]\label{def:ndim}
		Let $A$ and $B$ be $C^*$-algebras and let $\phi \colon A \to B$ be a cp map. The nuclear dimension of $\phi$, written $\ndim{\phi}$, is the smallest natural number $n \in \bN$ such that, for all finite subsets $\cF \subset A$ and $\epsilon >0$, there exist finite-dimensional $C^*$-algebras $F^0, F^1,\dots, F^n$ and cp maps
		\begin{equation}
			\begin{tikzcd}
				A \arrow{r}{\theta} & \bigoplus_{i=0}^n F^i \arrow{r}{\psi} &  B
			\end{tikzcd}
		\end{equation}
		such that
		\begin{enumerate}
			\item $\norm{\phi(a)-\psi(\theta(a))}<\epsilon$ for all $a \in \cF$,
			\item $\theta$ is contractive, and
			\item $\psi^i = \restrict{\psi}{F^i}$ is cpc$_\perp$ for $i=0, 1, \dots, n$.
		\end{enumerate}
		Iif no such $n$ exists, we write $\ndim{\phi}=\infty$. For a $C^*$-algebra $A$, $\ndim{A}$ is defined to be $\ndim{\mathrm{id}_A}$.
	\end{defi}
	
	We record the following simple fact for later use.  The analogous statement for algebras (i.e.\ when $\phi = \mathrm{id}_A$) is a special case of \cite[Proposition~2.3(ii)]{WZ10}.
	
	\begin{lem}\label{lem:nndim-stable}
		Suppose $A$ and $B$ are $C^*$-algebras and $\phi \colon A \rightarrow B$ is a cp map.  For an AF algebra $D$, $\ndim{\phi \otimes{\mathrm{id}_D}} \leq \ndim{\phi}$.
	\end{lem}
	
	\begin{proof}
		We may assume $n = \ndim{\phi} < \infty$. Let $\mathcal F \subseteq A$ and $\mathcal G \subseteq D$ be finite sets.  Let $\theta$ and $\psi$ and $F^0, \ldots, F^n$ be as in Definition~\ref{def:ndim}, and let $G \subseteq D$ be a finite dimensional C$^*$-algebra so that for every $c \in \mathcal G$, there is $d \in G$ such that $\|c - d\| < \epsilon$.  Let $E \colon D \rightarrow G$ be a conditional expectation and let $\iota \colon G \rightarrow D$ be the inclusion.  Set $\tilde F^i = F^i \otimes G$, $\tilde \theta = \theta \otimes E$, and $\tilde \psi = \psi \otimes \iota$.  Then $\tilde F^i$ is finite dimensional, $\tilde \theta$ is cpc, and $\tilde \psi|_{\tilde F^i}$ is cpc$_\perp$ for all $i = 0, \dots, n$ by \cite[Theorem~4.3]{WZ09}.  Further, for $a \in \mathcal F$ and $b \in \mathcal G$, we have 
		\begin{equation}
			\begin{split}
				\|(\phi \otimes \mathrm{id}_D)(a \otimes b) &- \tilde\psi(\tilde \theta(a \otimes b)\| \\
				&= \|\phi(a) \otimes b - \psi(\theta(a)) \otimes E(b) \| \\ 
				&\leq \|(\phi(a) \otimes (b-E(B)) \|+\|(\psi(\theta(a))-\phi(a))\otimes E(b) \|\\
				&< \|\phi\|\|a\|\epsilon + \|b\| \epsilon.
			\end{split}
		\end{equation}
		So $\ndim{\phi \otimes \mathrm{id}_D} \leq \ndim{\phi}$.
	\end{proof}

	Our next goal is to introduce the variant of nuclear dimension needed in the statement and proof of Theorem~\ref{intro-thm:extension-result}.  For this, we will need the theory of ultrapowers of $C^*$-algebras and maps, which we recall briefly here.  Let $\omega$ be a free ultrafilter on the natural numbers. For a $C^*$-algebra $A$, we will write $A_\omega$ for the $C^*$-algebra $\ell^\infty(A)/c_\omega(A)$, the algebra of bounded sequences in $A$ modulo the ideal of null sequences (where convergence is taken along the ultrafilter $\omega$). We will identify $A$ as a subalgebra of $A_\omega$ via constant sequences.  Note that a uniformly bounded sequence of cp maps $(\phi_n \colon A \rightarrow B)_{n=1}^\infty$ between $C^*$-algebras $A$ and $B$ induces a cp map $\phi_\omega \colon A \rightarrow B_\omega$.  Since we do not know if 
	$\ndim{\phi_\omega} \leq \sup_n \ndim{\phi_n}$, 
	we introduce the following notation.
	
	\begin{defi}\label{def:nndim}
		Let $A$ and $B$ be $C^*$-algebras and let $n \in \bN$.  For a cp map $\phi \colon A \rightarrow B_\omega$, write $\nndim{\phi} \leq n$ if $\phi$ lifts to a sequence of cp maps $(\phi_k \colon A \rightarrow B)_{k=1}^\infty$ such that $\ndim{\phi_k} \leq n$ for all $k \in \bN$.
	\end{defi}
	Note that if $A$ is nuclear, every cp map $\phi\colon A \to B_\omega$ has $\nndim{\phi}\leq \ndim{B}$. Nuclearity of $A$ is sufficient to lift $\phi$ to a sequence of cp maps $(\phi_k\colon A \to B)_{k=1}^\infty$ by the Choi--Effros theorem (\cite{CE76a}), and each of these maps will have their nuclear dimension bounded by $\ndim{B}$.

    In some previous articles (e.g. \cite{CETWW21, BGSW22}), the nuclear dimension of $\phi\colon A \to B_\omega$ was computed rather than the nuclear dimension of lifts of $\phi$. In order to allow us to focus on $\nndim{\,\cdot\,}$, we point out the following proposition.

    \begin{prop}\label{prop:nndim<ndim}
        Let $\phi\colon A \to B_\omega$ be a cp map with $A$ separable. Then $\nndim{\phi}\leq \ndim{\phi}$.
    \end{prop}

    \begin{proof}
        If $\ndim{\phi}=\infty$, there is nothing to show. Thus, we may assume that $\ndim{\phi}<\infty$. The result then follows from the projectivity of cones over finite-dimensional $C^*$-algebras (\cite{LP98}) and an application of Kirchberg's $\epsilon$-test (\cite[Lemma~3.1]{KR14}).
    \end{proof}
	
	As with all other results computing nuclear dimension of extensions, the nuclear dimension approximations will be given by adding approximations built from those of the ideal and the quotient.  To this end, we record how $\nndim{\,\cdot\,}$ behaves under sums of cp maps. 
	
	\begin{lem}\label{lem:sum}
		Suppose $A$ and $B$ are $C^*$-algebras such that $A$ is separable and let $\phi^0, \phi^1 \colon A \rightarrow B_\omega$ be cp maps.  Then 
		\begin{equation}
			\nndim{\phi^0 + \phi^1} \leq \nndim{\phi^0} + \nndim{\phi^1} + 1.
		\end{equation}
	\end{lem}
	
	\begin{proof}
		Let $\nndim{\phi^0}=m$ and $\nndim{\phi^1}=n$ and, without loss of generality, we may assume that $m\leq n<\infty$. Lift $\phi^0$ and $\phi^1$ to sequences of cp maps $(\phi^0_k\colon A \to B)_{k=1}^\infty$ and $(\phi^1_k\colon A \to B)_{k=1}^\infty$, respectively, so that $\ndim{\phi^0_k}\leq m$ and $\ndim{\phi^1_k}\leq n$ for all $k \in \bN$. We see that $(\phi^0_k+\phi^1_k)_{k=1}^\infty$ is a lift of $\phi^0+\phi^1$. We will show that $\ndim{\phi^0_k+\phi^1_k}\leq \ndim{\phi^0_k}+\ndim{\phi^0_k} +1$. To do this, let $\cF\subset A$ be a finite subset and let $\epsilon>0$. Then we may construct diagrams of the form
		\begin{equation}\label{eq:lem:sum}		
			\begin{tikzcd}
				A \arrow[rr, "\phi^0_k"] \arrow[dr, "\theta_k^0"'] & & B & A \arrow[rr, "\phi^1_k"] \arrow[dr, "\theta_k^1"'] & & B \\
				& \bigoplus_{j=0}^m F^j \arrow[ur, "\psi_k^0"']& & & \bigoplus_{j=m+1}^{n+m+1} F^j \arrow[ur, "\psi_k^1"']&			
			\end{tikzcd}	
		\end{equation}
		so that $\theta_k^0$ and $\theta_k^1$ are cpc, $\restrict{\psi_k^0}{F^j}$ is cpc$_\perp$ for $j=0,\dots, m$, $\restrict{\psi_k^1}{F^i}$ is cpc$_\perp$ for $i=m+1,\dots, n+m+1$, and $\norm{\phi^l_k(a)-\psi_k^l(\theta_k^l(a))}<\epsilon/2$ for all $a \in \cF$ and $l=0,1$. We then see that
		\begin{equation}		
			\begin{tikzcd}
				A \arrow[rr, "\phi^0_k+\phi^1_k"] \arrow[dr, "\theta_k^0+\theta_k^1"'] & & B \\
				& \bigoplus_{j=0}^{n+m+1} F^j \arrow[ur, "\psi_k^0+\psi_k^1"']&  &
			\end{tikzcd}	
		\end{equation}
		is a system of nuclear dimension approximations of $\phi_k^0+\phi_k^1$ for $\cF$ up to a tolerance of $\epsilon$. 
        The sum $\theta_k^0+\theta_k^1$ is cp and $\|\theta_k^0+\theta_k^1\| = \max\{\|\theta_k^0\|,\|\theta_k^1\|\} \leq 1$ since the ranges of $\theta_k^0$ and $\theta_k^1$ are orthogonal.
        Furthermore, $\restrict{(\psi_k^0+\psi_k^1)}{F^j} = \restrict{\psi_k^l}{F^j}$, for the appropriate choice of $l$ based on $j$, which is cpc$_\perp$ for $0\leq j \leq n+m+1$. We thus see that $\ndim{\phi^0_k+\phi^1_k}\leq n+m+1$.
	\end{proof}
	
	Our next goal is to strengthen the bound in Lemma~\ref{lem:sum} in the case when $\phi^0$ and $\phi^1$ have orthogonal ranges.  As a warm up, suppose $A$ and $B$ are $C^*$-algebras with $A$ separable and $\phi^0, \phi^1 \colon A \rightarrow B$ are cpc$_\perp$ maps with finite nuclear dimension and orthogonal ranges.  If $e \in A$ is a strictly positive element, then for $i \in \{0, 1\}$ the range of $\phi^i$ is contained in the hereditary subalgebra $\overline{\phi^i(e) B \phi^i(e)} \subseteq B$.  Indeed, writing $\Phi^i$ for the $^*$-homomorphism induced by $\phi^i$ and $z \in C_0(0, 1]$ for the element $z(t) = t$ (so $\phi^i(a) = \Phi^i(z \otimes a)$ for all $b \in A$), we have
	\begin{equation}\label{eq:strictly-positive}
		\lim_{n \rightarrow \infty} \phi^i(e)^{1/n} \phi^i(a) = \lim_{n \rightarrow \infty} \Phi^i(z^{1 + 1/n} \otimes e^{1/n} a) = \Phi^i(z \otimes a) = \phi^i(a), \quad a \in A.
	\end{equation}
	Then arguing as in the proof of Lemma~\ref{lem:sum}, one can show
	\begin{equation}
		\ndim{\phi^0 + \phi^1} \leq \max \{\ndim{\phi^0}, \ndim{\phi^1} \}.
	\end{equation}
	We will need the analogous estimate for $\nndim{\,\cdot\,}$; see Proposition~\ref{prop:orthogonal-sum}.  It is unclear if it can be reduced to the case of $\ndim{\,\cdot\,}$, but a very similar proof works.  The following is modeled on the proof of \cite[Proposition~1.6]{BGSW22}, which, in turn, is based on \cite[Proposition~2.5]{WZ10}.
	
	\begin{lem}\label{lem:corestrict}
		Let $A$ and $B$ be $C^*$-algebras, and let $\phi \colon A \rightarrow B$ be a cp map with $\ndim{\phi} < \infty$. Given $\epsilon > 0$ and a finite set $\mathcal F \subseteq A$, there exists $\delta > 0$ and a positive contraction $e \in A$ such that the following holds:
        for every positive contraction $h \in B$ satisfying
		\begin{equation}
			 \|h\phi(e) - \phi(e) \| < \delta,
        \end{equation}
		there is a cp map $\hat\phi \colon A \rightarrow \overline{hBh}$ such that $\ndim{\hat\phi} \leq \ndim{\phi}$ and
		\begin{equation}
			\|\phi(a) - \hat\phi(a) \| < \epsilon, \quad a \in \mathcal F.
		\end{equation}
	\end{lem}
	
	\begin{proof}
		Set $n = \ndim{\phi} < \infty$.  
        Let $\epsilon > 0$ and  $\mathcal F \subseteq A$ be a finite set.
        Without loss of generality, we may assume that $0 < \epsilon < 1$ and $\mathcal F$ consists of positive contractions.
        
        Choose $\delta \in (0,1)$ sufficiently small such that
        \begin{align}
            3\delta^{1/2} &< \frac{1}{16}, \label{eq:our-delta-1}\\
            13(n+1)\delta^{1/8} &< \epsilon. \label{eq:our-delta-2}
        \end{align}
        Take $e \in A$ to be a positive contraction such that $\|e a - a\| < \delta$ for all $a \in \mathcal F$. Note that such an element must exist by considering an approximate unit for $A$.
        
        Suppose now that $h \in B$ is a positive contraction such that 
        \begin{equation}\label{eq:approx-unit-delta}
			 \|h\phi(e) - \phi(e) \| < \delta,
        \end{equation}
        Since $\ndim{\phi} = n$, there exist finite dimensional $C^*$-algebras $F^0, \ldots, F^n$ and cp maps 
		\begin{equation}
			A \overset\theta\longrightarrow F^0 \oplus \cdots \oplus F^n \overset\psi\longrightarrow B
		\end{equation}
         such that $\theta$ is cpc, for each $i \in \{0\ldots,n\}$ the component map $\psi^i = \psi|_{F^i}$ is cpc$_\perp$, and
        \begin{equation}\label{eq:nd-estimate}
			 \|\phi(a) - \psi(\theta(a)) \| < \delta,
        \end{equation}         
        Writing $\phi_0 = \psi \circ \theta$, it follows from \eqref{eq:approx-unit-delta} and \eqref{eq:nd-estimate} that 
		\begin{equation}\label{eq:approx-unit-3delta}
			\| h \phi_0(e) - \phi_0(e) \| < 3\delta.
		\end{equation} 
		
        We now observe that $\theta(e)$ is a positive contraction in a finite-dimensional C$^*$-algebra. Hence, the spectrum of $\theta(e)$ is finite. Set $p = \chi_{[\delta^{1/2}, 1]}(\theta(e)) \in F$, i.e.\ $p$ is the spectral projection of $\theta(e)$ corresponding to the interval $[\delta^{1/2}, 1]$. 
		By construction, $p \leq \delta^{-1/2} \theta(e)$, and so
		\begin{equation}\label{eq:spectral-proj-estimate}
			\psi(p) \leq \delta^{-1/2} \phi_0(e). 
		\end{equation}
        Write $p = p_0 \oplus \cdots \oplus p_n$ where $p_i \in F^i$ for $i=\{0,\ldots,n\}$. Then it follows from \eqref{eq:approx-unit-3delta} and \eqref{eq:spectral-proj-estimate} that
		\begin{equation}\label{eq:oz-perturbation-estimate}
			\begin{split}
				\|\psi^i(p_i) (1 - h)\|^2 
				&= \|(1 - h) \psi^i(p_i)^2 (1 - h) \|  \\
				&\leq \|(1 - h) \psi^i(p_i) (1 - h)\| \\
				&\leq \|(1 - h) \psi(p)(1 - h)\| \\
				&\leq \delta^{-1/2} \|(1 - h) \phi_0(e)\| \\
				&< 3\delta^{1/2}.
			\end{split}
		\end{equation} 
        
        Let $\hat F^i = p_iF^i p_i$.
        We now apply \cite[Lemma~3.6]{KW04} using the estimates \eqref{eq:oz-perturbation-estimate} and \eqref{eq:our-delta-1}.
        We obtain a cpc$_\perp$ map $\hat\psi^i \colon \hat F^i \rightarrow \overline{hBh}$ such that
		\begin{equation}\label{eq:KW-estimate}
			\| \hat\psi^i(x) - \psi^i(x) \| \leq 8 \delta^{1/8} \|x\|.
		\end{equation}
        for all $x \in \hat F^i$.
        Let $\hat F = \hat F^0 \oplus \dots \oplus \hat F^n = pFp$.
		Define $\hat\theta = p \theta(\,\cdot\,)p$,  $\hat\psi = \hat\psi^0 \oplus \cdots \oplus \hat\psi^n$, and $\hat\phi = \hat \psi \circ \hat \theta$. 
         
        Let $a \in \mathcal F$. By the triangle inequality, we have
		\begin{equation}\label{eq:sam-estimate-step1}
				\|\phi_0(a) - \hat\phi(a) \| 
                \leq \|\psi(\theta(a)) - \psi(p\theta(a)p) \| + \|\psi(p\theta(a)p) - \hat\psi(p\theta(a)p) \|. 
		\end{equation}
        We estimate the first term on the right-hand side of \eqref{eq:sam-estimate-step1}. Using the definition of $e \in A$ and the fact that $\psi$ is the sum of $n+1$ cpc maps, we have
		\begin{equation}\label{eq:sam-estimate-step2}
			\begin{split}
				\| \psi(\theta(a)) - \psi(p \theta(a) p) \| 
				&\leq 2\|\psi\| \|(1 - p) \theta(a) \| \\
				&\leq 2\|\psi\| \|(1 - p) \theta(ea)\| + 2\|\psi\|\delta\\
                &\leq 2(n+1)\|(1 - p) \theta(ea)\| +2(n+1)\delta.
			\end{split}
		\end{equation}
		Using the $C^*$-identity and the Schwarz inequality, we obtain
		\begin{equation}\label{eq:sam-estimate-step3}
			\begin{split}
				\|(1 - p) \theta(ea) \|^2 &= \|(1 - p) \theta(ea) \theta(ea)^* (1 - p) \| \\
				&\leq \| (1 - p) \theta(ea^2 e) (1 - p) \| \\
				&\leq \|(1 - p) \theta(e)(1 - p)\| \\
				&\leq \delta^{1/2},
			\end{split}
		\end{equation}
        where the last inequality follows from the fact that $p$ was chosen to be the spectral projection of $\theta(e)$ corresponding to the interval $[\delta^{1/2}, 1]$.
        We now estimate the second term on the right-hand side of \eqref{eq:sam-estimate-step1}. Using \eqref{eq:KW-estimate} and the fact that $\theta$ and each $\psi^i$ are cpc, we obtain
		\begin{equation}\label{eq:sam-estimate-step4}
			\begin{split}
				 \|\psi(p\theta(a)p) - \hat\psi(p\theta(a)p) \| 
				&=  \sum_{i=0}^n \|\psi^i(p_i \theta(a) p_i) - \hat\psi^i(p_i \theta(a)p_i) \| \\
				&\leq 8(n+1)\delta^{1/8}.
			\end{split}
		\end{equation}
		Substituting \eqref{eq:sam-estimate-step2}, \eqref{eq:sam-estimate-step3}, and \eqref{eq:sam-estimate-step4} back into \eqref{eq:sam-estimate-step1}, we obtain
		\begin{equation}\label{eq:sam-estimate-step5}
			\begin{split}
				\|\phi_0(a) - \hat\phi(a) \| &\leq 2 (n + 1) \delta^{1/4} + 2(n+1)\delta + 8(n+1)\delta^{1/8} \\
			\end{split}
		\end{equation}
        Hence, using \eqref{eq:our-delta-2}, \eqref{eq:nd-estimate}, \eqref{eq:sam-estimate-step5} and that $\delta \in (0,1)$, we have 
        \begin{equation}
			\begin{split}
				\|\phi(a) - \hat\phi(a)\| &< 12(n+1)\delta^{1/8} +\delta\\ 
                &< 13(n+1)\delta^{1/8} \\
                &< \epsilon \qedhere
			\end{split}
		\end{equation}
	\end{proof}
	
	\begin{prop}\label{prop:orthogonal-sum}
		Suppose $A$ and $B$ are $C^*$-algebras such that $A$ is separable and $\phi^0, \phi^1 \colon A \rightarrow B_\omega$ are cpc$_\perp$ maps with orthogonal ranges.  Then 
		\begin{equation}
			\nndim{\phi^0 + \phi^1} \leq \max\{ \nndim{\phi^0}, \nndim{\phi^1} \}.
		\end{equation}
	\end{prop}
	
	\begin{proof}
		We may assume $\nndim{\phi^0}$ and $\nndim{\phi^1}$ are finite since the result is trivial otherwise.  Set $n=\max(\nndim{\phi^0},\nndim{\phi^1})$. For $i \in \{0, 1$\}, choose a sequence of cp maps $(\phi_m^i \colon A \rightarrow B)_{m=1}^\infty$ lifting $\phi^i$ such that each $\phi_m^i$ has nuclear dimension at most $n$.  
        
        Let $e \in A$ be a strictly positive element. Note that $\phi^0(e)$ and $\phi^1(e)$ are orthogonal positive contractions.  
        By \cite[Corollary~3.8]{LP98}, 
        we can find lifts $(b_{0,m})_{m=1}^\infty$ of $\phi^0(e)$ and  $(b_{0,1})_{m=1}^\infty$ of $\phi^1(e)$ consisting of positive contractions such that $b_{0, m} b_{1, m} = 0$ for all $m \geq 1$. 	
		
		Let $(\mathcal F_\ell)_{\ell = 1}^\infty$ be an increasing sequence of finite subsets of $A$ with dense union and let $(\epsilon_\ell)_{\ell = 1}^\infty$ be a decreasing sequence of strictly positive real numbers that converges to 0.  For each $\ell \in \bN$, let $\delta_\ell > 0$ and $e_\ell \in A$ be the output of Lemma~\ref{lem:corestrict} for $\cF=\cF_\ell$ and $\epsilon=\epsilon_\ell$. We may further arrange for $\epsilon_\ell\geq \delta_\ell > \delta_{\ell + 1}$ for all $\ell \geq 1$. Let $S_\ell \subseteq \mathbb N$ denote the set of $m \in \mathbb N$ such that $m \geq \ell$ and 
		\begin{equation}\label{eq:S-ell}
			 \inf_{k \in \mathbb N} \max_{i \in \{0, 1\}}\| b_{m, i}^{1/k} \phi_m^i(e_\ell) - \phi_m^i(e_\ell) \| < \delta_\ell.
		\end{equation}
		Then $S_{\ell+1} \subseteq S_{\ell}$ for all $\ell\geq 1$ and $\bigcap_{\ell = 1}^\infty S_\ell = \emptyset$
        
		We claim $S_\ell \in \omega$ for all $\ell \geq 1$.  To this end, first note that, as in \eqref{eq:strictly-positive},
		\begin{equation}
			\lim_{k \rightarrow \infty} \phi^i(e)^{1/k} \phi^i(a) = \phi^i(a)
		\end{equation}
        for all $a \in A$ and $i \in \{0,1\}$.
		In particular, for each $\ell \in \bN$, there is $k_\ell \in \mathbb N$ such that for $i \in \{0,1\}$
		\begin{equation}
			\| \phi^i(e)^{1/k_\ell} \phi^i(e_\ell) - \phi^i(e_\ell) \| < \delta_\ell.
		\end{equation}
		Hence, for $i \in \{0,1\}$, we have
		\begin{equation}
			\lim_{m \rightarrow \omega} \| b_{i, m}^{1/k_\ell} \phi_m^i(e_\ell) - \phi_m^i(e_\ell) \|  < \delta_\ell,
		\end{equation}
		so \eqref{eq:S-ell} holds for $\omega$-many $m \in \mathbb N$.  Hence $S_\ell \in \omega$.
		
        For $\ell \in \mathbb N$,  $m \in S_{\ell + 1} \setminus S_\ell$, and $i \in \{ 0, 1\}$, by  Lemma~\ref{lem:corestrict} with $h = b_{i, m}^{1/k_\ell}$, there is a cp map $\psi_m^i \colon A \rightarrow \overline{b_{i, m}Bb_{i, m}}$ such that $\ndim{\psi^i_m} \leq \ndim{\phi^i_m} \leq n$ and
		\begin{equation}
			\| \phi_m^i(a) - \psi_m^i(a) \| < \epsilon_\ell, \quad a \in \mathcal F_\ell,
		\end{equation}
        where we have used the fact that $b_{i, m}^{1/k_\ell}$ and $b_{i, m}$ generate the same herditary subalgebra.
        For $m \in \mathbb N \setminus S_1$ and $i \in \{ 0, 1\}$, set $\psi_m^i = 0$. 
		Since $(\mathcal F_\ell)_{\ell=1}^\infty$ is increasing and $(\mathcal \epsilon_\ell)_{\ell=1}^\infty$ is decreasing, we also note that this inequality holds for all $m \in S_\ell$, and hence for $\omega$-many $m$.  Therefore, $(\psi_m^i)_{m=1}^\infty$ is a lift of $\phi^i$.  As $\psi_m^0$ and $\psi_m^1$ have ranges in orthogonal hereditary subalgebras, $\psi_m^0 + \psi_m^1$ has nuclear dimension at most $n$.  Hence $\nndim{\phi^0  + \phi^1} \leq n$.  
	\end{proof}
	
	We record the following fact relating $\nndim{\cdot}$ to $\ndim{\cdot}$, which is how the nuclear dimension bounds in Theorems~\ref{intro-thm:graphs-with-dim-one}, \ref{intro-thm:dim-at-most-two}, and \ref{intro-thm:extension-result} will be computed. 
	
	\begin{prop}\label{prop:nndimisndim}
		Let $A$ be a $C^*$-algebra and let $\iota_A\colon A \to A_\omega$ denote the diagonal inclusion. Then $\ndim{A}=\nndim{\iota_A}$.
	\end{prop}
	
	\begin{proof}
		Let $n=\nndim{\iota_A}<\infty$. Since $\iota_A$ lifts to the sequence $(\mathrm{id}_A)_{k=1}^\infty$, we clearly have $\ndim{\mathrm{id}_A}\geq\nndim{\iota_A}$. We show the reverse inequality.
		
		Suppose that $\iota_A$ lifts to a sequence of cp maps $(\theta_k\colon A \to A)_{k=1}^\infty$ with $\ndim{\theta_k}\leq n$ for all $k$. Let $\cF \subset A$ be a finite set, and let $\epsilon>0$. Since $\lim_{k \to \omega}\norm{\theta_k(a)-a}=0$ for all $a \in A$, we may find a $k \in \bN$ such that $\norm{\theta_k(a)-a}<\epsilon/2$ for all $a \in \cF$. Since $\ndim{\theta_k}\leq n$, we may find finite-dimensional $C^*$-algebras $F^0,\dots, F^n$ a cpc map $\phi\colon A \to \bigoplus_{i=0}^n F^i$, and a cp map $\psi\colon \bigoplus_{i=0}^n F^i \to A$ such that $\norm{\psi(\phi(a))-\theta_k(a)}<\epsilon/2$ and $\restrict{\psi}{F^i}$ is cpc$_\perp$ for $i=0,\dots, n$. We then see that $\norm{\psi(\phi(a))-a}<\epsilon$, so $\ndim{\mathrm{id}_A}\leq n$.
	\end{proof}
	
	Our arguments focus on finding approximate unitary equivalences between maps into $A_\omega$. We give a straightforward proof of the following for clarity.
    
	\begin{prop}\label{prop:unitary-invariant}
		Let $A$ and $B$ be C$^*$-algebras with $A$ separable. Suppose that $\phi, \psi \colon A \to B_\omega$ are cpc$_\perp$ maps that are approximately unitarily equivalent via unitaries in the minimal unitization of ${B_\omega}$. Then $\nndim{\phi}= \nndim{\psi}$.	
	\end{prop}
	
	\begin{proof}
		By Kirchberg's $\epsilon$-test, $\phi$, $\psi$ are unitarily equivalent. Say $\psi = \mathrm{ad}(v) \circ \phi$ where $v$ is a unitary in the minimal unitization of ${B_\omega}$. Lift $v$ to a sequence $(v_n)_{n=1}^\infty$ of unitaries in the minimal unitization of ${B}$ and $\phi$ to a sequence of cp maps $(\phi_n)_{n=1}^\infty$ with $\ndim{\phi_n}\leq \nndim{\phi}$. Then the sequence $(\mathrm{ad}(v_n)\circ \phi_n)_{n=1}^\infty$ is a lift of $\psi$, and nuclear dimension does not increase when post-composing by a $^*$-homomorphism. Symmetrizing the argument shows that $\nndim{\phi}=\nndim{\psi}$.
	\end{proof}
	
	The following is the variation of nuclear dimension for $C^*$-algebras that will be used in the statement and proof of Theorem~\ref{intro-thm:extension-result}.
	
	\begin{defi}\label{def:fndim}
		Let $A$ be a separable $C^*$-algebra. We say that $A$ has \emph{faithful nuclear dimension} at most $n$, denoted $\fndim{A} \leq n$, if there exist cpc$_\perp$ maps $\phi^0, \phi^1 \colon A \to A_\omega$ such that
		\begin{enumerate}
			\item \label{it:def:fndim1}$\phi^0(a)+\phi^1(a)=a$ for all $a \in A$,
			\item \label{it:def:fndim2}$\phi^0$ is totally faithful,
			\item \label{it:def:fndim3}$\nndim{\phi^0} = 0$, and
			\item \label{it:def:fndim4}$\nndim{\phi^1} \leq n - 1$.
		\end{enumerate}
	\end{defi}
	
	We always have the following inequality, though we will see that it can be strict. Note that it does not rely upon faithfulness of $\phi^0$, as this simply follows from some earlier facts about $\nndim{\, \cdot \,}$.
	
	\begin{prop}\label{prop:ndim<fndim}
		If $A$ is a separable $C^*$-algebra, then $\ndim{A}\leq \fndim{A}$.
	\end{prop}

    \begin{proof}
        This follows from Lemma~\ref{lem:sum} and Proposition~\ref{prop:nndimisndim}
    \end{proof}
	
	For graph algebras, we will need to carefully compute the faithful nuclear dimensions of extensions by minimal ideals which arise as stabilizations of $C(\bT)$ or of Kirchberg algebras. For the $C(\bT)\otimes \cK$ ideals, it is more straightforward to compute faithful nuclear dimension approximations of $C(\bT)$ itself and then prove $\fndim{C(\bT)\otimes \cK}\leq \fndim{C(\bT)}$. Of course, there is nothing special about $C(\bT)$ or $\cK$ in this tensor product, so we prove the following, more general, statement.

	\begin{prop}\label{prop:fndimstable}
		Let $D$ be an AF algebra and $A$ a separable $C^*$-algebra. Then $\fndim{A \otimes D} \leq \fndim{A}$.
	\end{prop}

    \begin{proof}
		Let $n = \fndim{A}$. Without loss of generality, we may assume that $n<\infty$ and in particular that $A$ is nuclear. Let $\phi^0, \phi^1 \colon A \to A_\omega$ be such that $\phi^0+\phi^1=\iota_A$, $\phi^0$ is totally faithful, $\nndim{\phi^0}=0$, and $\nndim{\phi^1}\leq n-1$. 
        Lift $\phi^0$ and $\phi^1$ to sequences of cp maps $(\phi^i_k\colon A \to A)_{k=1}^\infty$ for $i=0,1$ with $\ndim{\phi^0_k}=0$ and $\ndim{\phi^1_k}\leq n-1$. 
        
        Note for all $k \in \bN$, $\ndim{\phi^0_k\otimes \mathrm{id}_D}\leq \ndim{\phi^0_k}=0$ and $\ndim{\phi^1_k \otimes \mathrm{id}_D}\leq \ndim{\phi^1_k}\leq n-1$ by Lemma~\ref{lem:nndim-stable}. Let $\psi^i\colon A \otimes D \to (A \otimes D)_\omega$ denote the map induced by the sequences $(\phi^i_k \otimes \text{id}_D)_{k=1}^\infty$ for $i=0,1$. Then $\nndim{\psi^0}=0$ and $\nndim{\psi^1}\leq n-1$. 
        Since  $\phi^0+\phi^1=\iota_A$, we have $\lim_{k\to\omega} (\phi^0_k(a)\otimes d+\phi^1_k(a)\otimes d) = a \otimes d$ for all $a \in A$ and $d \in D$, so $\psi^0+\psi^1=\iota_{A \otimes D}$.
        It remains to show that $\psi^0$ is totally faithful.

        To see this, let $\Phi^0\colon C_0(0,1] \otimes A \to A_\omega$ be the faithful $^*$-homomorphism dilating $\phi^0$. Since $A$ is nuclear, we may use the Choi--Effros lifting theorem (see \cite{CE76a}) to obtain a sequence of cpc maps $(\Phi_k^0 \colon C_0(0, 1] \otimes A \rightarrow A)_{k=1}^\infty$ lifting $\Phi^0$. We have
		\begin{align}
			\label{it:fndimaftensor1} \lim_{k \to \omega}\|\Phi_k^0(z \otimes a) - \phi_k^0(a)\|&=0, \quad a \in A, \\
			\label{it:fndimaftensor2} \lim_{k \to \omega}\|\Phi^0_k(\xi_1\xi_2) - \Phi^0_k(\xi_1)\Phi^0_k(\xi_2)\|&=0, \quad \xi_1, \xi_2 \in C_0(0,1]\otimes A.
		\end{align}
		Let $\Psi^0_k\colon C_0(0,1]\otimes A \otimes D \to A \otimes D$ be given by $\Phi^0_k \otimes \text{id}_D$. Then by \eqref{it:fndimaftensor1} and \eqref{it:fndimaftensor2} the map $\Psi^0\colon C_0(0,1]\otimes A \otimes D \to (A\otimes D)_\omega$ induced by the sequence $(\Psi^0_k)_{k=1}^\infty$ is the $^*$-homomorphism dilating $\Psi^0$ to the cone. 
        Supose $\Psi^0$ is not faithful. Then,
        by Kirchberg's slice lemma (\cite[Lemma~4.1.9]{RordamClassificationBook}), there must be a non-zero simple tensor $\xi\otimes d \in \ker \Psi^0$ for some $\xi \in C_0(0,1]\otimes A$ and $d \in D$. However, if $d\neq 0$, this implies that $\lim_{k \to \omega}\|\Phi^0_k(\xi)\|=0$. By the faithfulness of $\Phi^0$, this implies that $\xi=0$, and so we conclude that $\psi^0$ is totally faithful.
    \end{proof}
	
	The rest of this section is devoted to computing the value of $\fndim{A}$ for a series of examples.  For the applications to graphs, the most important examples are $C^*$-algebras stably isomorphic to $ C(\bT)$ (Theorem~\ref{thm:nndim:abelian}) and $\cO_\infty$-stable $C^*$-algebras (Example~\ref{ex:nndim:purelyinfinite}).  In both cases, the faithful nuclear dimension is 1, as it is the usual nuclear dimension approximations for these $C^*$-algebras that realize the upper bound on the faithful nuclear dimension. 
	
	\begin{ex}\label{ex:nndim:0}
		No non-zero $C^*$-algebra $A$ has faithful nuclear dimension zero, where we consider $\phi^1 = 0$ in this case. To see this, note that Definition~\ref{def:fndim}\ref{it:def:fndim1} forces $\phi^0 = \iota_A$.  Since $\phi^0$ is a $^*$-homomorphism, the corresponding $^*$-homomorphism $\Phi^0 \colon C_0(0,1]\otimes A \to A_\omega$ is given by $\iota_A \circ \mathrm{ev}_1$, which is clearly not faithful.  Thus $\fndim{A}\geq1$ for all non-zero $C^*$-algebras $A$.
	\end{ex}
	
	\begin{ex}\label{ex:nndim:finitedimensional}
		Suppose $F$ is a non-zero finite-dimensional $C^*$-algebra.  Then $\fndim{F}=\infty$. To see this, note that the unit ball of a finite-dimensional algebra is compact, so $F_\omega \cong F$. An easy dimension counting argument shows that there can be no faithful map from $C_0(0,1]\otimes F \to F_\omega \cong F$.
	\end{ex}
	
	\begin{ex}\label{ex:nndim:simplezstab}
		If $A$ is a unital simple separable nuclear $\cZ$-stable $C^*$-algebra, then $\fndim{A}=1$.  Indeed, the first factor inclusion $\mathrm{id}_A \otimes 1_\cZ\colon A \into A \otimes \cZ$ is approximately unitarily equivalent to an isomorphism $\psi \colon A \to A \otimes \cZ$ via some sequence of unitaries $(u_n)_{n=1}^\infty$ in $A \otimes \cZ$.  Let $u \in (A \otimes \cZ)_\omega$ be the unitary defined by the sequence $(u_n)_{n=1}^\infty$.  Then 
		\begin{equation}
			(\psi^{-1})_\omega\circ \mathrm{ad}(u)\circ (\mathrm{id}_A \otimes 1_\cZ) _\omega = \iota_A \colon A \into A_\omega.
		\end{equation}
		Let $h_0$ be a positive contraction of full spectrum in $\cZ$, and let $h_1=1-h_0$. As part of the proof of Theorem~H in \cite{CETWW21}, it is proven that the maps $\mathrm{id}_A \otimes h_0$ and $\mathrm{id}_A \otimes h_1$ have nuclear dimension zero. 
        Write $\Phi^i:C_0(0,1] \otimes A \rightarrow A \otimes \cZ$ for the $^*$-homomorphism induced by the cpc$_\perp$ map $\mathrm{id}_A \otimes h_i$ where $i \in \{0,1\}$. 
        Then, for $\xi \in C_0(0,1]$ and $a \in A$, we have $\Phi^i(\xi \otimes a) = a \otimes \xi(h_i)=0$ implies $a=0$ or $\xi=0$ (since $\sigma(h_i)=[0,1]$). By Kirchberg's slice lemma (\cite[Lemma~4.1.9]{RordamClassificationBook}), it follows that $\Phi^i$ is faithful, so $\mathrm{id}_A \otimes h_i$ is totally faithful. The compositions 
		\begin{equation}
			\phi^i\coloneq (\psi^{-1})_\omega\circ \mathrm{ad}(u)\circ (\mathrm{id}_A \otimes h_i) _\omega, \quad i = 0, 1,
		\end{equation}
		show $\fndim{A} \leq 1$ by Proposition~\ref{prop:nndim<ndim}, and Example~\ref{ex:nndim:0} gives the reverse inequality.
	\end{ex}
	
	\begin{ex}\label{ex:nndim:zstabaf}
		If $A$ is a $\mathcal Z$-stable AF algebra, then $\fndim{A} = 1$.  Indeed, by \cite[Theorem~2.2]{TW07}, there is a unital embedding $\mathcal Z \subseteq (A_\omega \cap A')/ (A_\omega \cap A^\perp)$.  Fix a positive contraction $\bar h \in \mathcal Z$ with $\sigma(\bar h) = [0, 1]$ and lift $\bar h$ to a positive contraction $h \in A_\omega \cap A'$.  Define $\phi^0, \phi^1 \colon  A \rightarrow A_\omega$ by 
		\begin{equation}
			\phi^0(a) = ha \quad \text{and} \quad \phi^1(a) = (1-h)a, \qquad a \in A
		\end{equation}
		and note that $\phi^0$ and $\phi^1$ are cpc$_\perp$ maps.  Then $\nndim{\phi^0}=\nndim{\phi^1}=0$ since $\ndim{A} = 0$, and $\phi^0(a) + \phi^1(a) = a$ for all $a \in A$.  So it suffices to show $\phi^0$ is totally faithful.  Let 
		\begin{equation}
			\Phi^0 \colon C_0(0, 1] \otimes A \rightarrow A_\omega
		\end{equation}
		be the $^*$-homomorphism induced by $\phi^0$.  For non-zero $\xi \in C_0(0,1]$ and $a \in A$, note that $\xi(\bar h) \in (A_\omega \cap A') / (A_\omega \cap A^\perp)$ is full since $\mathcal Z$ is simple.  So $a$ is in the ideal generated by $\xi(h) a + A_\omega \cap A^\perp$.  As $a \notin A_\omega \cap A^\perp$ (since $a^*a \neq 0)$, we have $\Phi(\xi \otimes a) = \xi(h)a \neq 0$.  By Kirchberg's slice lemma (\cite[Lemma~4.1.9]{RordamClassificationBook}), $\Phi^0$ is faithful.
	\end{ex}
	
	\begin{ex}\label{ex:nndim:purelyinfinite}
		All nuclear $\cO_\infty$-stable $C^*$-algebras $A$ have $\fndim{A} = 1$, following from the work in \cite[Theorem~A]{BGSW22}, which states $\ndim{A}  = 1$. The argument is essentially the same as Example~\ref{ex:nndim:zstabaf}, except that the domain is no longer AF.  Thus, $\phi^0$ and $\phi^1$ have nuclear dimension zero by \cite[Theorem~D]{BGSW22}, recalling that $C_0(0, 1] \otimes A$ is quasidiagonal by \cite{V91}.
	\end{ex}
	
	In order to compute the nuclear dimension of graph $C^*$-algebras, we will need the computation $\fndim{C(\mathbb T)} = 1$ (Theorem~\ref{thm:nndim:abelian}).  We give a more general computation of $\fndim{C(X)}$ for compact metrizable spaces $X$ with a suitable connectivity hypothesis: namely, connected spaces with no totally disconnected open sets. 	
	The rough idea is to take a sequence of $\dim(X)+1$-colorable finite covers of $X$ with the diameter of the sets in each cover converging to $0$. Then, letting $h_n$ be the sum of one color of a partition of unity subordinate to the $n$th cover of $X$, set $h=(h_n)\in C(X)_\omega$. The maps $\phi^0(a) = ha$ and $\phi^1(a) =(1-h)a$ will be ultraproducts of maps with nuclear dimension $0$ and $\dim(X)-1$ respectively, and they clearly sum to the canonical inclusion $C(X) \into C(X)_\omega$. However, in general, $\phi^0$ will fail to be totally faithful. Extra care must be taken in the choice of color used to form $h_n$; specifically, we will arrange $h_n$ to have range $[0,1]$ on each element of the cover on which it is supported. 
	
	The following lemma constructs partitions of unity of the appropriate form. 
	
	\begin{lem}\label{lem:pointsettopology1}
		Let $X$ be a connected compact Hausdorff space such that no non-empty open subset is totally disconnected.  If $\cU = \bigsqcup_{i=0}^{\dim(X)}\cU^i$ is a $(\dim X + 1)$-colored finite open cover of $X$ by non-empty proper open sets, then there is a partition of unity $(\xi_U)_{U \in \cU} \subseteq C(X)$ subordinate to $\cU$ with $\xi_U(U) = [0, 1]$ for every $U \in \cU^0$.
	\end{lem}
	
	\begin{proof}
		Let $(\xi_U')_{U \in \cU}$ be a partition of unity subordinate to $\cU$, and fix $U \in \cU^0$.  Since $X$ is connected, $\mathrm{supp}(\xi_U')$ is a proper subset of $U$.  Therefore, $U \setminus \mathrm{supp}(\xi_U')$ is a non-empty open subset of $X$ and hence is not totally disconnected.  So there are distinct points $x_U, y_U \in U \setminus \mathrm{supp}(\xi_U)$ which belong to the same connected component of $U$.  By Tietze's extension theorem, there is $\eta_U \in C(X)_+$ with $\eta_U(x_U) = 1$, $\eta_U(y_U) = 0$, and $\mathrm{supp}(\eta_U) \subseteq U$.  Using Tietze's extension theorem again, there is $\zeta \in C(X)_+$ such that $\zeta(x_U) = 0$ for all $U \in \cU^0$ and $\zeta(x) = 1$ for all $x \in \bigcup_{U \in \cU^0} g_U^{-1}(0)$.
		
		Let $\xi \in C(X)$ be given by
		\begin{equation}
			\xi = \sum_{U \in \cU^0} (\xi_U' + \eta_U) + \sum_{U \in \cU\setminus \cU^0} \xi_U' \zeta.
		\end{equation}
		Note that if $\eta_U(x) > 0$ for some $U \in \cU^0$, then $\xi(x) > 0$, and if $\eta_U(x) = 0$ for all $U \in \cU^0$, then $\zeta(x) = 1$. Thus $\xi_U'(x)\zeta(x)>0$ for some $U \in \cU$.
		Hence $\xi(x) > 0$ for all $x \in X$.  Define
		\begin{alignat}{2}
			\xi_U &= (\xi_U' + \eta_U)/\xi, &U \in \cU^0 \\
			\intertext{and}
			\xi_U &= \xi_U' \zeta / \xi,  &U \in \cU \setminus \cU^0.
		\end{alignat}
		Then $(\xi_U)_{U \in \cU}$ is a partition of unity subordinate to $\cU$.
		
		Let $U \in \cU^0$.  Note that $\xi_U'(x_U) = 0$ and $\eta_U(x_U) = 1$.  If $V \in \cU^0 \setminus \{U\}$, then $\xi'_V(x_U) = \eta_V(x_U) = 0$ since $x_U \in U$, $U \cap V = \emptyset$ and $\xi'_V$ and $\eta_V$ are supported on $V$.  Furthermore, $\zeta(x_U) = 0$,  so $\xi_U(x_U) = 1$ since $\xi(x_U)=1$.  Also note that $\xi_U(y_U) = 0$ since $\xi_U(y_U) = \eta_U(y_U) = 0$.  Since $x_U$ and $y_U$ belong to the same connected component of $U$, we have $\xi_U(U) = [0, 1]$.
	\end{proof}
	
	\begin{thm}\label{thm:nndim:abelian}
		Let $X$ be a connected compact metric space such that no non-empty open set is totally disconnected. Then $\fndim{C(X)} = \dim(X)$.
	\end{thm}
	
	\begin{proof}
		First note that 
		\begin{equation}
			\dim(X) \leq \ndim{C(X)} \leq \fndim{C(X)}
		\end{equation}
		by Proposition~\ref{prop:ndim<fndim} and \cite[Proposition~2.4]{WZ10}, so it suffices to show $\fndim{C(X)} \leq \dim(X)$.
		
		For each $n \in \mathbb N$, let $\cU_n = \bigsqcup_{i=0}^{\dim(X)} \cU_n^i$ be a $(\dim(X)+1)$-colored finite open cover of $X$ consisting of non-empty proper open sets of diameter at most $1/n$.  We may further arrange that for each $n \in \mathbb N$ and $U \in \cU \setminus \cU^0_n$, there is $V \in \cU^0_n$ such that $U \cap V \neq \emptyset$.  Indeed, for fixed $n \in \mathbb N$, if there exists $1 \leq i \leq \dim(X)$ and $U \in \mathcal \cU^i$ with $V \cap U = \emptyset$ for all $V \in \cU^0$, we may replace $\cU_n^i$ with $\cU_n^i \setminus \{U\}$ and $\cU_n^0$ with $\cU_n^0 \cup \{U\}$.  As $\cU_n$ is finite, iterating this argument yields the desired partition $\cU_n = \bigsqcup_{i=0}^{\dim(X)} \cU_n^i$.
		
		Let $(\xi_{n, U})_{n=1}^\infty$ be a partition of unity of $X$ subordinate to $\cU^n$ as in Lemma~\ref{lem:pointsettopology1}.  For each $i$, define
		\begin{equation}
			h_n^i = \sum_{U \in \cU_n^i} \xi_{n,U} \in C(X)
		\end{equation}
		and note the $h_n^i$ are positive contractions summing to $1_{C(X)}$ such that for each $U \in \cU_n^0$, $h_n^0(U) = [0, 1]$.  Let $h_\omega^i \in C(X)_\omega$ denote the element represented by the sequence $(h_n^i)_{n=1}^\infty$.
		
		It suffices to show the maps 
		\begin{equation}
			\phi_\omega^i \colon C(X) \rightarrow C(X)_\omega \colon \eta \mapsto h_\omega^i \eta
		\end{equation}
		satisfy $\nndim{\phi_\omega^i} = 0$ and $\phi_\omega^0$ is totally faithful.  Indeed, in this case, the maps $\phi_\omega^0$ and $\sum_{i=1}^{\dim(X)} \phi_\omega^i$ realize the inequality $\fndim{C(X)} \leq \dim(X)$ by Lemma~\ref{lem:sum}.
		
		For each $n \in \mathbb N$, and $U \in \cU_n$, fix a point $x_{n, U}$ and, for each $0 \leq i \leq \dim(X)$, define
		\begin{equation}
			\phi_n^i \colon C(X) \rightarrow C(X) \colon \eta \mapsto \sum_{U \in \cU_n^i} \eta(x_{n, U}) h_n^i
		\end{equation}
		Then $\phi_n^i$ has nuclear dimension 0.  Fix $0 \leq i \leq \dim(X)$, $\eta \in C(X)$, and $\epsilon > 0$.  Let $N \in \mathbb N$ be such that for all $x, y \in X$ with $d(x, y) < 1/N$, we have $|\eta(x) - \eta(y)| < \epsilon$.  Let $n \in \mathbb N$ with $n \geq N$ and $x \in X$.  If $x \in U$ for some $U \in \cU_n^i$, then $d(x, x_{n, U}) < 1/N$, and 
		\begin{equation}
			|\phi_n^i(\eta)(x) - h_n^i(x) \eta(x)| \leq |\eta(x_{n, U}) - \eta(x)| < \epsilon, 
		\end{equation}
		and otherwise, $\phi_n^i(\eta)(x) = h_n^i(x) \eta(x_{n, U}) = 0$.  Hence
		\begin{equation}
			\|\phi_n^i(\eta) - h_n^i \eta \| < \epsilon, \quad n \geq N.
		\end{equation}
		It follows that the sequence $(\phi_n^i)_{n=1}^\infty$ lifts $\phi_\omega^i$, so $\nndim{\phi_\omega^i} = 0$.
		
		Now we'll show $\phi^0_\omega$ is totally faithful.  
        Let $\Phi^0 \colon C_0(0, 1] \otimes C(X) \rightarrow C(X)_\omega$ denote the induced $^*$-homomorphism and note that $\Phi^0(\zeta \otimes \eta) = \zeta(h_\omega^i) \eta$ for all $\zeta \in C_0(0, 1]$ and $\eta \in C(X)$ by Remark~\ref{remark:multiplicationoncone}.  
        By Kirchberg's slice lemma (\cite[Lemma~4.1.9]{RordamClassificationBook}), it suffices to show $\Phi^0(\zeta \otimes \eta) = 0$ implies $\zeta \otimes \eta = 0$.  To this end, suppose $\Phi^0(\zeta \otimes \eta) = 0$ and $\eta \neq 0$ 
        Fix $\epsilon > 0$.  
        Let $n \in \mathbb N$ be such that $\|\zeta(h_n^0) \eta\| < \epsilon \|\eta\|/3$ and for all $x, y \in X$, $d(x, y) < 1/n$ implies $|\eta(x) - \eta(y)| < \|\eta\|/3$.  Let $x \in X$ be such that $\eta(x) = \|\eta\|$ and let $U \in \cU_n$ be such that $x \in U$.  
        Choose $V  \in \cU_{n}^0$ such that $V \cap U \neq \emptyset$ and fix $y \in U \cap V$.  
        Since $U$ and $V$ have diameter at most $1/n$, for all $z \in V$, we have $d(x,y) < 1/n$ and $d(y, z) < 1/n$, so by the triangle inequality and the choice of $n$, $\abs{\eta(z)} \geq \|\eta\|/3$.  
        Further, since $\|\zeta(h_n^0) \eta\| < \epsilon \|\eta\|/3$, we have
		\begin{equation}
			|\zeta(h_n^0(z))| = |\zeta(h_n^0)(z)| < \epsilon, \quad z \in V.
		\end{equation}
		By the choice of $h_n^0$, we have $h_n^0(V) = [0, 1]$, so $\|\zeta\| < \epsilon$.  
        Since $\epsilon > 0$ was arbitrary, $\zeta = 0$, so $\Phi^0$ is faithful.
	\end{proof}
	
	\section{Proof of Theorem~\ref{thm:D}}\label{sec:proof-of-thm-D}
	
	Theorem~\ref{thm:D} is an important tool needed in the nuclear dimension computations stated in the introduction.  The proof is the culmination of several decades of work on non-stable extension theory, going back to Brown, Douglas, and Fillmore's highly influential paper \cite{BDF77}, and with more recent developments largely built in connection with the classification and regularity theories for simple nuclear $C^*$-algebras.  We recall the concepts and results needed in this section.
	
	Consider an extension
	\begin{equation}\label{eq:ext-section3}
		\begin{tikzcd}
			\mathsf e \colon 0 \arrow{r} & I \arrow{r} & A \arrow{r}{q} & B \arrow{r} & 0
		\end{tikzcd}
	\end{equation}
	of separable $C^*$-algebras with $I$ stable.  Letting $\mathcal M(I)$ and $\mathcal Q(I) = \mathcal M(I) / I$ denote the multiplier and corona algebras of $I$, respectively, there is a morphism of extensions
	\begin{equation}
		\begin{tikzcd} 
			0 \arrow{r} & I \arrow[equals]{d} \arrow{r} & A \arrow{d}{\lambda} \arrow{r}{q} & B \arrow{d}{\bar \lambda} \arrow{r} & 0\phantom{.} \\
			0 \arrow{r} & I \arrow{r} & \mathcal M(I) \arrow{r}{\bar q} & \mathcal Q(I) \arrow{r} & 0,
		\end{tikzcd}
	\end{equation}
	where $\lambda$ is the canonical map determined by $\lambda(a)x=ax$ for all $a \in A$ and $x \in I$. The extension $\mathsf e$ is \emph{full} if its \emph{Busby invariant} $\bar\lambda$ is full in the sense that $\bar\lambda(b)$ generates $\mathcal Q(I)$ as an ideal for all non-zero $b \in B$.
	
	It is a standard fact, using that $I$ is stable, that if $\mathsf e$ is full, then for all $a \in A_+ \setminus I$, $x \in I_+$, and $\epsilon > 0$, there exist $n \in \mathbb N$ and $v_1,\ldots,v_n \in A$  such that
	\begin{equation}
		\Big\|\sum_{i=1}^n v_i^* a v_i - x \Big\| < \epsilon.
	\end{equation}
	To see this, one argues as follows: Without loss of generality $x \neq 0$. Since $\bar\lambda(q(a)) \in \mathcal Q(B)_+$ generates $\mathcal Q(B)$ an an ideal, there are $n \in \mathbb N$ and $w_1, \ldots, w_n \in \mathcal \mathcal \cM(I)$ such that
	\begin{equation}
		1 - \sum_{i=1}^n w_i^* \lambda(a)  w_i \in I. 
	\end{equation}
	Then since $I$ is stable, there is an isometry $s \in \mathcal M(I)$ with 
	\begin{equation}
		\Big\| s^* \Big(1 - \sum_{i=1}^n w_i^* \lambda(a) w_i \Big) s \Big\| < \frac{1}{\|x\|}\epsilon.
	\end{equation}
	Then the elements $v_i = w_i s x^{1/2} \in I \subseteq A$ are as required.
	
	Elliott and Kucerovsky identified a condition on an extension allowing one to take $n = 1$ in the previous paragraph, in analogy with the fact that if $D$ is a simple purely infinite $C^*$-algebra, $a \in D_+ \setminus \{0\}$, $x \in D_+$, and $\epsilon > 0$, then there exists $v \in A$ such that $\|v^* a v - x \| < \epsilon$.  Specifically, they say the extension $\mathsf e$ is \emph{purely large} if, for every $a \in A_+ \setminus \{0\}$, there is a stable $C^*$-algebra $C \subseteq \overline{aIa}$ which generates $I$ as an ideal.  The following result is implicitly contained in Elliott and Kucerovsky's work in \cite{EK01}.  Indeed, the forward direction is \cite[Lemma~7]{EK01}, and the reverse direction amounts to noting that this is the only use of pure largeness in their proof that pure largeness implies nuclear absorption.  A self-contained proof is given in \cite{CGSTW26}.
	
	\begin{thm}[{\cite{CGSTW26}}]\label{pure-large}
		Let $\mathsf e \colon 0 \to I \to A \to B \to 0$ be an extension of separable $C^*$-algebras with $I$ stable.  Then $\mathsf e$ is purely large if and only if for every $a \in A_+ \setminus I$, $x \in I_+$, and $\epsilon > 0$, there exists $v \in A$ such that $\|v^*a v - x\| < \epsilon$.
	\end{thm}
	
	Elliott and Kucerovsky showed purely large extensions satisfy a suitable version of Voiculescu's non-commutative Weyl--von Neumann theorem (\cite{V76}).  We will use pure largeness in the following form, which combines Elliott and Kucerovsky's absorption theorem with  Szab\'o's $KK$-uniqueness theorem in \cite{S26}, which in turn improves (and depends on) the Dadarlat--Eilers stable uniqueness theorem from \cite{Dadarlat-Eilers1}.  Since these result will be applied to $^*$-homomorphisms out of cones, the $KK$-theoretic obstructions vanish by the homotopy invariance of $KK$-theory.
	
	\begin{thm}\label{thm:unique-lift}
		Suppose $0 \to I \to A \xrightarrow q B \to 0$ is a purely large extension of separable $C^*$-algebras with $I$ stable and let $D$ be a separable nuclear $C^*$-algebras.  If $\phi, \tilde\phi \colon D \rightarrow A$ are cpc$_\perp$ maps with $q \circ \phi = q \circ \tilde\phi$ and this common map is totally faithful, then there is a continuous path of unitaries $(u_t)_{t \geq 0}$ in the minimal unitization of $I$ such that
		\begin{equation}\label{eq:unique-lift}
			\lim_{t \rightarrow \infty} \|u_t \phi(d) u_t^* - \tilde\phi(d)\| = 0, \quad d \in D.
		\end{equation}
	\end{thm}
	
	\begin{proof}
		Set $\psi = q\circ \phi$ and let $\Psi \colon C_0(0, 1] \otimes D \rightarrow B$ be the induced $^*$-homomorphism. Consider the diagram 
		\begin{equation}
			\begin{tikzcd}
				\tilde{\mathsf e} \colon 0 \arrow{r} & I \arrow[equals]{d} \arrow{r} & E \arrow{d}{\tilde\Psi}\arrow{r}{\tilde q} \arrow{r} & C_0(0, 1] \otimes D \arrow{d}{\Psi} \arrow{r} & 0\phantom{,} \\
				\mathsf e \colon 0 \arrow{r} & I \arrow[equals]{d} \arrow{r} & A \arrow{d}{\lambda} \arrow{r}{q} & B \arrow{d}{\bar\lambda} \arrow{r} & 0\phantom{,} \\
				\bar{\mathsf e}\colon 0 \arrow{r} & I \arrow{r}& \mathcal M(I) \arrow{r}{\bar q}  & \mathcal Q(I) \arrow{r} & 0,
			\end{tikzcd}
		\end{equation}
		where $\tilde e$ is the pullback extension of $\mathsf e$ along $\Psi$, so
		\begin{equation}
			E = \{(a, d) \in A \oplus (C_0(0, 1] \otimes D) : q(a) = \Psi(d) \}
		\end{equation}
		and $\tilde\Psi$ and $\tilde q$ are the projection maps.
		
		Note that $\tilde{\mathsf e}$ is purely large.  Indeed, suppose $x \in I_+$, $e \in E_+ \setminus I$, and $\epsilon > 0$.  As $\Psi$ is faithful, $\tilde\Psi(e) \in A_+ \setminus I$, so there exists $v \in A$ such that $\|v^* \tilde\Psi(e) v - x \| < \epsilon$ by Theorem~\ref{pure-large}.  After multiplying $v$ on the left by an approximate unit for $I$, we may further assume $v \in I$, so $\|\Psi(v^* e v - x)\| = \|v^* \tilde\Psi(e) v - x \|< \epsilon$. As $\psi$ is totally faithful, $\Psi$ is isometric, so $\|v^* e v - x\| < \epsilon$. Hence Theorem~\ref{pure-large} implies $\tilde e$ is purely large.  By Elliott and Kucerovsky's absorption theorem (\cite[Corollary~16]{EK01} -- cf.\ \cite{G16}), $\mathsf{\tilde e}$ is an absorbing extension.  Note that the Busby invariant of $\tilde{\mathsf e}$ is $\bar \lambda \circ \Psi$, so $\bar \lambda \circ \Psi$ is absorbing as a morphism into $\mathcal Q(I)$, and hence any lift of $\bar \lambda \circ \Psi$ to a $^*$-homomorphism $C_0(0, 1] \otimes D \rightarrow \mathcal M(I)$ is absorbing; see \cite[Definition~5.9 and Corollary~5.11]{CGSTW23}, for example.
		Let $\Phi, \tilde\Phi \colon C_0(0, 1] \otimes A \rightarrow  A$ be the $^*$-homomorphisms induced by $\phi$ and $\tilde \phi$.  Then $\lambda \circ \Phi$ and $\lambda \circ \tilde \Phi$ are both lifts of $\bar \lambda \circ \Psi$, so they are absorbing.  Further, $\lambda \circ \Phi$ and $\lambda \circ \tilde \Phi$ induce an element $[\lambda \circ \Phi, \lambda \circ \tilde \Phi]$ in $KK(C_0(0, 1] \otimes D, I))$.  This $KK$-group vanishes as $C_0(0, 1] \otimes A$ is contractible, and, in particular, $[\lambda \circ \Phi, \lambda \circ \tilde\Phi] = 0$.  Then \cite[Theorem~A]{S26} implies there is a continuous path of unitaries $(u_t)_{t \geq 0}$ in the minimal unitization of $I$ such that
		\begin{equation}
			\lim_{t \rightarrow \infty} \|u_t \lambda(\Phi(\xi)) u_t^* - \lambda(\tilde\Phi(\xi))\| = 0, \quad \xi \in C_0(0, 1] \otimes D.
		\end{equation}
		Using that $\lambda|_I$ is the identity map and restricting to elements of the form $z \otimes d \in C_0(0, 1] \otimes D$ for $d \in D$ gives \eqref{eq:unique-lift}.
	\end{proof}
	
	Note that for an extension as in \eqref{eq:ext-section3}, there is an induced extension 
	\begin{equation}
		\begin{tikzcd}
			\mathsf e_\omega \colon 0 \arrow{r} & I_\omega \arrow{r} & A_\omega \arrow{r}{q_\omega} & B_\omega \arrow{r} & 0.
		\end{tikzcd}
	\end{equation}
	Indeed, since contractions in $B$ lift to contractions in $A$, it is clear that $q_\omega$ is surjective.  Also, $I_\omega$ is naturally identified with an ideal of $A_\omega$ and is contained in $\ker(q_\omega)$.  Finally, if $a \in \ker(q_\omega)_+$, lift $a$ to a sequence of positive elements $(a_n)_{n=1}^\infty \subseteq A$.  Then $\|q(a_n)\| \rightarrow 0$, so $((a_n - \|q(a_n)\|)_+)_{n=1}^\infty \subseteq A$ represents $a$, and $(a_n - \|q(a_n)\|)_+\in I$ for all $n \in \mathbb N$.  So $a \in I_\omega$.
	
	To prove Theorem~\ref{thm:D}, we will need to apply Theorem~\ref{thm:unique-lift} to the ultrapower of the given extension, which will be highly non-separable. The following separabilization result from \cite{CGSTW26} will allow us to do this.  The proof is standard using the characterization of pure largeness in Theorem~\ref{pure-large}.
	\begin{prop}\label{prop:pure-large-elementary}
		If $\mathsf e \colon 0 \to I \to A \to B \to 0$ is a purely large extension and $I_0 \subseteq I_\omega$, $A_0 \subseteq A_\omega$, and $B_0 \subseteq B_\omega$ are separable $C^*$-algebras, then there is a separable subextension $\mathsf e_1 \subseteq \mathsf e_\omega$ as in  
		\begin{equation}
			\begin{tikzcd}
				\mathsf e_1 \colon 0 \arrow{r} & I_1 \arrow[tail]{d} \arrow{r} & A_1\arrow[tail]{d} \arrow{r} & B_1 \arrow[tail]{d} \arrow{r} & 0 \\
				\mathsf e_\omega \colon 0 \arrow{r} & I_\omega \arrow{r} & A_\omega \arrow{r} & B_\omega \arrow{r} & 0
			\end{tikzcd}
		\end{equation}
		such that $I_1$ is stable, $\mathsf e_1$ is purely large, $I_0 \subseteq I_1$, $A_0 \subseteq A_1$, and $B_0 \subseteq B_1$.
	\end{prop}
	
	Combining Theorem~\ref{thm:unique-lift} and Proposition~\ref{prop:pure-large-elementary} gives the following result on the uniqueness of order zero lifts.
	
	\begin{cor}\label{cor:unique-lift}
		Suppose $0 \to I \to A \xrightarrow q B \to 0$ is a purely large extension of separable $C^*$-algebras with $I$ stable and let $D$ be a separable nuclear $C^*$-algebra.  If $\phi, \tilde\phi \colon D \rightarrow A_\omega$ are cpc$_\perp$ maps with $q_\omega \circ \phi = q_\omega \circ \tilde\phi$ and this common map is totally faithful, then there is a unitary $u$ in the minimal unitization of $I_\omega$ such that $\mathrm{ad}(u) \circ \phi = \tilde\phi$.
	\end{cor}
	
	\begin{proof}
		By Proposition~\ref{prop:pure-large-elementary}, there is a purely large extension as in the top row of the diagram
		\begin{equation}
			\begin{tikzcd}
				0 \arrow{r} & I_1 \arrow[tail]{d} \arrow{r} & A_1 \arrow[tail]{d} \arrow{r}{q_1} & B_1 \arrow[tail]{d} \arrow{r} & 0\phantom{,}\\
				0 \arrow{r} & I_\omega \arrow{r} & A_\omega \arrow{r}{q_\omega} & B_\omega \arrow{r} & 0,
			\end{tikzcd}
		\end{equation}
		with $I$ stable and with $\phi(D) \cup \tilde\phi(D) \subseteq A_1$. Note that $q_1 \circ \phi = q_1\circ \tilde\phi$ and this common map is totally faithful, so Theorem~\ref{thm:unique-lift} shows $\phi$ and $\tilde \phi$ are approximately unitarily equivalent as maps into $A_1$ with unitaries in the minimal unitization of $I_1$. A reindexing argument gives the result; see the proof of \cite[Lemma~3.15]{CGSTW23}.
	\end{proof}
	
	Kucerovsky and Ng introduced the \emph{corona factorization property} in \cite[Definition~2.1]{KN06} as a method for verifying pure largeness.  In particular, they show that a full extension of separable $C^*$-algebras is purely large whenever the ideal is stable and has the corona factorization property.  Ortega, Perera, and R{\o}rdam showed in \cite{OPR12} that the corona factorization property follows from various weak forms of strict comparison, including $n$-comparison for some $n \in \mathbb N$.  This includes $C^*$-algebras with finite nuclear dimension as shown by Robert in \cite{R11}, leading to the following theorem.
	
	\begin{thm}[{\cite{KN06} and \cite[Corollary~3.5]{R11}}]\label{thm:CFP}
		If $\mathsf e \colon 0 \to I \to A \to B \to 0$ is a full extension of separable $C^*$-algebras such that $I$ is stable and has finite nuclear dimension, then $\mathsf e$ is purely large. 
	\end{thm}
	
	Note that Theorem~\ref{thm:D} is an immediate consequence of Corollary~\ref{cor:unique-lift} and Theorem~\ref{thm:CFP}. We point out that $\ndim{I}<\infty$ in Theorem~\ref{thm:D} is a stronger condition than necessary; we only need that $I$ has the corona factorization property so as to access pure largeness and absorption.
	
	\section{Proof of Theorem~\ref{intro-thm:extension-result}}\label{sec:maintechnicalresult}
	
	Our main goal is to prove Theorem~\ref{intro-thm:extension-result}, but in order to set the stage for the proofs of Theorems~\ref{intro-thm:graphs-with-dim-one} and~\ref{intro-thm:dim-at-most-two}, some of the preliminary results, specifically Lemma~\ref{lem:key-lemma}, will work with a weaker condition than pure largeness.  
	For an extension
	\begin{equation}
		\begin{tikzcd}
			0 \arrow{r} & I \arrow{r} & A \arrow{r} & B \arrow{r} & 0,	
		\end{tikzcd}
	\end{equation}
	let $I^\perp \trianglelefteq A$ denote the two-sided annihilator of $I$ in $A$.  We will make the identifications
	\begin{equation}
		I = (I + I^\perp)/I^\perp \trianglelefteq A/I^\perp \quad \text{and} \quad I^\perp = (I + I^\perp) / I \trianglelefteq B.
	\end{equation}
	Then there is an induced extension
	\begin{equation}
		\begin{tikzcd}
			0 \arrow{r} & I \arrow{r} & A/I^\perp \arrow{r} & B/I^\perp \arrow{r} & 0.
		\end{tikzcd}
	\end{equation}
	
	The following is a very special case of Gabe's ideal-preserving lifting theorem from \cite{G22}.
	
	\begin{lem}[{\cite[Theorem~5.7]{G22}}]\label{lem:ideal-splitting}
		If $0 \to I \to A \xrightarrow{q} B \rightarrow 0$ is an extension of separable $C^*$-algebras and $B$ is nuclear, then there is a cpc map $\mu \colon B \rightarrow A$ such that $q \circ \mu  = \mathrm{id}_B$ and $\mu(I^\perp) = I^\perp$.
	\end{lem}
	
	\begin{proof}
		Take $X = \{0, 1\}$ to be the two-point topological space with only one non-trivial open set $U = \{0\}$.  View $I$, $A$, and $B$ as $X$-$C^*$-algebras with $I(U) = 0$, $A(U) = I^\perp$, and $B(U) = I^\perp$.  Then apply \cite[Theorem~5.7]{G22}.
	\end{proof}
	
	The following result is a standard consequence of the existence of quasicentral approximate units and Kirchberg's $\epsilon$-test (\cite[Lemma~3.1]{KR14}).
	
	\begin{prop}\label{prop:qcapproxunit}
		Let $I$ be an ideal in a separable $C^*$-algebra $A$.  If $S \subseteq A_\omega$ is a separable set, then there is a positive contraction $g \in I_\omega \cap S'$ such that $gx = x$ for all $x \in I_\omega \cap S$.
	\end{prop}
	
	We will need a technical variation of this proposition, which chooses the element $g$ with extra properties (Lemma~\ref{lem:constructing-g}).  The following stability result is a standard consequence of the Winter--Zacharias structure theorem for order zero maps (\cite{WZ09}) and Loring and Pedersen's (\cite{LP98}) result on projectivity of cones over finite-dimensional $C^*$-algebras.
	
	\begin{thm}[\cite{LP98, WZ09}]\label{thm:order-zero-stability}
		For every finite dimensional $C^*$-algebra $F$ and $\epsilon > 0$, there is $\delta > 0$ such that for every $C^*$-algebra $A$ and cpc map $\psi\colon F \to A$, if $\norm{\psi(a)\psi(b)}<\delta$ for all positive contractions $a,b \in F$ with $ab=0$, then there is a cpc$_\perp$ map $\phi \colon F \rightarrow A$ such that $\|\phi - \psi\| < \epsilon$.
	\end{thm}
	
	\begin{proof}
		Suppose the result fails for some finite dimensional $C^*$-algebra $F$ and $\epsilon > 0$.  For each $n \in \mathbb N$, set $\delta_n = 1/n$ and choose a counterexample $A_n$ and $\phi_n \colon F \rightarrow A_n$ for $\delta = \delta_n$.  Let $A_\omega = \prod_\omega A_n$ denote the ultrapower of the $A_n$ and let $\phi_\omega \colon F \rightarrow A_\omega$ denote the map induced by the sequence $(\phi_n)_{n=1}^\infty$.  By \cite[Corollary~3.8]{LP98} and \cite[Corollary~4.1]{WZ09}, there is a cpc$_\perp$ map $\psi \colon A \rightarrow \prod A_n$ lifting $\phi_\omega$.  Letting $\psi_n \colon F \rightarrow A_n$ denote the $n$th coordinate of $\psi$, we have $\psi_n$ is cpc$_\perp$ and
		\begin{equation}
			\lim_{n \rightarrow \omega} \|\phi_n(a) - \psi_n(a)\| = 0, \quad a \in  F.
		\end{equation}
		As the unit ball of $F$ is compact, for each $n \in \mathbb N$, there is a contraction $a_n \in F$ such that
		\begin{equation}
			\|\phi_n(a_n) - \psi_n(a_n)\| = \|\phi_n - \psi_n\|.
		\end{equation}
		Using compactness again, let $a = \lim_{n \rightarrow \omega} a_n \in F$.  Then for each $n \in \bN$,
		\begin{equation}
			\|\phi_n - \psi_n\| = \|\phi_n(a_n) - \psi_n(a_n)\| \leq \|\phi_n(a) - \psi_n(a) \| + 2\|a - a_n\| \rightarrow 0.
		\end{equation}
		So for some $n \in \mathbb N$, we have $\|\phi_n - \psi_n\| < \epsilon$, contradicting the choice of $\phi_n$ and $A_n$.
	\end{proof}
	
	The following will be needed to build an appropriate element $g \in I_\omega \cap A'$ as used in the definition of $\phi^0$ and $\phi^1$ in the sketch of the proof in the introduction.
	
	\begin{lem}\label{lem:constructing-g}
		Suppose $0 \to I \to A \xrightarrow q B \to 0$ is an extension of separable $C^*$-algebras and $\mu \colon B \rightarrow A$ is a cpc spitting of $q$.  
		Given a separable $C^*$-algebra $D$, a separable set $S \subseteq A_\omega$, 
		and cpc$_\perp$ maps $\psi^0, \psi^1 \colon D \rightarrow B_\omega$ with $\nndim{\psi^i} < \infty$ for $i \in \{0, 1\}$, there is a positive contraction $g \in I_\omega \cap S'$ such that $gx = x$ for all $x \in I_\omega \cap S$ and the maps
		\begin{equation}
			\phi^i \colon D \rightarrow A_\omega \colon d \mapsto (1 -g)\mu_\omega(\psi^i(d))
		\end{equation}
		are cpc$_\perp$ with $\nndim{\phi^i} \leq \nndim{\psi^i}$ for $i \in \{0, 1\}$.
	\end{lem}
	
	\begin{proof}
		Enlarging $S$ if necessary, we may assume 
		$C^*(\mu_\omega(\psi^0(D) \cup \psi^1(D))) \subseteq S$.  Let $(s_m)_{m=1}^\infty \subseteq S$ and $(t_m)_{m=1}^\infty \subseteq I_\omega \cap S$ be dense sequences and represent each $s_m$ and $t_m$ by sequences $(s_{m, n})_{n=1}^\infty \subseteq A$ and $(t_{m, n})_{n=1}^\infty \subseteq I$.
		
		Set $k_i = \nndim{\psi^i} < \infty$ and choose cp maps $(\psi_n^i \colon D \rightarrow B)_{n=1}^\infty$ representing $\psi^i$ with $\ndim{\psi^i} \leq k_i$.  Perturbing the $\psi_n^i$, we may assume there are finite dimensional $C^*$-algebras $F_n^i = F_n^{i,0} \oplus \cdots \oplus F_n^{i,k_i}$, cpc maps $\alpha_n \colon D \rightarrow F_n^i$, and cpc$_\perp$ maps $\beta_n^{i,j} \colon F_n^{i,j} \rightarrow B$ such that $\psi_n^i = \beta_n^i \circ \alpha_n^i$, where $\beta_n^i = \beta_n^{i,0} \oplus \cdots \beta_n^{i,k_i}$.  Let $\delta_n^{i, j} > 0$ be given by applying Theorem~\ref{thm:order-zero-stability} with $F = F_n^{i, j}$ and $\epsilon = 1/n$.
		
		Let $(e_l)_{l=1}^\infty \subseteq I$ be an approximate unit for $I$ which is quasicentral in $A$ with each $e_l$ positive and contractive.  For each $n \in \mathbb N$, $i\in\{ 0,1\}$, $0 \leq j \leq k_i$, and $x, y \in (F_n^{i,j})_+$ with $xy = 0$, we have
		\begin{equation}
			\begin{split}
				\limsup_{l \rightarrow \infty} \|(1 - e_l)^{1/2} \mu(\beta_n^{i, j}(x))& (1 - e_l) \mu(\beta_n^{i, j}(y)) (1 - e_l)^{1/2}\| \\
				&=\limsup_{l \rightarrow \infty} \|(1-e_l)^2 \mu(\beta_n^{i, j}(x)) \mu(\beta_n^{i, j}(y)) \| \\
				&\leq \limsup_{l \rightarrow \infty} \|(1-e_l) \mu(\beta_n^{i, j}(x)) \mu(\beta_n^{i, j}(y)) \| \\
				&= \|q(\mu(\beta_n^{i, j}(x))
				\mu(\beta_n^{i, j}(y)))\| \\
				&= \|\beta_n^{i, j}(x) \beta_n^{i, j}(y)\| = 0.
			\end{split}
		\end{equation}
		So after passing to a suitable subsequence of $(e_{l_n})_{n=1}^\infty$, setting $g_n = e_{l_n}$, and defining
		\begin{equation}
			\mu_n \colon B \rightarrow A \colon b \mapsto (1 - g_n)^{1/2} \mu(b) (1 - g_n)^{1/2},
		\end{equation}
		we can arrange that for all $n \in \mathbb N$,
		\begin{enumerate}
			\item\label{it:f1} $\|g_n s_{m, n} - s_{m, n} g_n\| < 1/n$ for $1 \leq m \leq n$,
			\item\label{it:f2} $\|g_n t_{m, n} - t_{m, n}\| < 1/n$ for $1 \leq m \leq n$,
			\item\label{it:f3} $\|\mu_n(\beta_n^{i, j}(x)) \mu_n(\beta_n^{i, j}(y)) \| < \delta_n^{i, j}$ for all $i \in\{ 0, 1\}$, $0 \leq j \leq k_i$, and $x, y \in (F_n^k)_+$ with $\norm{x},\norm{y}\leq 1$ and $xy = 0$.
		\end{enumerate}
		
		The sequence $(g_n)_{n=1}^\infty \subseteq I$ defines a positive contraction $g \in I_\omega$.  For all $m \in \mathbb N$, \ref{it:f1} implies $g s_m = s_m g$ and \ref{it:f2} implies $gt_m = t_m$, so $g \in I_\omega \cap S'$ and $gx = x$ for all $x \in I_\omega \cap S$.  Since $g$ commutes with the image of $\mu_\omega \circ \psi^i$, $\phi^i$ is cpc.  Note that if $a, b \in D_+$ satisfy $ab = 0$, then $\mu_\omega(\psi^i(a)) \mu_\omega(\psi^i(b)) \in I_\omega\cap S$, so
		\begin{equation}
			\begin{split}
				\phi^i(a) \phi^i(b) &= (1 - g) \mu_\omega(\psi^i(a)) (1 -g) \mu_\omega(\psi^i(b)) \\ &= (1 -g)^2 \mu_\omega(\psi^i(a)) \mu_\omega(\psi^i(b))\\ & = 0,
			\end{split}
		\end{equation}
		so $\phi^i$ is cpc$_\perp$.
		
		For each $n\in \mathbb N$, $0 \leq j\leq k_i$, and $i = 0, 1$, by \ref{it:f3} and the choice of $\delta_n^{i, j}$, there is a cpc$_\perp$ map $\tilde \beta_n^{i, j} \colon F_n^{i,j} \rightarrow A$ such that $\|\tilde\beta_n^{i,j} - \mu_n \circ \beta_n^{i,j}\| < 1/n$.  Define $\tilde\beta_n^i  = \tilde\beta_n^{i, 0} \oplus \cdots \oplus \tilde\beta_n^{i, k_i}$.  Then $\phi_n^i = \tilde\beta_n^i \circ \alpha_n \colon D \rightarrow A$ has nuclear dimension at most $k_i$.  It suffices to show the sequence $(\phi_n^i)_{n=1}^\infty$ represents $\phi^i$.  To this end, fix a contraction $a \in A$ and note that
		\begin{equation}
			\begin{split}    
				(1 - g_n)^{1/2} \mu(\psi_n^i(a)) (1 - g_n)^{1/2}
				&=  \mu_n(\beta_n^i(\alpha_n^i(a))) \\
				&\approx_{k_i/n} \tilde\beta_n^i(\alpha_n^i(a)) \\
				&= \phi_n^i(a),
			\end{split}
		\end{equation}
		completing the proof.
	\end{proof}
	
	\begin{rmk}
		If the conclusion of Lemma~\ref{lem:constructing-g} is changed to $\ndim{\phi} \leq \ndim{\psi}$, the proof becomes much easier.  Indeed, let $C \subseteq B_\omega$ be a separable $C^*$-algebra containing the range of $\psi$ such that the corestriction of $\psi$ to $C$ has nuclear dimension $\ndim{\psi}$.  Let $T = C^*(S \cup \mu_\omega(C)) \subseteq A_\omega$ and choose $g \in I_\omega \cap T'$ such that $gx = x$ for all $x \in I_\omega \cap T$ by Proposition~\ref{prop:qcapproxunit}. For $a, b \in C^*(\psi(D))_+$ with $ab = 0$, $q_\omega(\mu_\omega(a) \mu_\omega(b)) = ab= 0$, so $\mu_\omega(a)\mu_\omega(b) \in I_\omega \cap T$.   Hence 
		\begin{equation}
			(1 -g) \mu_\omega(a) (1 - g)\mu_\omega(b) = (1 - g)^2 \mu_\omega(a) \mu_\omega(b) = 0.
		\end{equation}
		So $(1 - g)\mu_\omega|_{C^*(\psi(D))} \colon C^*(\psi(D)) \rightarrow B_\omega$ has order zero, and hence $\ndim{\phi} \leq \ndim{\psi}$.  It is not immediately clear, however, if postcomposing with an order zero map decreases the value of $\nndim{\,\cdot\,}$, since the cpc$_\perp$ map $(1-g)\mu_\omega(\,\cdot\,)$ need not lift to a sequence of cpc$_\perp$ maps.
	\end{rmk}
	
	The next result is the main technical lemma behind the proof of Theorem~\ref{intro-thm:extension-result}, which also feeds into the proofs of Theorems~ \ref{intro-thm:graphs-with-dim-one} and~\ref{intro-thm:dim-at-most-two}.  The construction of the maps $\tilde\psi^0$ and $\tilde\psi^1$ are motivated by the proof that the Cuntz--Toeplitz algebras have nuclear dimension one in \cite{15A20}.
	
	\begin{lem}\label{lem:key-lemma}
		Let $0 \to I \to A \xrightarrow{q} B \to 0$ be an extension of separable $C^*$-algebras such that $I$ is stable, $B$ is nuclear, and $0 \to I \to A / I^\perp \to B / I^\perp \to 0$ is a purely large extension.  Suppose $\psi^0, \psi^1 \colon B \rightarrow B_\omega$ are cpc$_\perp$ maps such that
		\begin{enumerate}
			\item $\psi^0(b) + \psi^1(b) = b$ for all $b \in B$, 
			\item $\psi^0(I^\perp) \subseteq (I^\perp)_\omega$, and 
			\item the induced map $\bar\psi^0 \colon B/I^\perp \rightarrow (B/I^\perp)_\omega$ is totally faithful.
		\end{enumerate}
		Suppose further $f \in I_\omega \cap A'$ is such that $fx = x$ for all $x \in I$.
		Then there are cpc$_\perp$ maps $\tilde\psi^0, \tilde\psi^1 \colon A \rightarrow A_\omega$ such that
		\begin{enumerate}[resume]
			\item\label{it:key-lemma4} $\tilde\psi^0(a) + \tilde\psi^1(a) = (1 - f)a$ for all $a \in A$,
			\item\label{it:key-lemma5} $f\tilde\psi^1(a) = \tilde\psi^1(a) f = 0$ for all $a \in A$,
			\item\label{it:key-lemma6} $\tilde\psi^0(I^\perp) \subseteq (I^\perp)_\omega$
			\item\label{it:key-lemma7} $q_\omega \circ \tilde\psi^i = \psi^i \circ q$ for $i = 0, 1$, and
			\item\label{it:key-lemma8} 
			$\nndim{\tilde\psi^i} \leq \nndim{\psi^i}$ for $i = 0 , 1$.
		\end{enumerate}
	\end{lem}
	
	\begin{proof}
		Let $\mu \colon B \rightarrow A$ be a cpc splitting of $q$ with $\mu(I^\perp) = I^\perp$ as in Lemma~\ref{lem:ideal-splitting} and let $\mu_\omega \colon B_\omega \rightarrow A_\omega$ denote the induced cpc map.  
		Define the separable $C^*$-algebra
		\begin{equation}
			S=C^*(A \cup \mu_\omega(\psi^0(B))\cup \mu_\omega(\psi^1(B))\cup \{f\})\subset A_\omega.
		\end{equation}
		Apply Lemma~\ref{lem:constructing-g} to find $g \in I_\omega \cap S'$ such that $gx = x$ for all $x \in I_\omega \cap S$ and so that the maps
		$\phi^0, \tilde\psi^1 \colon A \rightarrow A_\omega$ given by
		\begin{align}
			\phi^0(a) &= (1-g) \mu_\omega(\psi^0(q(a)))
			\shortintertext{and}
			\tilde\psi^1(a) &= (1 -g) \mu_\omega(\psi^1(q(a)))
			\intertext{for $a \in A$ are cpc$_\perp$ with $\nndim{\phi^0}\leq \nndim{\psi^0}$ and $\nndim{\tilde\psi^1}\leq \nndim{\psi^1}$.  Then define}
			\tilde\psi^0(a) &= (g-f)a + \phi^0(a).
		\end{align}
		We claim that $\tilde\psi^0(a) + \tilde\psi^1(a) = (1-f)a$ for all $a \in A$. To see this, first note
		\begin{equation}
			\tilde\psi^0(a)+\tilde\psi^1(a)= (g-f)a+(1-g)\mu_\omega(q(a)).
		\end{equation}
		Since $q_\omega(\mu_\omega(q(a)))=q(a)$, we have $\mu_\omega(q(a))-a \in I_\omega \cap S$, so $(1-g)(\mu_\omega(q(a))-a)=0$. Therefore $(1-g)\mu_\omega(q(a))=(1-g)a$, from which it follows that
		\begin{equation}
			(g-f)+(1-g)\mu_\omega(q(a))=(g-f)a+(1-g)a=(1-f)a
		\end{equation}
		and \ref{it:key-lemma4} holds.  Since $a \mapsto (1-f)a$ is cpc$_\perp$, so is $\tilde\psi^0$.  Conditions \ref{it:key-lemma5}--\ref{it:key-lemma7} hold by construction and $\nndim{\tilde\psi^1} \leq \nndim{\psi^1}$ by the choice of $g$.  It remains to show $\nndim{\tilde\psi^0} \leq \nndim{\psi^0}$.  By construction $\nndim{\phi^0} \leq \nndim{\psi^0}$.  We will show $\phi^0$ and $\tilde\psi^0$ are unitarily equivalent, and then the result will follow from Proposition~\ref{prop:unitary-invariant}.  
		
		First, we define the quotient maps as in the diagram below.
		\begin{equation}\label{eq:pullback-annihilator}
			\begin{tikzcd}
				& & 0 \arrow{d} & 0 \arrow{d} & \\
				& & I^\perp \arrow{d} \arrow[equals]{r} & I^\perp \arrow{d} & \\
				0 \arrow{r} & I \arrow[equals]{d} \arrow{r} &
				A \arrow{r}{q} \arrow{d}[swap]{p} \arrow{dr}{r} & B \arrow{d}{\bar p} \arrow{r} & 0 \\
				0 \arrow{r} & I \arrow{r} & A/I^\perp \arrow{d} \arrow{r}[swap]{\bar q} & B/ I^\perp \arrow{d} \arrow{r} & 0 \\
				& & 0 & 0 & 
			\end{tikzcd}
		\end{equation}
		In this notation, we have
		\begin{equation}
			\bar\psi^0 \circ \bar p = \bar p_\omega \circ \psi^0.		
		\end{equation}
		
		Note that $p_\omega \circ \tilde\psi^0$ and $p_\omega\circ \phi^0$ vanish on $I + I^\perp$ since $\mu_\omega \circ \psi^0$ takes $I^\perp$ into $(I^\perp)_\omega$, 
        so there are cpc$_\perp$ maps $\tilde\theta, \theta \colon B/I^\perp \rightarrow (A/I^\perp)_\omega$ such that $p_\omega \circ \tilde\psi^0 = \tilde \theta \circ r$ and $p_\omega \circ  \phi^0 = \theta \circ r $.  
		Then $\tilde q_\omega \circ \theta = \tilde q_\omega \circ \tilde \theta = \bar\psi^0$.  As the lower row of \eqref{eq:pullback-annihilator} is purely large and $\bar \psi^0$ is totally faithful, there is a unitary $u$ in the minimal unitization of $I_\omega$ such that $\mathrm{ad}(u) \circ \theta = \tilde\theta$ by Corollary~\ref{cor:unique-lift}.  Viewing $u$ as a unitary in the minimal unitization of $A_\omega$, we have
		\begin{equation}
			p_\omega \circ \mathrm{ad}(u) \circ \phi^0 = \mathrm{ad}(p_\omega(u)) \circ p_\omega \circ \phi^0 = \mathrm{ad}(p_\omega(u)) \circ \theta \circ r =   \tilde \theta \circ r = p_\omega \circ \tilde\psi^0 
		\end{equation}
        and
		\begin{equation}
			q_\omega \circ \mathrm{ad}(u) \circ  \phi^0 = \mathrm{ad}(q_\omega(u)) \circ q_\omega \circ  \phi^0 = q_\omega \circ  \phi^0 = q_\omega \circ \tilde\psi^0.
		\end{equation}
		It follows that $\mathrm{ad}(u) \circ \phi^0 = \tilde\psi^0$.
	\end{proof}
	
	The following result proves the first part of Theorem~\ref{intro-thm:extension-result}.
	
	\begin{thm}[{cf.\ Theorem~\ref{intro-thm:extension-result}}]\label{body-thm:extension-result-1}
		If $0 \to I\to A\to B \to 0$ is a full extension of separable $C^*$-algebras with $I$ stable, then
		\begin{equation}
			\ndim{A}\leq\max\{\ndim{I}+1, \fndim{B}\}.
		\end{equation}
	\end{thm}
	
	\begin{proof}
		If $\ndim{I} = \infty$, there is nothing to prove.  If $\ndim{I} < \infty$, then the extension is purely large by Theorem~\ref{thm:CFP}. We will apply Lemma~\ref{lem:key-lemma}, noting that $I^\perp = 0$.	
		
		Let $\psi^0, \psi^1 \colon B \rightarrow B_\omega$ be cpc$_\perp$ maps with $\psi^0(b) + \psi^1(b) = b$, $\psi^0$ totally faithful, $\nndim{\psi^0} = 0$, and $\nndim{\psi^1} = \fndim{B} - 1$.  Further, let $f \in I_\omega \cap A'$ be such that $fx = x$ for all $x \in I$.  Let $\tilde\psi^0, \tilde\psi^1 \colon A \rightarrow A_\omega$ be the maps from Lemma~\ref{lem:key-lemma}.  Define $\tilde\phi^0, \tilde\phi^1 \colon A \rightarrow A_\omega$ by 
		\begin{equation}
			\tilde\phi^0(a) = \tilde\psi^0(a) \quad \text{and} \quad \tilde\phi^1(a) = fa + \tilde\psi^1(a), \qquad a \in A.
		\end{equation}
		Note that $\nndim{a \mapsto fa}\leq \ndim{I}$.  Using that the range of $\tilde\psi^1$ is orthogonal to $f$, Lemma~\ref{lem:sum}, and Proposition~\ref{prop:orthogonal-sum}, this implies
		\begin{equation}
			\begin{split} 
				\ndim{A} 
				&= \nndim{\iota_A} \\
				&\leq \max \{ \ndim{I}, \nndim{\tilde\psi^1}\} + \nndim{\tilde\psi^0} + 1 \\
                &\leq \max \{ \ndim{I}, \nndim{\psi^1}\} + \nndim{\psi^0} + 1\\
                &= \max \{ \ndim{I},\fndim{B} -1\} + 1\\
				&= \max\{\ndim{I} + 1, \fndim{B}\},
			\end{split}
		\end{equation}
		as required.
	\end{proof}
	
	When $I$ is a stable Kirchberg algebra, we will be able to obtain a stronger bound.  We isolate the following lemma from the proof for later use in Section~\ref{sec:no-sources}. The approach here is similar to that in \cite[Section~3.3]{ENSW25}.
	
	\begin{lem}\label{lem:kirchberg}
		Suppose $A$ is a separable nuclear $C^*$-algebra and $I \trianglelefteq A$ is a stable Kirchberg algebra.  If $S \subseteq I_\omega$ is a separable set, there are a positive contraction $f \in I_\omega \cap S'$ such that $fx = x$ for all $x \in I_\omega \cap S$ and cpc$_\perp$ maps $\tilde\theta^0, \tilde\theta^1 \colon A \rightarrow I_\omega$ such that
		\begin{enumerate}
			\item\label{it:kirchberg1} $\tilde\theta^0(a) + \tilde\theta^1(a) = fa$ for all $a \in A$,
			\item\label{it:kirchberg2} $f\tilde\theta^0(a) = \tilde\theta^0(a) f = \tilde\theta^0(a)$ for all $a \in A$, and
			\item\label{it:kirchberg3} $\ndim{\tilde\theta^0} = \ndim{\tilde\theta^1} = 0$.
			\item\label{it:kirchberg4} $\tilde\theta^0|_I$ is totally faithful.
		\end{enumerate}
	\end{lem}
	
	\begin{proof}
		Let $e \in I_\omega \cap A'$ be such that $ex = x$ for all $x \in I$ as in Proposition~\ref{prop:qcapproxunit}. Set $S = C^*(eA \cup\{e\})$. Since $I$ is $\mathcal O_\infty$-stable by \cite[Theorem~3.15]{KP00} and $S$ is separable, there is a unital embedding 
		\begin{equation}
			\mathcal O_\infty \into \frac{I_\omega \cap S'}{I_\omega \cap S^\perp};
		\end{equation}
		see \cite[Theorem~2.2]{TW07}. 
        
        Let $\bar h$ be a positive contraction in $\mathcal O_\infty$ with $\sigma(\bar h) = [0, 1]$ and let $h \in I_\omega \cap S'$ be a positive contraction lifting $\bar h$.  Let $T = C^*(S \cup A \cup \{e, h\}) \subseteq A_\omega$ and apply Proposition~\ref{prop:qcapproxunit} to find a positive contraction $f \in I_\omega \cap T'$ such that $fx = x$ for all $x \in I_\omega \cap T$.
		
		Define $\tilde\theta^0,\tilde\theta^1\colon A \to I_\omega$ by
		\begin{gather}
			\tilde\theta^0(a)=eha
			\shortintertext{and}
			\tilde\theta^1(a)=e(1-h)a+(f-e)a
		\end{gather}
		Clearly, \ref{it:kirchberg1} and \ref{it:kirchberg2} hold.
		
		In order to prove \ref{it:kirchberg3}, we will first show $\theta^0 = \tilde\theta^0|_I$ and $\theta^1 = \tilde\theta^1|_I$ are totally faithful; 
        note that this will also show \ref{it:kirchberg4}.  
        Exploiting the symmetry between $h$ and $1- h$, it suffices to show $\theta^0$ is totally faithful. 
        Let $\Theta^0 \colon C_0(0, 1] \otimes I \rightarrow I_\omega$ denote the $^*$-homomorphism induced by $\theta^0$.  
        Explicitly, $\Theta^0(f \otimes x) = f(h) x$ for all $f \in C_0(0, 1]$ and $x \in I$ by Remark~\ref{remark:multiplicationoncone}.  
        By Kirchberg's slice lemma (\cite[Lemma~4.1.9]{RordamClassificationBook}) it suffices to show that if $f \in C_0(0,1]$ and $x \in I$ are non-zero, then $\Theta^0(f \otimes x) \neq 0$.  
        Since $\bar h$ has full spectrum $f(\bar h) \neq 0$, and as $f(\bar h)$ is contained in the simple $C^*$-algebra $\mathcal O_\infty$, it is full in $(I_\omega \cap S')/(I_\omega \cap S^\perp)$.  
        In particular, since $x = ex \in S$, the ideal generated by $\Theta^0(f \otimes x) = f(h)x \in I_\omega \cap S'$ contains an element of the form $1- y$ for some $y \in I_\omega \cap S^\perp$.  
        As $x \in I$ and $I = eI \subseteq S$, we have $x(1- y) = x$, so $x$ is in the ideal generated by $\Theta^0(f \otimes x)$.  
        In particular,  $\Theta^0(f \otimes x) \neq 0$.
		
		Let $i \in \{0, 1\}$ be given. Note that $\tilde\theta^i$ vanishes on $I^\perp$ since $e, f \in I_\omega$ and hence induces a cpc$_\perp$ map $\bar \theta^i \colon A / I^\perp \rightarrow I_\omega$, which in turn induces a $^*$-homomorphism $\bar \Theta^i \colon C_0(0, 1] \otimes A/I^\perp \rightarrow I_\omega$.   Note that $\bar\Theta^i|_{C_0(0, 1] \otimes I}=\Theta^i$.  As $C_0(0, 1] \otimes I$ is an essential ideal in $C_0(0, 1] \otimes (A/I^\perp)$, $\bar\Theta^i$ is faithful.  Note that $\bar\Theta^i$ is full and $\mathcal O_\infty$-stable since $I$ is a Kirchberg algebra (\cite[Corollary~9.8]{Gabe2}).  As $C_0(0, 1] \otimes (A/I^\perp)$ is quasidiagonal (\cite{V91}), it follows from \cite[Theorem~D]{BGSW22}, that $\ndim{\bar\Theta^i} = 0$, so $\ndim{\tilde\theta^i} = 0$, proving \ref{it:kirchberg3}.
	\end{proof}
	
	We now use the previous lemma to obtain our improved bound in the second part of Theorem~\ref{intro-thm:extension-result}.
	
	\begin{thm}[{cf.\ Theorem~\ref{intro-thm:extension-result}}]\label{body-thm:extension-result-2}
		If $0 \to I \to A \to B \to 0$ is an essential extension of separable $C^*$-algebras such that $I$ is a stable Kirchberg algebra, then
		\begin{equation}
			\ndim{A} \leq \fndim{B}.
		\end{equation}
	\end{thm}
	
	\begin{proof}
		If $\fndim{B}=\infty$, there is nothing to prove. Thus we may assume $\fndim{B}<\infty$, so $B$ is nuclear. Since $I$ is a Kirchberg algebra, it follows that $A$ is nuclear. Essentiality of the extension implies the Busby map $A \rightarrow \mathcal Q(I)$ is faithful.  Since $\mathcal Q(I)$ is simple by \cite[Theorem~2.8]{Lin91}, the Busby map, and hence the extension, is full.  As $I$ is stable and has finite nuclear dimension, the extension is purely large by Theorem~\ref{thm:CFP}.
		
		Fix a positive contraction $e \in I_\omega \cap A'$ such that $ex = x$ for all $x \in I$, as in Proposition~\ref{prop:qcapproxunit}, and set $S = C^*(A \cup \{e\})$.  Let $f \in I_\omega \cap S'$ and $\tilde\theta^0, \tilde\theta^1 \colon A \rightarrow I_\omega$ be as in Lemma~\ref{lem:kirchberg}\ref{it:kirchberg1}.  By the definition of faithful nuclear dimension, there are cpc$_\perp$ maps $\psi^0, \psi^1 \colon B \rightarrow B_\omega$ adding to the canonical inclusion with $\psi^0$ totally faithful, $\nndim{\psi^0} = 0$, and $\nndim{\psi^1} = \fndim{B} - 1$.   Let $\tilde\psi^0$ and $\tilde\psi^1$ be the maps obtained by applying Lemma~\ref{lem:key-lemma}, noting that $I^\perp = 0$.
		
		Define $\phi^0, \phi^1 \colon A \rightarrow A_\omega$ by
		\begin{align}
			\phi^0(a) &= \tilde\theta^0(a) + \tilde\psi^0(a)
			\shortintertext{and}
			\phi^1(a) &= \tilde\theta^1(a) + \tilde\psi^1(a)
		\end{align}
		for all $a \in A$. Since $f$ acts as the unit on $\tilde\theta^0(a)$ and $\tilde\psi^0(a) \leq (1-f)a$ for all $a \in A$, $\tilde\theta^0$ and $\tilde\psi^0$ have orthogonal ranges.  So Proposition~\ref{prop:orthogonal-sum} implies
		\begin{equation}
			\nndim{\phi^0} =\max\{ \nndim{\tilde\theta^0}, \nndim{\tilde\psi^0} \}= 0.
		\end{equation}    
		Also note that $\tilde\theta^1(a) \leq fa$ and $(1-f)$ acts as the unit on $\tilde\psi^1(a)$ for all $a \in A$, so $\tilde\theta^1$ and $\tilde\psi^1$ have orthogonal ranges.  Then another application of Proposition~\ref{prop:orthogonal-sum} yields
		\begin{equation}
			\nndim{\phi^1} =\max\{ \nndim{\tilde\theta^1} , \nndim{\tilde\psi^1} \}= \fndim{B}-1.
		\end{equation}    
		Since $\phi^0(a) + \phi^1(a) = a$ for all $a \in A$, the result follows from Lemma~\ref{lem:sum} and Proposition~\ref{prop:nndimisndim}. 
	\end{proof}
	
	This completes the proof of the second part of Theorem~\ref{intro-thm:extension-result}.
	
	\section{Graph algebra prerequisites}\label{sec:graph-cstar-algebras}
	
	This section recalls some basic results and definitions from graph algebras, as well as proving some results needed to obtain Theorems~\ref{intro-thm:graphs-with-dim-one} and~\ref{intro-thm:dim-at-most-two}.  Our general reference for graph $C^*$-algebras is Raeburn's monograph \cite{R05}, and we follow the edge conventions from there.
	
	A \emph{directed graph} is a $4$-tuple $E=(E^0, E^1, r,s)$ where $E^0$ is the set of vertices, $E^1$ the set of edges, and $r,s \colon E^1 \to E^0$ are the range and source maps.  Paths will be written right to left \`a la function composition: $ef$ is a path if $r(f)=s(e)$.  For $n \geq 2$, let $E^n$ be the set of all paths $e_1 \cdots e_n$, where $e_1, \ldots, e_n \in E^1$ with $s(e_i) = r(e_{i+1})$ for $1 \leq i < n$, and let $E^* = \bigcup_{n=0}^\infty E^n$ be the set of finite paths.  We extend the range and source maps to maps of finite paths $r, s \colon E^* \rightarrow E^0$ in the obvious way.  For $v \in E^0$, let $E^*v = \{\mu \in E^* : s(\mu) = v \}$.
	
	Given a directed graph $E$, $C^*(E)$ is the universal $C^*$-algebra generated by mutually orthogonal projections $(p_v)_{v \in E^0}$ and partial isometries $(s_e)_{e \in E^0}$ such that \begin{enumerate}
		\item\label{it:CK1} $s_e^*s_e=p_{s(e)}$ and $s_e s_e^* \leq p_{r(e)}$ for all $e \in E^1$, and
		\item\label{it:CK2} $\displaystyle p_v=\sum_{e \in r^{-1}(v)} s_es_e^*$ whenever $r^{-1}(v)$ is non-empty and finite.
	\end{enumerate}
	If $v \in E^0$, we write $s_v = p_v \in C^*(E)$, and if $\mu = e_1 \cdots e_n \in E^n$ with $n \geq 2$, write $s_\mu = s_{e_1} \cdots s_{e_n} \in C^*(E)$.
	
	We start by describing the combinatorial conditions in Theorem~\ref{intro-thm:graphs-with-dim-one}.  The following definitions are standard, but there seems to be competing terminology in the literature, so we give the formal definitions we use here.
	
	\begin{defi}\label{def:isolatedcycles}
		Let $E$ be a directed graph.
		\begin{enumerate}
			\item A \emph{cycle} in $E$ is a path $e_1 \cdots e_n$ in $E$ with $n \geq 1$ such that $r(e_1) = s(e_n)$ and the vertices $s(e_1), \ldots, s(e_n)$ are distinct.
			\item A \emph{return path} in $E$ is a path $e_1 \cdots e_n$ in $E$ with $n \geq 1$ such that $r(e_1) = s(e_n)$ and $s(e_i) \neq s(e_n)$ for all $i = 1, \ldots, n-1$.
			\item A cycle $e_1 \cdots e_n$ is \emph{isolated} if the only return path $\mu$ with $s(\mu) = s(e_n)$ is $e_1 e_2 \cdots e_{n}$.
		\end{enumerate}
	\end{defi}
	
	Informally, both cycles and return paths are paths in the graph from a vertex $v \in E^0$ to itself.  A cycle cannot repeat any vertex in $E^0$, whereas a return path can repeat cycles and edges, but it does not return to $v$ before the end of the path.  An isolated cycle is a cycle with no other return path based on the cycle.  Note that with this terminology, Condition~$(K)$ states that there are no isolated cycles.
	
	We will be particularly interested in the ideal structure of graph algebras, but only in the gauge-invariant ideals. To this end, recall that for a directed graph $E$, there is a canonical action $\rho \colon \mathbb T \curvearrowright C^*(E)$, called the \emph{gauge action}, defined  by
	\begin{equation}
		\rho_\lambda(p_v)=p_v \quad \text{and} \quad \rho_\lambda(s_e)=\lambda s_e, \qquad v \in E^0,\ e \in E^1,\ \lambda \in \mathbb T.
	\end{equation}
	We will only need the ideal structure for $C^*$-algebras of \emph{row-finite graphs}; i.e.\ those $E$ for which $r^{-1}(v)$ is finite for all $v \in E^0$.
	
	\begin{defi}
		For a row-finite directed graph $E$, a set $H \subseteq E^0$ is said to be 
		\begin{enumerate}
			\item \emph{hereditary} if for every path $\mu$ in $E^*$ with $r(\mu) \in H$, we have $s(\mu) \in H$.
			\item \emph{saturated} if $v \in H$ whenever $r^{-1}(v)\neq\emptyset$ and $s(r^{-1}(v))\subseteq H$.
		\end{enumerate}
	\end{defi}
	
	Two subgraphs are associated to each hereditary and saturated subset $H \subseteq E^0$, namely $E_H$ and $E\setminus H$.  Their vertex and edge sets are as follows:
	\begin{equation}
		\begin{aligned}
			E_H^0&=H; &\qquad (E\setminus H)^0&=E^0\setminus H; \\
			E_H^1&=\{e \in E^1: r(e) \in H); &\qquad (E\setminus H)^1&=\{e \in E^1: s(e) \notin H\}.
		\end{aligned}
	\end{equation}
	The range and source maps for $E_H$ and $E\setminus H$ are simply the restrictions of $r$ and $s$ to $E^1_H$ and $(E\setminus H)^1$. These are well defined because $H$ is hereditary.
	
	For a hereditary and saturated subset $H\subseteq E^0$, let $I_H$ denote the ideal of $C^*(E)$ generated by $\{p_v: v\in H\}$.
	
	\begin{thm}[{\cite[Theorem~6.6]{KPRR97}, cf.\ \cite[Theorem~4.9]{R05}, combined with \cite[Theorem~2.8]{B77}}]\label{thm:graph-ideals}
		For a countable row-finite graph $E$, the assignment $H \mapsto I_H$ is a bijection between the sets of hereditary and saturated subsets of $E^0$ and gauge-invariant ideals of $C^*(E)$. Furthermore, $I_H$ is stably isomorphic to $C^*(E_H)$, and $C^*(E)/I_H \cong C^*(E\setminus H)$.
	\end{thm}
	
	We now consider the ideal generated by a source.  Recall that a vertex $v \in E^0$ is a \emph{source} if $r^{-1}(v) = \emptyset$.   For a set $S$, let $(\delta_s)_{s \in S}$ be the standard orthonormal basis for the Hilbert space $\ell^2(S)$ and for $\xi, \eta \in \ell^2(S)$, let $\xi\eta^* \in \mathcal K(\ell^2(S))$ denote the rank one operator $\xi\eta^*(\zeta) = \langle \zeta, \eta \rangle \xi$.  The following result is standard. See {\cite[Section~4, Lemma~B]{FC24}}, for example.
	
	\begin{prop}\label{prop:source}
		If $E$ is a row-finite graph and $v \in E^0$ is a source, then there is an isomorphism
		\begin{equation}
			\phi \colon \mathcal K(\ell^2(E^*v)) \overset\cong\longrightarrow I_v
		\end{equation}
		given by $\phi(\delta_\mu \delta_\nu^*) = s_\mu s_\nu^*$ for $\mu, \nu \in E^*v$.
	\end{prop}
	
	We next describe the ideal generated by an entryless cycle.  More precisely, let $\mu = e_1 \cdots e_n$ be an entryless cycle in the graph $E$ and let $I_\mu$ denote the ideal of $C^*(E)$ generated by $p_{s(\mu)}$.  Let $E_\mu$ denote the graph given by removing the final edge of $\mu$, so $E_\mu^0 = E^0$ and $E_\mu^1 = E^1 \setminus \{e_1\}$, with the range and source maps being the restrictions of the corresponding maps for $E$.  The following result is certainly known to experts, but we have been unable to find a suitable reference.
	
	\begin{prop}\label{prop:entrylesscyle}
		If $E$ is a row-finite graph with an entryless cycle $\mu$, then there is an isomorphism
		\begin{equation}
			\phi \colon  C(\mathbb T) \otimes \mathcal K(\ell^2(E_\mu^* s(\mu)))  \overset\cong\longrightarrow I_\mu
		\end{equation}
		given by $\phi(z \otimes \delta_{s(\mu)}\delta_{s(\mu)}^*) = s_\mu$ and $\phi(1_{C(\mathbb T)} \otimes \delta_{\nu} \delta_{\nu'}^*) = s_\nu s_{\nu'}^*$ for $\nu, \nu' \in E_\mu^*s(\mu)$, where $z \in  C(\mathbb T)$ is the inclusion $z \colon \mathbb T\hookrightarrow \mathbb C$.
	\end{prop}
	
	\begin{proof}
		We define a graph $F$ as follows.  Set 
		\begin{equation}
			F^0 = E_\mu^* s(\mu) \quad \text{and} \quad F^1 = \{ (e, \nu) \in E^1 \times E_\mu^*s(\mu) : s(e) = r(\nu) \} \cup \{\mu\}
		\end{equation}
		with range and source given by
		\begin{equation}
			r(e, \nu) = e\nu, \quad s(e, \nu) = \nu, \quad \text{and} \quad r(\mu) = s(\mu)
		\end{equation}
		where $r(\mu)=s(\mu)$ is the usual range and source vertex of $\mu$ in $E$, instead considered as a length-zero path in $E_\mu^*s(\mu)$.
        
        We will construct isomorphisms
		\begin{equation}
			\begin{tikzcd}
				C(\mathbb T) \otimes \mathcal K(\ell^2(E_\mu^* s(\mu))) & C^*(F) \arrow{l}{\cong}[swap]{\theta} \arrow{r}{\psi}[swap]{\cong} & I_\nu
			\end{tikzcd}
		\end{equation}
		and define $\phi = \psi \circ \theta^{-1}$.
		
		Define $\psi$ by
		\begin{equation}
			\psi(p_\nu) = s_\nu s_\nu^*, \quad \psi(s_{(e, \nu)}) = s_e s_\nu s_\nu^*, \quad \text{and} \quad \psi(s_\mu) = s_\mu.
		\end{equation}
		for all $\nu \in F^0$ and $(e, \nu) \in F^1\setminus \{\mu\}$.
		It is easy to see that $\psi$ is well-defined, noting that $s_\mu s_\mu^*=s_\mu^*s_\mu=p_{s(\mu)}$ since $\mu$ is entryless.  Note that $\psi(p_\nu) \neq 0$ for all $\nu \in F^0$, $F$ has a unique entryless cycle $\mu$, and the spectrum of $\psi(s_\mu)$ contains $\mathbb T$.  By Szyma\'{n}ski's uniqueness theorem, $\psi$ is injective.  
		
		To show surjectivity, let $H$ be the saturated hereditary subset of $E^0$ generated by $s(\mu)$, and note that $I_\mu = I_H$.  By \cite[Lemma~4.3]{BPRS00}, $I_H$ is the closed span of elements of the form $s_\alpha s_\beta^*$ for $\alpha, \beta \in E^*$ with $s(\alpha) = s(\beta)\in H$; as this common source vertex is in $H$, there is at least one path from $s(\mu)$ to $s(\alpha)$ in $E_\mu$. Since $E_H$ only contains one cycle, $\mu$, the hereditary and saturated subset $\tilde{H}\subseteq F^0$ generated by $s(\mu)$ does not contain any cycles. Thus, there are only finitely many such paths, say $\nu_1, \ldots, \nu_k \in E_\mu^*$ from $s(\mu)$ to $s(\alpha)$.  Then, since $\mu$ is an entryless cycle,
		\begin{equation}
			s_\alpha s_\beta^* = \sum_{i = 1}^k s_{\alpha \nu_i} s_{\beta \nu_i}^* \in C^*(E).
		\end{equation}
		So $I_\mu$ is the closed span of elements of the form $s_\alpha s_\beta^*$ with $\alpha, \beta \in E^*$ with $s(\alpha) = s(\beta) = s(\mu)$. 
		
		Let $\nu = e_1 \cdots e_k \in E^*$ be a path with $s(\nu) = s(\mu)$, and suppose that $\nu$ does not factor as $\nu'\mu$. For $i = 1, \ldots, k$, set $\nu_i = e_i \cdots e_k$, so $e_i \nu_{i+1} = \nu_i$ for all $i = 1, \ldots, k - 1$.  Note that
		\begin{equation}
			\psi(s_{(e_i, \nu_{i+1})}) = s_{e_i} s_{\nu_{i+1}} s_{\nu_{i+1}}^* = s_{\nu_i} s_{\nu_{i+1}}^*, \quad \text{and} \quad s(p_{e_k}) = s_{\nu_k} s_{\nu_k}^*.
		\end{equation}
		for $i = 1, \ldots, k -1$.  Hence
		\begin{equation}\label{eq:compute-psi}
			\psi(s_{(e_1, \nu_2)}) \cdots \psi(s_{(e_{k-1}, \nu_k)}) \psi(p_{e_k}) = s_{\nu_1} s_{\nu_2}^* \cdots s_{\nu_{n-1}} s_{\nu_k}^* s_{\nu_k} s_{\nu_k^*} = s_\nu.
		\end{equation}
		Now for any $\nu \in E^*$ with $s(\nu)=s(\mu)$, we may write $\nu=\nu'\mu^k$ for some $k \in \bN$, where $\nu'$ does not factor as $\nu''\mu$. Then $s_{\nu}=s_{\nu'}s^k_\mu$ is in the range of $\psi$. We conclude that $\psi$ is surjective.

		Define $\theta$ by
		\begin{equation}
			\theta(p_\nu) = 1_{C(\mathbb T)} \otimes \delta_\nu \delta_\nu^*,\ \theta(s_{(e, \nu)}) = 1_{C(\mathbb T)} \otimes \delta_{e\nu} \delta_{\nu}^*,\ \text{and}\ \theta(s_\mu) = z \otimes \delta_{s(\mu)} \delta_{s(\mu)}^*.
		\end{equation}
		Again, it is easy to check that $\theta$ is well-defined, and using Szyma\'{n}ski's uniqueness theorem as above, $\theta$ is injective. Surjectivity of $\theta$ is clear.
		
		By construction, the isomorphism $\phi = \psi \circ \theta^{-1}$ satisfies $\phi(z \otimes \delta_{s(\mu)} \delta_{s(\mu)}^*) = s_{\mu}$.  Fix $\nu \in E_\mu^*$ with $s(\nu) = s(\mu)$ and write $\nu = e_1 \cdots e_k$. As before, for $1 \leq i \leq k$, let $\nu_i  = e_i \cdots e_k$.  Using \eqref{eq:compute-psi}, we have
		\begin{equation}
			\psi^{-1}(s_{\nu}) = s_{(e_1, \nu_2)} s_{(e_2, \nu_3)} \cdots s_{e_{(k-1}, \nu_k)} p_{e_k} \in C^*(F). 
		\end{equation}
		Applying $\theta$ yields $\phi^{-1}(s_\nu) = 1_{C(\bT)} \otimes \delta_{\nu}\delta_{s(\mu)}^*$.  Then for all $\nu, \nu' \in E_\mu^*s(\mu)$, we have
		\begin{equation}
			\phi(1_{ C(\bT)} \otimes \delta_\nu \delta_\nu'^*) = \phi(1_{ C(\bT)} \otimes \delta_\nu \delta_{s(\mu)}^*)\phi(1_{C(\bT)} \otimes \delta_{\nu'} \delta_{s(\mu)}^*)^* = s_{\nu} s_{\nu'}^*,
		\end{equation}
		so $\phi$ has the required properties.
	\end{proof}
	
	The following result gives a complete description of the possible minimal gauge-invariant ideals in a row-finite directed graph.
	
	\begin{prop}\label{prop:minimalideals}
		If $E$ is a finite graph and $I \trianglelefteq E$ is a minimal gauge-invariant ideal, then either $I$ is a Kirchberg algebra, $I$ is stably isomorphic to $\mathbb C$, or $I$ is stably isomorphic to $C(\mathbb T)$. Furthermore, $I$ is either stable or unital.
	\end{prop}
	
	\begin{proof}
		By Theorem~\ref{thm:graph-ideals}, there is a minimal hereditary saturated set $H \subseteq E^0$ such that $I = I_H$ and $I$ is stably isomorphic to $E_H$.  If $H$ contains a source $v \in E^0$, then $v$ must be the unique source in $H$, so $I_H \cong \mathcal K(\ell^2(E^*_{s(\mu)}v))$ by Proposition~\ref{prop:source}.  If $H$ contains an entryless cycle $\mu$, $I_H$ must equal $I_\mu$, and so $I_H \cong C(\mathbb T) \otimes \mathcal K(\ell^2(E^*_\mu s(\mu)))$ by Proposition~\ref{prop:entrylesscyle}.  Suppose $H$ contains no source or entryless cycle.  As $H$ is finite and has no source, $E_H$ contains a cycle.  Since $E_H$ has no gauge-invariant ideals, it must be cofinal in the sense of \cite{BPRS00}, and every cycle has an entry, so $C^*(E_H)$ is simple and purely infinite by \cite[Propositions~5.1 and 5.3]{BPRS00}. Graph algebras are nuclear by \cite{KP99} (also see \cite[Corollary~4.5.4]{BrownOzawa}), and $C^*(E_H)$ is unital and separable since $E_H$ is finite. It follows that $C^*(E_H)$ is $\cO_\infty$-stable by \cite[Theorem~3.15]{KP00}, and so is $I_H$ by \cite[Corollary~3.2]{TW07}. Since nuclearity and simplicity are also preserved by stable isomorphism, it follows that $I_H$ is a Kirchberg algebra. In all three cases, observe that $I_H$ is either stable or unital, with the Kirchberg case being \cite[Theorem~1.2(i)]{Z92}.
	\end{proof}
	
	Our next goal is to show that if $E$ is a finite graph and $I \trianglelefteq  C^*(E)$ is a minimal gauge-invariant ideal such that $I$ is stable, then $I$ is purely large in $C^*(E) / I^\perp$.  The difficult part is when $I \cong C(\mathbb T) \otimes \mathcal K$, and this case requires some further preliminaries.
	
	Fix a compact metrizable space $X$ and an essential extension
	\begin{equation} 
		\begin{tikzcd}
			\mathsf e \colon 0 \arrow{r} & C(X) \otimes \mathcal K \arrow{r} & A \arrow{r} & B \arrow{r} & 0.
		\end{tikzcd}
	\end{equation}
	For $x \in X$, the evaluation map $\mathrm{ex}_x \otimes \mathrm{id}_{\mathcal K}$ extends canonically to a $^*$-homomorphism $\widetilde{ \mathrm{ev}}_x \colon \mathcal M(C(X) \otimes \mathcal K) \rightarrow \mathcal M(\mathcal K)$, giving a diagram
	\begin{equation}
		\begin{tikzcd}
			\phantom{(x)}\mathsf e \colon 0 \arrow{r} & C(X) \otimes \mathcal K \arrow{d}{\mathrm{ev_x} \otimes \mathrm{id}_{\mathcal K}} \arrow{r} & A \arrow{d}{\widetilde{\mathrm{ev}}_x} \arrow{r} & B \arrow[equals]{d} \arrow{r} & 0\phantom{.} \\
			\mathsf e(x) \colon 0 \arrow{r} & \mathcal K \arrow{r} & \widetilde{\mathrm{ev}}_x(A) \arrow{r} & B \arrow{r} & 0.
		\end{tikzcd}
	\end{equation}
	
	This type of extension was originally studied by Pimsner, Popa, and Voiculescu in \cite{PPV79}. In particular, they said that an essential extension of this form is \emph{homogeneous} if $\mathsf{e}(x)$ is essential for all $x \in X$. Pimsner, Popa, and Voiculescu proved that these homogeneous extensions are nuclearly absorbing. The major breakthrough in \cite{EK01} was that nuclear absorption is equivalent to pure largeness. As a consequence, homogeneous extensions must be purely large -- a fact that was proven directly in \cite[Theorem~17(i)']{EK01}. We record this result below.
	
	\begin{thm}[{\cite[Theorem~17(i)']{EK01}, cf.\ \cite{PPV79}}]\label{thm:PPV}
		Let $X$ be a finite-dimensional compact metrizable space and let 
		\begin{equation} 
			\begin{tikzcd}
				\mathsf e \colon 0 \arrow{r} & C(X) \otimes \mathcal K \arrow{r} & A \arrow{r} & B \arrow{r} & 0.
			\end{tikzcd}
		\end{equation}
		be an essential extension of separable C$^*$-algebras.  If $\mathsf e(x)$ is essential for all $x \in X$, then $\mathsf e$ is purely large.  
	\end{thm}
	
	The following result will allow us to quote Lemma~\ref{lem:key-lemma} to lift nuclear dimension approximations from $C^*(E) / I$ to $C^*(E)$.
	
	\begin{thm}\label{thm:purelylargeannihilator}
		Let $E$ be a finite graph and $I \trianglelefteq E$ be a minimal gauge-invariant ideal.  If $I$ is stable, then
		\begin{equation}\label{eq:thm:purelylargeannihilator}
			\begin{tikzcd}
				\mathsf{e}\colon 0 \arrow{r} & I \arrow{r} & C^*(E) / I^\perp \arrow{r} & C^*(E)/(I + I^\perp) \arrow{r} & 0
			\end{tikzcd}
		\end{equation}
		is purely large.
	\end{thm}
	
	\begin{proof}
		By Proposition~\ref{prop:minimalideals}, we know that $I$ is isomorphic to either $\cK$, $C(\bT)\otimes \cK$, or a stable Kirchberg algebra. In the Kirchberg case, note that the corona $\cQ(I)$ is simple by \cite[Theorem~2.8]{Lin91}; the same is true for $I=\cK$. Since the extension $\mathsf{e}$ is automatically essential, the Busby map is injective and therefore full by simplicity of $\cQ(I)$. Thus, by Theorem~\ref{thm:CFP}, it follows that $\mathsf{e}$ is purely large when $I=\cK$ or is a stable Kirchberg algebra.
		
		For the remaining case, where $I\cong C(\bT)\otimes \cK$, we need to verify the point-wise essentiality of $\mathsf{e}$. We achieve this through a careful analysis of the graph. Since $I^\perp$ is gauge-invariant, $I^\perp=I_{H^\perp}$ where $H^\perp=\{v \in E^0: p_v \in I^\perp\}$. Then $C^*(E)/I^\perp \cong C^*(E\setminus H^\perp)$, and $H\subseteq (E \setminus H^\perp)^0$ since $H \cap H^\perp=\emptyset$. Therefore, without loss of generality, we may assume that $I^\perp=0$. Let $\mu=e_1e_2\cdots e_n$ be the entryless cycle that generates $I=I_\mu$; as a consequence of the fact that $I_\mu$ is essential in $C^*(E)$, each cycle in $E$ can be reached from $\mu$ by \cite[Proposition~3.4]{BFPR22}. We now construct a collection of Cuntz--Krieger $E$-families in $C^*(E_\mu)$, where the collection is indexed by the circle $\bT$. Let $s_e$ and $p_v$ denote the usual generators for $C^*(E_\mu)$, and let $\lambda \in \bT$. Then define
		\begin{equation}\label{eq:ckfamily}
			q_v=p_v \qquad t_{e_1}= \lambda s_{e_n}^*\cdots s_{e_2}^* \qquad t_{e}=s_e \text{ if $e\neq e_1$}.
		\end{equation}
		The main observation that verifies this to be a Cuntz--Krieger $E$-family is that $t_{e_1}t_{e_1}^*=p_{r(e_1)}$ and $t_{e_1}^*t_{e_1}=p_{s(e_1)}$, which is a straightforward computation relying upon the fact that $\mu$ is entryless. We thus obtain a $^*$-homomorphism $\pi_\lambda\colon C^*(E)\to C^*(E_\mu)$ given by the Cuntz--Krieger $E$-family in \eqref{eq:ckfamily}. 
		
		We now wish to verify that $\pi_\lambda$ agrees with $\mathrm{ev}_\lambda \otimes \mathrm{id}_\cK$. We do this via the isomorphism $\phi\colon  C(\bT)\otimes \cK\to I_\mu$ constructed in Proposition~\ref{prop:entrylesscyle} and the isomorphism $\psi\colon \cK\to I_{s(\mu)}$ given in Proposition~\ref{prop:source}. In particular, we seek to verify that the following diagram commutes:
		\begin{equation}\label{eq:circleidealdiagram}
			\begin{tikzcd}
				C(\bT)\otimes \cK(\ell^2(E_\mu^* s(\mu))) \arrow[d, "\mathrm{ev}_\lambda\otimes \mathrm{id}_\cK"'] \arrow[r, "\phi"] & C^*(E)\phantom{.} \arrow[d, "\pi_\lambda"] \\
				\cK(\ell^2(E_\mu^* s(\mu))) \arrow[r, "\psi"'] & C^*(E_\mu).
			\end{tikzcd}
		\end{equation}
		To this end, let $\nu, \nu' \in E_\mu^*s(\mu)$. Since neither $\nu$ nor $\nu'$ can contain the edge $e_1$, we see that
		\begin{equation}
			\begin{tikzcd}
				1_{C(\bT)} \otimes \delta_{\nu}\delta_{\nu'}^* \arrow[d, mapsto, "\mathrm{ev}_\lambda\otimes \mathrm{id}_\cK"'] \arrow[r, mapsto, "\phi"] & s_{\nu}s_{\nu'}^*\phantom{.} \arrow[d, mapsto, "\pi_\lambda"] \\
				\delta_{\nu}\delta_{\nu'}^* \arrow[r, mapsto, "\psi"'] & s_{\nu}s_{\nu'}^*.
			\end{tikzcd}
		\end{equation}
		This confirms that $\pi_\lambda$ restricts to the standard isomorphism $1 \otimes \cK \cong \cK$. On the other hand, we see that
		\begin{equation}
			\begin{tikzcd}
				z \otimes \delta_{s(\mu)}\delta_{s(\mu)}^* \arrow[d, mapsto, "\mathrm{ev}_\lambda\otimes \mathrm{id}_\cK"'] \arrow[r, mapsto, "\phi"] & s_\mu\phantom{,} \arrow[d, mapsto, "\pi_\lambda"] \\
				\lambda \delta_{s(\mu)}\delta_{s(\mu)}^* \arrow[r, mapsto, "\psi"'] & \lambda p_{s(\mu)},
			\end{tikzcd}
		\end{equation}
		using the fact that $\pi_\lambda(s_\mu)=(\lambda s_{e_n}^*\cdots s_{e_2}^*)(s_{e_2}\cdots s_{e_n})=\lambda p_{s(\mu)}$. Since we have checked this on a set of generators for $C(\bT)\otimes \cK$, we conclude that \eqref{eq:circleidealdiagram} commutes. Thus, to show that $I_\mu \trianglelefteq C^*(E)$ is purely large, we need only verify that $\pi_\lambda(I_\mu)=I_{s(\mu)}\trianglelefteq C^*(E_\mu)$ is essential for all $\lambda \in \bT$. However, this is clear since $s(\mu)$ is the unique source in $E_\mu$, and a path may be found from $s(\mu)$ to any other vertex in $E_\mu$. We conclude by Theorem~\ref{thm:PPV} that $I_\mu$ is purely large in $C^*(E)$.
	\end{proof}

	\section{Finite graphs with no sources}\label{sec:no-sources}
	
	The goal of this section is to prove Theorems~\ref{intro-thm:graphs-with-dim-one} and~\ref{intro-thm:dim-at-most-two} in the special case that $E$ is a graph with no sources.  In particular, for such graphs $E$, we will prove that $\fndim{C^*(E)} \leq 2$.  Since $E$ is a finite graph, $E$ has finitely many gauge-invariant ideals, the number of which will be inducted upon. The induction will require a technical strengthening of faithful nuclear dimension.  In short, given a minimal ideal $I \trianglelefteq C^*(E)$ with $I$ stable and $C^*(E) / I \cong C^*(F)$ for a suitable subgraph $F$ of $E$, we know the extension
	\begin{equation}
		\begin{tikzcd}
			0 \arrow{r} & I \arrow{r} & C^*(E)/I^\perp \arrow{r} & C^*(F)/ I^\perp \arrow{r} & 0
		\end{tikzcd}
	\end{equation}
	is purely large by Theorem~\ref{thm:purelylargeannihilator}.  This allows us to use Lemma~\ref{lem:key-lemma} to lift nuclear dimension approximations for $C^*(F)$ to $C^*(E)$ assuming one color of the approximations for $C^*(F)$ is totally faithful modulo $I^\perp$.  Since we can't specify which ideal of $C^*(E)$ will play the role of $I^\perp$ at the next stage of the induction, we strengthen our notion of faithful nuclear dimension to require that one color of the approximation is totally faithful modulo every gauge-invariant ideal.  
	
	In the following, the role of the group $G$ is only to restrict the collection of ideals considered.  The case of interest here is when $G = \mathbb T$, $A$ is a graph $C^*$-algebra, and the action $\rho\colon G \curvearrowright A$ is the canonical gauge action.
	
	\begin{defi}\label{defi:nicealgebraapproximations}
		Let $G$ be a locally compact group, and let $\cC_G(n)$ be the class of $G$-$C^*$-algebras $(A,\rho)$ with cpc$_\perp$ maps $\phi^0, \phi^1 \colon A \to A_\omega$ such that
		\begin{enumerate}
			\item $\phi^0(a)+\phi^1(a)=a$,
			\item $\phi^0$ preserves all $G$-invariant ideals of $A$ in the sense that $\phi^0(J)\subseteq J_\omega$ for all $G$-invariant ideals $J \trianglelefteq A$,
			\item $\phi^0$ descends to a totally faithful map $\bar\phi^0\colon A/J \to A_\omega/J_\omega$ for all $G$-invariant ideals $J$,
			\item $\nndim{\phi^0}=0$,
			\item $\nndim{\phi^1}\leq n-1$.
		\end{enumerate}
	\end{defi}
	
	Note that if $A$ is a $C^*$-algebra and there is a locally compact group $G$ and an action of $G$ on $A$ such that $A \in \cC_G(n)$, then $\fndim{A} \leq n$, and, in particular, $\ndim{A} \leq n$.  We will show that if $E$ is a finite graph with no sources, then $C^*(E) \in \cC_\mathbb T(2)$, and further, if there is no path from an isolated cycle to another cycle in $E$, then $C^*(E) \in \cC_\mathbb T(1)$.  
	
	Our first goal is to show these classes are invariant under direct sums.

	\begin{prop}\label{prop:classdirectsum}
		Let $G$ be a locally compact group and let $k, m\in \mathbb N$.  Suppose $A \in \cC_G(k)$ and $B \in \cC_G(m)$. Then $(A \oplus B, \rho_A \oplus \rho_B) \in \cC_G(\max\{k,m\})$.
	\end{prop}
	
	\begin{proof}	
		Since $\cC_G(k_1)\subset \cC_G(k_2)$ for $k_1\leq k_2$, without loss of generality we may assume that $k=m$. Let $\phi^i=(\phi^i_n)_{n=1}^\infty$ and $\psi^i=(\psi^i_n)_{n=1}^\infty$ for $i=0,1$ be approximations for $A$ and $B$, respectively, realizing them as algebras in $\cC_G(m)$. Then $\nndim{\phi^0\oplus\psi^0}=0$ and $\nndim{\phi^1\oplus\psi^1}\leq m-1$ by Proposition~\ref{prop:orthogonal-sum}.  It is straightforward to verify that these maps sum to the inclusion $A \oplus B \into (A \oplus B)_\omega$.  The map $\phi^0 \oplus \psi^0$ will preserve $G$-invariant ideals of $A \oplus B$ and is totally faithful modulo each such ideal since every $G$-invariant ideal in $A \oplus B$ has the form $I \oplus J$ for $G$-invariant ideals $I \trianglelefteq A$ and $J \trianglelefteq B$.
	\end{proof}
	
	Compatibility with extensions will take more work. To prepare for this, we need the following lemma to extend the domain of nuclear dimension approximations for an ideal.  Note that the lemma is similar to Lemma~\ref{lem:kirchberg}.  The most important difference is that map $\tilde\theta^1$ satisfies $\fndim{\tilde\theta^1} \leq \fndim{I}$, whereas when $I$ is a stable Kirchberg algebra, we can improve this bound to  $\fndim{\tilde\theta^1} = \fndim{I} - 1 = 0$.  If this improved bound was possible for $I = C(\mathbb T) \otimes \mathcal K$, our proof of Theorem~\ref{intro-thm:dim-at-most-two} will show all graph $C^*$-algebras have nuclear dimension at most one.
	
	\begin{lem}\label{lem:fndimidealapproximations}
		Let $A$ be a separable $C^*$-algebra and let $I \trianglelefteq A$ be an ideal with $\fndim{I} < \infty$.  For every separable set $S \subseteq A_\omega$, there are a positive contraction $f \in I_\omega \cap S'$ such that $fx = x$ for all $x \in I_\omega \cap S$ and cpc$_\perp$ maps $\tilde\theta^0, \tilde\theta^1 \colon A \rightarrow I_\omega$ such that
		\begin{enumerate}
			\item\label{it:extend1} $\tilde\theta^0(a) + \tilde\theta^1(a) = fa$ for all $a \in A$,
			\item\label{it:extend2} $f\tilde{\theta^0}(a) = \tilde\theta^0(a)f = \tilde\theta^0(a)$ for all $a \in A$,
			\item\label{it:extend3} $\nndim{\tilde\theta^0} = 0$ and $\nndim{\tilde\theta^1} \leq \fndim{I} $, and
			\item\label{it:extend4} $\tilde\theta^0|_I$ is totally faithful.
		\end{enumerate}
	\end{lem}
	
	\begin{proof}
		Enlarging $S$ if necessary, we may assume $A \subseteq S$.  Proposition~\ref{prop:qcapproxunit} provides a positive contraction $e \in I_\omega \cap S'$ such that $ex = x$ for all $x \in I_\omega \cap S$.  Fix increasing sequences $(\mathcal F_n)_{n=1}^\infty$, $(\mathcal G_n)_{n=1}^\infty$ and $(\mathcal H_n)_{n=1}^\infty$ of finite subsets of $I$, $A$, and $C_0(0, 1] \otimes I$, respectively, with dense unions.
		
		We will construct cp maps $\theta_n^i \colon I \rightarrow I_\omega$ and $\Theta_n^0 \colon C_0(0, 1] \otimes I \rightarrow I_\omega$ such that
		\begin{enumerate}\setcounter{enumi}{4}
			\item\label{it:extend6} $\|\theta_n^0(x) + \theta_n^1(x) - x\| < 1/n$ for all $x \in e_n^{1/2} \cG_n e_n^{1/2}$, 
			\item\label{it:extend7} $\|\Theta_n^0(z \otimes x) - \theta_n^0(x)\| < 1/n$ for all $x \in \cF_n$, where $z(t) = t$ for $t \in (0, 1]$,
			\item\label{it:extend8} $\|\Theta^0_n(\xi_1\xi_2) - \Theta^0_n(\xi_1)\Theta^0_n(\xi_2)\|< 1/n$ for all $\xi_1, \xi_2 \in \mathcal H_n$,
			\item\label{it:extend9} $\|\Theta^0_n(\xi)\| \geq \|\xi\| - 1/n$ for all $\xi \in \mathcal H_n$,
			\item\label{it:extend10} $\ndim{\theta_n^0} = 0$ and $\ndim{\theta_n^1} \leq \fndim{I} - 1$.
		\end{enumerate}
		Using the definition of faithful nuclear dimension, let $\phi^0, \phi^1 \colon I \rightarrow I_\omega$ be cpc$_\perp$ maps with $\phi^0(x) + \phi^1(x) = x$ for all $x \in I$, $\nndim{\phi^0} = 0$, $\nndim{\phi^1} = \fndim{I} -1$, and such that the induced $^*$-homomorphism $\Phi^0 \colon C_0(0, 1] \otimes I \rightarrow I_\omega$ is faithful.  By definition, $\phi^0$ and $\phi^1$ can be lifted to sequences of cp maps $(\phi_k^0)_{k=1}^\infty$ and $(\phi_k^1)_{k=1}^\infty$ from $I$ to $I$ such that $\ndim{\phi_k^0} = 0$ and $\ndim{\phi_k^1} \leq \fndim{I} - 1$.  By the Choi--Effros theorem (\cite{CE76a}), there is a sequence of cpc maps $(\Phi_k^0 \colon C_0(0, 1] \otimes I \rightarrow I_\omega)_{k=1}^\infty$ lifting $\Phi^0$.  Then $\theta_n^0 = \phi_k^0$, $\theta_n^1 = \phi_k^1$, and $\Theta_n^0 = \Phi_k^0$ will satisfy \ref{it:extend6}--\ref{it:extend10} for some $k \in \mathbb N$.
		
		Let $(e_n)_{n=1}^\infty$ be a sequence of positive contractions in $I$ representing $e$ . For $n \in \mathbb N$ and $i \in \{0,1 \}$, define a cp map
		$\tilde\theta_n^i \colon A \rightarrow I$ by
		\begin{equation}
			\tilde \theta_n^0(a) = \theta_n^0(e_n^{1/2} a e_n^{1/2})
		\end{equation}
		Let $\tilde\theta^0$ denote the map induced by the sequence $(\tilde\theta^0_n)_{n=1}^\infty$, and observe that $\tilde\theta^0$ is cpc since each $\theta^0_n$ is cpc, and has order zero by \ref{it:extend6}. Let $\tilde\Theta^0 \colon C_0(0, 1] \otimes I \rightarrow I_\omega$ be the map induced by the sequence $(\Theta_n)_{n=1}^\infty$.  By \ref{it:extend8} and \ref{it:extend9}, $\Tilde\Theta^0$ is a faithful $^*$-homomorphism. Further, by \ref{it:extend7}, $\tilde\Theta^0(z \otimes x) = \tilde\theta^0(x)=\theta^0(x)$ for all $x \in I$. Hence $\tilde\theta^0|_I$ is totally faithful, which proves \ref{it:extend4}. 
		
		Consider the separable $C^*$-algebra
		\begin{equation}
			T = C^*(S \cup \{e\}\cup \tilde\Theta^0(C_0(0,1]\otimes A)) \subseteq A_\omega.
		\end{equation}
		A second application of Proposition~\ref{prop:qcapproxunit} gives $f\in I_\omega \cap T'$ such that $fx = x$ for all $x \in I_\omega \cap T$, so in particular, $f e = e$ and $f\tilde\theta^0(a)=\tilde\theta^0(a)f=\tilde\theta^0(a)$ for all $a \in A$, proving \ref{it:extend2}.  Now, define
		\begin{equation}
			\tilde\theta_n^1(a) = \theta_n^1(e_n^{1/2} a e_n^{1/2}) + ((f_n - e_n)_+)^{1/2} a ((f_n - e_n)_+)^{1/2}.
		\end{equation}
		Let $\tilde\theta^i$ denote the map induced by the sequence $(\tilde\theta^n_n)_{n=1}^\infty$.  It suffices to prove \ref{it:extend1} and \ref{it:extend3} since \ref{it:extend1} will imply $\tilde\theta^1$ is cpc$_\perp$.
		
		For $n \in \mathbb N$, $a \in \mathcal G_n$, and $m \geq n$, we have
		\begin{equation}
			\theta^0_m(e_m^{1/2}ae_m^{1/2}) + \theta^1_m(e_m^{1/2}ae_m^{1/2}) \approx_{1/m} e_m^{1/2} a e_m^{1/2}.
		\end{equation}
		Since the union of the sets $\cG_n$ is dense in $A$, this implies
		\begin{equation}
			\tilde\theta^0(a) + \tilde\theta^1(a) = e^{1/2} a e^{1/2} + (f - e)_+^{1/2} a (f - e)_+^{1/2} = fa,
		\end{equation}
		so \ref{it:extend1} holds.
		
		Since $\theta_n^0$ has nuclear dimension zero and the map $A \rightarrow I \colon a \mapsto e_n^{1/2} a e_n^{1/2}$ is cpc, the composition $\tilde\theta_n^0$ has nuclear dimension zero.  Also, 
		\begin{equation}
			\nndim{\tilde\theta^1} \leq \ndim{I} \leq \fndim{I}.
		\end{equation}
		Hence \ref{it:extend3} holds.  
	\end{proof}
	
	The following result constructing suitable nuclear dimension approximations for extensions will enable the induction through finite graphs with no sources.
	
	\begin{thm}\label{thm:superniceapproximations}
		Let $G$ be a locally compact group and let $0 \to I \to A \to B \to 0$ be an extension of separable $G$-$C^*$-algebras such that $I$ is stable and has no $G$-invariant ideals.  Suppose the extension $0 \to I \to A/I^\perp \to B / I^\perp \to 0$ is purely large, $\fndim{I}<\infty$, and $B \in \cC_G(m)$ for some $m \in \mathbb N$. Then
		\begin{enumerate}
			\item\label{it:supernice1} $A \in \cC_G(\max\{\fndim{I}+1,m\})$.
			\item\label{it:supernice2} If $I$ is a Kirchberg, $A \in \cC_G(m)$.
		\end{enumerate}
	\end{thm}
	
	\begin{proof}
		We will prove both parts simultaneously.  Let $e \in I_\omega \cap A'$ be such that $e x = x$ for all $x \in I$ and define $S = C^*(A \cup \{e\})$.  Apply Lemma~\ref{lem:fndimidealapproximations} in case \ref{it:supernice1} or Lemma~\ref{lem:kirchberg} in case \ref{it:supernice1} to obtain the $f \in A_\omega \cap S'$ with $fx =x$ for all $x \in A \cap S'$ and the suitable maps $\tilde\theta^0, \tilde\theta^1 \colon A \rightarrow I_\omega$.  We note that the only difference in this cases is $\nndim{\tilde\theta^1} \leq \fndim{I}$ in the first case and $\nndim{\tilde\theta^1} = 0$ in the second case.
		
		Let $\psi^0, \psi^1 \colon B \rightarrow B_\omega$ be the maps realizing $B \in \cC_G(m)$ and let $\tilde\psi^0, \tilde\psi^1 \colon A \rightarrow A_\omega$ be the maps obtained from applying Lemma~\ref{lem:key-lemma}.  Define
		$\phi^i \colon A \rightarrow A_\omega$ by 
		\begin{equation}
			\phi^i(a) = \tilde\theta^i(a) + \tilde\psi^i(a), \qquad a \in A. 
		\end{equation} 
		Then $\phi^0$ and $\phi^1$ are cpc$_\perp$ maps adding to the inclusion $A \rightarrow A_\omega$.
		
		As in the proof of Theorem~\ref{body-thm:extension-result-2}, $\tilde\theta^i$ and $\tilde\psi^i$ have orthogonal ranges for $i \in \{0, 1\}$, so using Proposition~\ref{prop:orthogonal-sum},
		\begin{equation}
			\nndim{\phi^i} = \max \{\nndim{\tilde\theta^i}, \nndim{\tilde \psi^i}\}
		\end{equation}
		So $\nndim{\phi^0} = 0$.  In case \ref{it:supernice2}, $\nndim{\phi^1} \leq m - 1$, and in case \ref{it:supernice1}, 
		\begin{equation}
			\nndim{\phi^1} \leq \max \{\fndim{I}, m - 1\}.
		\end{equation}
		It remains to show $\phi^0$ preserves $G$-invariant ideals and is totally faithful modulo these ideals.
		
		Let $J \trianglelefteq A$ be a $G$-invariant ideal.  Then $IJ \trianglelefteq I$ is a $G$-invariant ideal, so $IJ = 0$ or $IJ = I$.  That is, either $J \subseteq I^\perp$ or $I \subseteq J$.  Letting $q\colon A \to B$ denote the quotient map, note that, in either case, for $a \in J$, 
		\begin{equation}
			q_\omega(\phi^0(a)) = \psi^0(q(a)) \in \psi^0(q(J)) \subseteq q(J)_\omega = q_\omega(J_\omega).
		\end{equation}
		So $\phi^0(a) \in J_\omega + I_\omega$.  When $I \subseteq J$, we have $J_\omega + I_\omega = J_\omega$, so $\phi^0(J) \subseteq J_\omega$.  When $J \subseteq I^\perp$, we have $\phi^0(a) \in \phi^0(I^\perp) \subseteq (I^\perp)_\omega$, so $\phi^0(J) \subseteq (J_\omega + I_\omega) \cap (I^\perp)_\omega = J_\omega$.
		
		We conclude by showing $\phi^0$ is totally faithful modulo each $G$-invariant ideal.  Fix a $G$-invariant ideal $J \trianglelefteq A$.  Let $\Phi^0, \Psi^0, \Theta^0$ denote the $^*$-homomorphisms induced by $\phi^0$, $\psi^0$, and $\theta^0$, and let $Q = \mathrm{id}_{C_0(0,1]} \otimes q$.  If $a \in C_0(0, 1] \otimes A$ such that $\Phi^0(a) \in J_\omega$, then \begin{equation}
			\Psi^0(Q(a)) = q_\omega (\Phi^0(a)) \in q_\omega(J_\omega) = q(J)_\omega.
		\end{equation}
		Since $\Psi^0$ is faithful modulo $q(J)$, we have $Q(a) \in Q(C_0(0, 1] \otimes J)$.  Therefore,
		\begin{equation}
			(\Phi^0)^{-1}(J_\omega) \subseteq C_0(0, 1] \otimes (I + J).
		\end{equation}
		
		When $I \subseteq J$, it follows immediately that $\Phi^0$ is faithful modulo $J$.  Suppose $J \subseteq I^\perp$.  It suffices to show
		\begin{equation}
			(\Phi^0)^{-1}((I^\perp)_\omega) \subseteq C_0(0, 1] \otimes I^\perp
		\end{equation}
		since $(I + J) \cap I^\perp = J$.  That is, it suffices to show the map
		\begin{equation}
			\bar \Phi^0 \colon C_0(0, 1] \otimes (A/I^\perp) \rightarrow (A/I^\perp)_\omega 
		\end{equation}
		is totally faithful.  Recall that $\tilde\psi^0(a) \leq (1-f)a$ for all $a \in A$ and $fx = x$ for all $x \in I$.  Hence $\tilde\psi^0|_I = 0$, so $\bar\phi^0|_I = \theta^0$.  Then $\bar \Phi^0 = \Theta^0$ is faithful on $C_0(0, 1] \otimes I$.  Since $C_0(0, 1] \otimes I$ is an essential ideal of $C_0(0, 1] \otimes (A/I^\perp)$, $\bar\Phi^0$ is totally faithful.
	\end{proof}
	
	With the extension results above and the graph algebra results from the previous section, we are now able to obtain Theorems~\ref{intro-thm:graphs-with-dim-one} and \ref{intro-thm:dim-at-most-two} in the case of finite graphs with no sources.
	
	\begin{thm}\label{thm:nosources}
		If $E$ is a finite graph with no sources, then $C^*(E) \in \cC_\mathbb T(2)$.  Further, if there is no isolated cycle with a path to another cycle in $E$, then $C^*(E) \in \cC_\mathbb T(1)$.
	\end{thm}
	
	\begin{proof}
		We induct on the number of proper gauge-invariant ideals in $C^*(E)$.  In the gauge-simple case, $C^*(E)$ is either a unital Kirchberg algebra or is isomorphic to $M_n(C(\mathbb T))$ for some $n \in \mathbb N$.  In either case, $C^*(E) \in \cC_\bT(1) \subseteq \cC_\bT(2)$.  Indeed, $\fndim{C(\mathbb T)} = 1$ by Theorem~\ref{thm:nndim:abelian}, and hence $C(\mathbb T) \in \cC_\mathbb T(1)$ as $C(\mathbb T)$ has no $\mathbb T$-invariant ideals.  For matrix amplifications, this is Proposition~\ref{prop:fndimstable}.  For a Kirchberg algebra $A$, since $A$ is simple and $\fndim{A} = 1$ by Example~\ref{ex:nndim:purelyinfinite}, we have $A \in \cC_\mathbb T(1)$.
		
		Assume the result holds whenever $C^*(E)$ has at most $n$ gauge-invariant ideals and let $E$ be a graph such that $C^*(E)$ has exactly $n + 1$ gauge-invariant ideals.  Let $I \trianglelefteq C^*(E)$ be a minimal gauge-invariant ideal and note that $C^*(E) / I \cong C^*(F)$ for a subgraph $F$ of $E$.  By the inductive hypotheses, we have $C^*(F) \in \cC_\mathbb T(2)$.  Note that $I$ is either a Kirchberg algebra or is stably isomorphic to $C(\mathbb T)$.  In either case, $I \in \cC_\mathbb T(1)$ by Example~\ref{ex:nndim:purelyinfinite} or Theorem~\ref{thm:nndim:abelian} and Proposition~\ref{prop:fndimstable}.  Further, $I$ is either unital or stable, and in the latter case, the canonical extension
		\begin{equation}
			\begin{tikzcd}
				0 \arrow{r} & I \arrow{r} & C^*(E) / I^\perp \arrow{r} & C^*(F) / I^\perp \arrow{r} & 0.
			\end{tikzcd}
		\end{equation}
		is purely large by Theorem~\ref{thm:purelylargeannihilator}.  If $I$ is unital, $C^*(E) \in \cC_\mathbb T(2)$ by Proposition~\ref{prop:classdirectsum}, and if $I$ is stable, this holds by Theorem~\ref{thm:superniceapproximations}\ref{it:supernice1}.
		
		Now assume no isolated cycle in $E$ with a path to another cycle.  Then the same holds for $F$, so $C^*(F) \in \cC_\mathbb T(1)$ by the alternative inductive hypothesis.  If $I$ is unital, then $C^*(E) \cong I \oplus C^*(F) \in \cC_\mathbb T(1)$ by Proposition~\ref{prop:classdirectsum}.  If $I$ is stable, then $I$ is a Kirchberg algebra by Proposition~\ref{prop:minimalideals}.  Hence Theorem~\ref{thm:superniceapproximations}\ref{it:supernice2} implies $C^*(E) \in \cC_\mathbb T(1)$.
	\end{proof}
	
	\section{Proof of Theorems~\ref{intro-thm:graphs-with-dim-one} and~\ref{intro-thm:dim-at-most-two}}
	
	We now have the machinery needed to prove the main result of the paper.  
	
	\begin{thm}[Theorems~\ref{intro-thm:graphs-with-dim-one} and~\ref{intro-thm:dim-at-most-two}]\label{thm:graphs}
		Let $E$ be  a directed graph.  Then 
		\begin{equation}
			\ndim{C^*(E)} \leq 2.
		\end{equation}
		Further, if $E$ has no path from an isolated cycle to another cycle, then
		\begin{equation}
			\ndim{C^*(E))} = 
			\begin{cases} 
				0 & $E$ \text{ has no cycles}; \\ 1 & $E$ \text{  has a cycle}. 
			\end{cases}
		\end{equation}
	\end{thm}
	
	The proof will be given in stages, which we separate out as lemmas below.
	The case of finite graphs with no sources follows from Theorem~\ref{thm:nosources}, and our next goal is to extend this to all finite graphs.  This will be proven via induction on the ideal generated by a projection corresponding to a source, noting that this ideal is either a direct summand or isomorphic to $\mathcal K$.  We will apply Lemma~\ref{lem:key-lemma} and the structure of the approximations obtained in Theorem~\ref{thm:nosources} to lift the approximations from the quotient along these extensions.
	
	\begin{lem}\label{lem:finite-graphs}
		Theorem~\ref{thm:graphs} holds for all finite graphs $E$.
	\end{lem}
	
	\begin{proof}
		When $E$ has no cycles, $C^*(E)$ is finite dimensional and hence has nuclear dimension zero.  Conversely, if $E$ has a cycle, then $E$ is not an AF algebra by (\cite[Remark~2.2(iii)]{WZ10}; cf.\ \cite[Example~4.1]{KW04}) and hence has nuclear dimension at least one by \cite[Proposition~2.3(iii)]{WZ10}.
		
		Assume $E$ has at least one cycle.  Let $k = 1$ if $E$ has no path from an isolated cycle to another cycle, and let $k = 2$ otherwise.  We'll show that $\ndim{C^*(E)} \leq k$.  When $E$ has no sources, this is immediate from Theorem~\ref{thm:nosources}, so we assume $E$ has at least one source, and we will prove the result by inducting through the sources in $E$.
		
		Let $v_1, \ldots, v_n \in E^0$ be the sources in $E$, and let $H$ be the hereditary saturated set they generate.  Then $I_H$ is the ideal generated by the projections $p_{v_1}, \ldots p_{v_n}$, and $C^*(E) / I_H \cong C^*(E \setminus H)$ as recalled in Theorem~\ref{thm:graph-ideals}. Specifically, $H$ consists of all vertices in $E$ that receive only finitely many paths from $v_1, \ldots, v_n$ and do not receive paths from any vertices on cycles.  Then $E \setminus H$ is a graph with no sources, and if no isolated cycle in $E$ has a path to a cycle, then the same is true for $E\setminus H$.
		
		For $m \in \{0, 1, \ldots, n\}$, let $H_m$ denote the hereditary saturated set generated by $v_{m+1}, \ldots, v_n$, so in particular, $H_0 = H$ and $H_n = \emptyset$.  We will show that for each $m = 0, \ldots n$, there are cpc$_\perp$ maps $\phi^0, \phi^1 \colon C^*(E \setminus H_m) \rightarrow C^*(E \setminus H_m)_\omega$ such that
		\begin{enumerate}
			\item\label{it:finite1} $\phi^0(a) + \phi^1(a) = a$ for all $a \in C^*(E \setminus H_m)$,
			\item\label{it:finite2} for all gauge-invariant ideals of $J \trianglelefteq C^*(E \setminus H_m)$ with $p_{v_{m+1}}, \ldots p_{v_m} \in J$, $\phi^0(J) \subseteq J_\omega$ and the induced map $\bar\phi^0 \colon C^*(E \setminus H_m) / J \rightarrow C^*(E \setminus H_m)_\omega / J_\omega$ is totally faithful, and
			\item\label{it:finite3} $\nndim{\phi^0} = 0$, and $\nndim{\phi^1} \leq k - 1$. 
		\end{enumerate}
		Note that the case $m = 0$ is the statement $C^*(E \setminus H) \in \cC_\mathbb T(k)$, which is true by Theorem~\ref{thm:nosources} when $C^*(E\setminus H)$ is equipped with the canonical gauge action.
		
		Suppose the claim holds for some $m \in \{0, \ldots, n - 1\}$, realized by maps 
		\begin{equation}
			\psi^0, \psi^1 \colon C^*(E \setminus H_m) \rightarrow C^*(E \setminus H_m)_\omega.
		\end{equation}
		Let $I \trianglelefteq C^*(E \setminus H_{m+1})$ denote the ideal generated by $p_{v_{m+1}}$  and note that 
		\begin{equation}
			C^*(E \setminus H_{m+1}) \cong C^*(E \setminus H_m) / I.    
		\end{equation} 
		By Proposition~\ref{prop:source}, either $I$ is finite dimensional or $I \cong \mathcal K$.  If $I$ is finite dimensional, then 
		\begin{equation}
			C^*(E \setminus H_m) \cong I \oplus C^*(E \setminus H_{m+1}),
		\end{equation}
		and the maps $\phi^0 = 0 \oplus \psi^0$ and $\phi^1 = \mathrm{id}_I \oplus \psi^1$ show the claim holds. 
		
		Assume that $I \cong K$ and  consider the extension
		\begin{equation}
			\begin{tikzcd}
				0 \arrow{r} & I \arrow{r} & C^*(E \setminus H_{m+1}) \arrow{r}{q} & C^*(E \setminus H_m) \arrow{r} & 0.
			\end{tikzcd}
		\end{equation}
		By Theorem~\ref{thm:purelylargeannihilator}, the induced extension
		\begin{equation}
			\begin{tikzcd}
				0 \arrow{r} & I \arrow{r} & C^*(E \setminus H_{m+1}) / I^\perp \arrow{r}{\bar q} & C^*(E \setminus H_m) / I^\perp \arrow{r} & 0 
			\end{tikzcd}
		\end{equation}
		is purely large.  Note that $p_{v_1}, \ldots, p_{v_m}$ are orthogonal to $p_{v_{m+1}}$ and hence belong to $I^\perp$.  Also, $I^\perp$ is gauge-invariant since $I$ is, so by assumption, $\bar\psi^0(I^\perp) \subseteq I^\perp$ and the induced map
		\begin{equation}
			\bar\psi^0 \colon C^*(E\setminus H_m)/I^\perp \rightarrow C^*(E \setminus H_m)_\omega / (I^\perp)_\omega
		\end{equation}
		is totally faithful.
		
		Let $f \in I_\omega \cap C^*(E \setminus H_{m+1})'$ be a positive contraction with $fx = x$ for all $x \in I$ as in Proposition~\ref{prop:qcapproxunit}, and let
		\begin{equation}
			\tilde\psi^0, \tilde\psi^1 \colon C^*(E \setminus H_{m+1}) \rightarrow C^*(E \setminus H_{m+1})_\omega
		\end{equation}
		be the maps given by Lemma~\ref{lem:key-lemma}.  Define $\phi^0, \phi^1 \colon C^*(E \setminus H_{m+1}) \rightarrow C^*(E \setminus H_{m+1})_\omega$ by
		\begin{equation}
			\phi^0(a) = \tilde\psi^0(a) \quad \text{and} \quad \phi^1(a) = fa + \tilde\psi^1(a), \qquad a \in A.
		\end{equation}
		Then \ref{it:finite1} holds by construction.  
		
		For \ref{it:finite3}, first note that 
		\begin{equation} 
			\nndim{\tilde\psi^0} \leq \nndim{\psi^0} = 0 \quad \text{and} \quad \nndim{\tilde\psi^1} \leq \nndim{\psi^1} \leq k - 1.
		\end{equation}
		Then $\nndim{\phi^0} = 0$ since $\phi^0 = \tilde\psi^0$.  Since the range of $\tilde\psi^1$ is orthogonal to $f$, Proposition~\ref{prop:orthogonal-sum} implies
		\begin{equation}
			\nndim{\phi^1} \leq \max\{\ndim{I}, \nndim{\tilde\psi^1} \} = k - 1,
		\end{equation}
		which proves \ref{it:finite3} holds.
		
		Now let $J \trianglelefteq C^*(E \setminus H_{m+1})$ be a gauge-invariant ideal containing $p_{v_1}, \ldots, p_{v_{m+1}}$ and note in particular that $I \subseteq J$.  Then $q(J) \trianglelefteq C^*(E \setminus H_m)$ is a gauge-invariant ideal containing  $p_{v_1}, \ldots, p_{v_m}$, so $\psi^0(q(J)) 
		\subseteq q(J)_\omega$ and the induced map $\bar\psi^0$ is faithful; i.e.\ $(\psi^0)^{-1}(q(J)_\omega) = q(J)$.  Since $q_\omega \phi^0 = \psi^0 \circ q$, if $a \in J$, \begin{equation}
			q_\omega(\phi^0(a)) \in \psi^0(q(J)) \subseteq q(J)_\omega = q_\omega(J_\omega),
		\end{equation}
		so $\phi^0(a) \in J_\omega + I_\omega = J_\omega$.  On the other hand, if $a \in C^*(E \setminus H_{m+1})$ and $\phi^0(a) \in J_\omega$, then
		\begin{equation}
			\psi^0(q(a)) = q_\omega(\phi^0(a)) \in q_\omega(J_\omega) = q(J)_\omega,
		\end{equation}
		so $q(a) \in q(J)$.  Since $I \subseteq J$, this implies $a \in J$, which shows \ref{it:finite2} and hence completes the induction.
		
		Now taking $m = n$, we have cpc$_\perp$ maps $\phi^0, \phi^1 \colon C^*(E) \rightarrow C^*(E)_\omega$ satisfying \ref{it:finite1}--\ref{it:finite3}.  In particular, by Lemma~\ref{lem:sum}, the diagonal inclusion $\iota_A \colon A \rightarrow A_\omega$ satisfies
		\begin{equation}
			\nndim{\iota_A} \leq \nndim{\phi^0} + \nndim{\phi^1} + 1 \leq k,
		\end{equation}
		and then Proposition~\ref{prop:nndimisndim} yields $\ndim{C^*(E)} \leq k$.
	\end{proof}

	It is known that the $C^*$-algebra of a countable graph can be approximated by $C^*$-algebras of finite graphs; see \cite[Lemma~1.2]{RS04}, \cite[Proposition~3.1]{JP02}, and \cite[Section~2]{FC24}), for example. In order to do this in a manner compatible with the structure of cycles, we first handle the case of countable row-finite graphs.  The following definition and proposition will aid in the construction of these approximations.
	
	\begin{defi}
		Let $F\subseteq E$ be directed graphs. We say $F$ is \emph{entrance-complete} if whenever $v \in E^0$ and $r^{-1}(v)\cap F^1 \neq \emptyset$, then $r^{-1}(v)\subseteq F^1$.
	\end{defi}
	
	The proof of the following proposition is straightforward; entrance-completeness is used to verify the second Cuntz--Krieger condition \ref{it:CK2} at vertices which receive only finitely many edges.
	
	\begin{prop}
		Let $E$ be a directed graph, and let $F$ be an entrance-complete subgraph of $E$. Then the inclusion $F \subseteq E$ canonically yields an inclusion of $C^*$-algebras $C^*(F)\subseteq C^*(E)$.
	\end{prop}
	
	In the row-finite case, modifying the construction from \cite[Section~2]{FC24} allows us to build such approximations while maintaining the structure of cycles.
	
	\begin{lem}\label{lem:row-finite-graphs}
		Theorem~\ref{thm:graphs} holds for countable row-finite graphs.
	\end{lem}
	
	\begin{proof}
		Recall that $E$ has no cycles if and only if $C^*(E)$ is an AF algebra by \cite[Theorem~2.4]{KPR98}, and further $C^*(E)$ has nuclear dimension zero if and only if $C^*(E)$ is AF by (\cite[Remark~2.2(iii)]{WZ10}; cf.\ \cite[Example~4.1]{KW04}). 
		
		It follows from \cite[Lemma~1.2]{RS04} that $C^*(E)$ is a direct limit of $C^*$-algebras associated to finite graphs.  Combining this with Lemma~\ref{lem:finite-graphs} and the lower semicontinuity of nuclear dimension under inductive limits (\cite[Proposition~2.3(iii)]{WZ10}) shows that $\ndim{C^*(E)} \leq 2$.
		
		Suppose now that $E$ has no path from an isolated cycle to another cycle.  Let $(F'_n)_{n=1}^\infty$ be a sequence of finite subgraphs of $E$ with $\cup_{n=1}^\infty F'_n=E$. We now enlarge the sequence $(F_n')_{n=1}^\infty$ to a sequence of finite subgraphs $(F_n)_{n=1}^\infty$ with $F_n'\subset F_n$ so that, for all $n \in \bN$,
		\begin{enumerate}
			\item\label{it:graphlimit1} $F_n \subset F_{n+1}$,
			\item\label{it:graphlimit2} $F_n$ is entrance-complete, so  $C^*(F_n)\subseteq C^*(F)$, and
			\item\label{it:graphlimit3} isolated cycles in $F_n$ do not have a path to other cycles of $F_n$.
		\end{enumerate}
		
		Assume $n \in \mathbb N$ and $F_{n-1}$ has been constructed (with $F_0$ being the empty graph).  Enlarging $F_n'$ if necessary, we may assume $F_{n-1} \subseteq F_n'$.  Consider the following three sets:  
		\begin{equation}
			\begin{split}
				X_n=&\{v \in (F_n')^0: \text{ $v$ does not lie on a cycle in $E$}\}; \\
				Y_n=&\{v \in (F_n')^0: \text{ $v$ lies on an isolated cycle in $E$}\}; \\
				Z_n=&\{v \in (F_n')^0: \text{ $v$ lies on a non-isolated cycle in $E$}\}.\\
			\end{split}
		\end{equation}
		Since each vertex $v \in Z_n$ lies on a non-isolated cycle in $E$, choose distinct return paths $\mu_v, \nu_v$ in $E$ for $v$. Let $F_n$ be the entrance-completion of $F_n' \cup \{\mu_v,\nu_v: v \in Z_n\}$.
		
		Note that \ref{it:graphlimit1} and \ref{it:graphlimit2} are trivially satisfied, so we proceed to verify \ref{it:graphlimit3}. Let $\mu$ be an isolated cycle in $F_n$. Since the only new vertices added in an entrance-completion are sources, we know $\mu$ must contain vertices in $F_n' \cup \{\mu_v,\nu_v: v \in Z_n\}$. However, following the proof of \cite[Section~2, Lemma~A]{FC24}, each vertex on $\{\mu_v,\nu_v: v \in Z_n\}$ has at least two distinct return paths, so no vertices of $\mu$ can lie on $\mu_v$ or $\nu_v$ for $v \in Z_n$. We conclude that $\mu$ is a cycle in $F_n'$. Since we added multiple return paths to vertices in $Z_n$, and vertices of $X_n$ don't lie on cycles in $E$, we conclude that $\mu$ consists solely of vertices in $Y_n$. But this implies that $\mu$ must be an isolated cycle in $E$; thus, it cannot connect to another cycle in $E$, nor another cycle in $F_n$. We conclude that $F_n$ satisfies \ref{it:graphlimit3}. Proceed inductively to construct a sequence of finite subgraphs $(F_n)$ satisfying \ref{it:graphlimit1}--\ref{it:graphlimit3}.
		
		Now, since the sequence $(F_n)$ is entrance-complete, we have inclusions $C^*(F_n)\into C^*(F_{n+1})$. Furthermore, since $\cup_{n=1}^\infty F_n=E$, we see that $C^*(E)=\varinjlim C^*(F_n)$. By Lemma~\ref{lem:finite-graphs}, $\ndim{C^*(F_n)}\leq 1$ for all $n \in \bN$, and so $\ndim{C^*(E)}\leq 1$ by another application of \cite[Proposition~2.3(iii)]{WZ10}.
	\end{proof}
	
	Drinen and Tomforde (\cite{DT05}) provided a constructive method to realize the $C^*$-algebra of a countable graph as a full corner of the $C^*$-algebra of a countable row-finite graph.  Since nuclear dimension is preserved by passing to full corners, we can extend to countable graphs by noting that their desingularization construction preserves the structure of cycles. We need to recall some properties of this construction in order to verify that our condition on isolated cycles is preserved.
	
	First, recall the desingularization process as follows: If $v_0 \in E^0$ is such that $r^{-1}(v_0)=\{e_i:1\leq i <\infty\}$, we construct a new graph $G$ where
	\begin{equation}
		\begin{aligned}
			G^0=&\phantom{.}E^0 \cup \{v_i:1\leq i < \infty\} \\
			G^1=&\phantom{.}E^1\setminus \{e_i:1\leq i < \infty\} \cup \{f_i, g_i: 1\leq i< \infty\}
		\end{aligned}
	\end{equation} 
	where $r_G$ and $s_G$ agree with $r$ and $s$ on $E^1$; for the new edges, we define
	\begin{equation}
		\begin{aligned}
			r_G(f_i)=&v_i \qquad& s_G(f_i)&=s(e_i) \\
			r_G(g_i)=&v_{i-1} \qquad& s_G(g_i)&=v_i.
		\end{aligned}
	\end{equation}    
	
	A Drinen--Tomforde desingularization for $E$ is a graph $F$ with a set of paths $M$ where this desingularization process has occurred at each infinite receiver in $E$; $M$ is the collection of paths added to each infinite receiver. This produces a row-finite graph $F$ such that $C^*(E)\subset C^*(F)$ as a full corner (\cite[Theorem~2.11]{DT05}). We next recall the map $\phi\colon E^* \to F^*$ defined by Drinen and Tomforde.
	
	\begin{lem}[{\cite[Lemma~2.6(a)]{DT05}}]\label{lem:desing-paths}
		There is a bijective map
		\begin{equation}
			\phi\colon E^* \to \{\alpha \in F^*: s(\alpha), r(\alpha) \in E^0\subset F^0\}
		\end{equation}
		that preserves the range and source, and, when restricted to $E^1$, is given by
		\begin{equation}
			\phi(e)=\begin{cases}
				e & \text{if $s(e)$ is a finite receiver}\\
				\alpha_j & \text{if $e=e_j$ for some $j$}
			\end{cases}
		\end{equation}
		where $\alpha_j\coloneq g_1g_2\cdots g_jf_j.$
	\end{lem}
	
	The case of an arbitrary countable graph essentially follows from the proof of Lemma~2.7 in \cite{DT05} and the fact that nuclear dimension is preserved by stable isomorphism.
	
	\begin{lem}\label{lem:countable-graphs}
		Theorem~\ref{thm:graphs} holds for all countable graphs $E$.
	\end{lem}
	
	\begin{proof} 
		All that we need to show is that the condition on isolated cycles is preserved when taking a Drinen--Tomforde desingularization $(F,M)$ since $C^*(F)\otimes \cK \cong C^*(E)\otimes \cK$ by \cite[Theorem~2.8]{B77} and \cite[Theorem~2.11]{DT05}, from which it follows that $\ndim{C^*(E)}=\ndim{C^*(F)}$ by \cite[Corollary~2.8(i)]{WZ10}. Using the map $\phi$ from Lemma~\ref{lem:desing-paths}, we will show the following two statements:
		\begin{enumerate}
			\item\label{it:desing1} If $\mu$ is an isolated cycle in $F$, then $\phi^{-1}(\mu)$ is an isolated cycle in $E$, and
			\item\label{it:desing2} If $\mu_1$ and $\mu_2$ are two cycles in $F$ that are connected by a path in $F$, then $\phi^{-1}(\mu_1)$ and $\phi^{-1}(\mu_2)$ are two cycles in $E$ that are connected by a path in $E$.
		\end{enumerate}
		Note that any cycle in $F$ necessarily contains vertices from $E$, justifying the use of $\phi^{-1}$ above. Combining \ref{it:desing1} and \ref{it:desing2} shows that if $F$ has an isolated cycle connecting to another cycle, so too does $E$, from which the result follows from Lemma~\ref{lem:row-finite-graphs}.
		
		Proving \ref{it:desing1} is essentially identical to the proof of \cite[Lemma~2.7(b)]{DT05}. Let $\mu \in F$ be an isolated cycle. Then $\phi^{-1}(\mu)$ is also a cycle in $E$. If $\phi^{-1}(\mu)$ were not isolated, a distinct return path $\nu$ could be constructed for a vertex $v$ of $\phi^{-1}(\mu)$, and $\phi(\nu)$ is a return path in $F$ for $v$ as well. By the injectivity of $\phi$, it follows that $\phi(\nu)\neq \mu$, and so $\mu$ would fail to be isolated.
		
		As for \ref{it:desing2}, suppose $\mu_1$ and $\mu_2$ are cycles in $F$, and $\gamma$ is a path with $s(\gamma) \in \mu_1$ and $r(\gamma)\in \mu_2$. By extending $\gamma$ in both directions along $\mu_1$ and $\mu_2$, we may assume that $s(\gamma)$ and $r(\gamma)$ belong to $E^0$. Then $\phi^{-1}(\gamma)$ is a path from $\phi^{-1}(\mu_1)$ to $\phi^{-1}(\mu_2)$. Therefore, if two cycles are connected in $F$, the corresponding cycles must be connected in $E$, completing the proof.
	\end{proof}

	All that remains is removing the countability assumption.  The following separabilization result allows us to reduce to countable graphs.
	
	\begin{prop}\label{lem:reductiontocountablegraph}
		If $E$ is graph and $S \subseteq C^*(E)$ is a separable set, then there is a countable subgraph $F \subseteq E$ such that the inclusion $F \hookrightarrow E$ induces an injective $^*$-homomorphism $C^*(F) \hookrightarrow C^*(E)$ whose range contains $S$.  Further, if no isolated cycle in $E$ has a path to a cycle, then $F$ can be chosen with the same properties.
	\end{prop}
	
	\begin{proof}
		Since $C^*(E)$ is generated by $p_v$ and $s_e$ for $v\in E^0$ and $e \in E^1$ and $S$ is separable, there is a countable subgraph $F_0 \subseteq E$ such that $S$ is contained in the $C^*$-algebra generated by $p_v$ and $s_e$ for $v \in F_0^0$ and $e \in F_0^1$.  
		
		We will inductively construct an increasing sequence of countable subgraphs $(F_n)_{n=1}^\infty$ of $E$ containing $F_0$ such that for each $n \geq 0$ and $v \in F_n^0$,
		\begin{enumerate}
			\item if $r^{-1}(v) \subseteq E^1$ finite, then  $r^{-1}(v) \subseteq F_{n+1}^1$, and 
			\item if $r^{-1}(v) \subseteq E^1$ is infinite, then $r^{-1}(v) \cap F_{n+1}^1$ is infinite
		\end{enumerate} 
		To this end, suppose $F_n$ has been constructed.  Decompose $F_n^0 = X_n \sqcup Y_n$, where
		\begin{equation}
			X_n = \{v \in F_n^0 : |r^{-1}(v)| < \infty \} \quad \text{and} \quad Y_n = \{v \in F_n^0 : |r^{-1}(v)| = \infty \}.
		\end{equation}
		For each $v \in Y_n$, let $Y_{n, v} \subseteq r^{-1}(v)$ be a countably infinite set.
		Then define
		\begin{align}
			F_{n+1}^1 &= F_n^1 \cup \bigcup_{v \in X_n} r^{-1}(v) \cup \bigcup_{v \in Y_n} Y_{n, v}
			\shortintertext{and}
			F_{n+1}^0 &= F_n^0  \cup r(F_{n+1}^1) \cup s(F_{n+1}^1).
		\end{align}
		This completes the induction.
		
		Let $F = \bigcup_{n=1}^\infty F_n$ and note that $F$ is a countable subgraph of $E$ containing $F_0$ and for each $v \in F^0$,
		\begin{enumerate}
			\item[(i')] if $r^{-1}(v) \subseteq E^1$ finite, then  $r^{-1}(v) \subseteq F^1$, and 
			\item[(ii')] if $r^{-1}(v) \subseteq E^1$ is infinite, then $r^{-1}(v) \cap F^1$ is infinite.  
		\end{enumerate}
		Then there is a $^*$-homomorphism $C^*(F) \rightarrow C^*(E)$ given by $p_v \mapsto p_v$ and $s_e \mapsto s_e$ for all $v \in F^0$ and $e \in F^1$.  By the gauge-invariant uniqueness theorem, this map is injective, and by construction, the image contains $S$.
		
		Now suppose $E$ has the property that no isolated cycle has a path to another cycle, and we will explain how to modify the construction above.  Construct the countable graph $F_0$ in the same fashion.  Assuming $F_n$ has been constructed, for each isolated cycle in $F_n$ with a path to a cycle in $E$ has a return path in $E$.  Enlarging $F_{n+1}$ in the above construction, we may arrange for each isolated cycle in $F_n$ with a path to a cycle in $E$ to have a return path in $F_{n+1}$.  Then again let $F = \bigcup_{n=1}^\infty F_n$.  Every cycle in $F$ belongs to some $F_n$ (where $n$ depends upon the choice of cycle).  By construction, this cycle is either not isolated in $F_{n+1}$, and hence is not isolated in $F$, or has no path to a cycle in $E$, and hence has no path to a cycle in $F$.
	\end{proof}
	
	Finally we prove the main result.
	
	\begin{proof}[Proof of Theorem~\ref{thm:graphs}]
		We first note that $\ndim{C^*(E)} \geq 1$ if $E$ has a cycle.  Indeed, let $\mu$ be a cycle. If $\mu$ is entryless, then $p_{s(\mu)} C^*(E) p_{s(\mu)} \cong C(\mathbb T)$.  On the other hand, if $\mu$ has an entry, then the Cuntz--Krieger relations imply $p_{s(\mu)}$ is an infinite projection.  In either case, $A$ is not locally finite dimensional, so $\ndim{C^*(E)} \geq 1$ by \cite[Remark~2.2(iii)]{WZ09}.
		
		Set $k = 0$ if $E$ has no cycles, set $k = 1$ if $E$ has a cycle but no path from an isolated cycle to another cycle, and set $k = 2$ otherwise.  We show $\ndim{C^*(E)} \leq k$.  To this end, let $\mathcal F \subset C^*(E)$ be a finite subset, and let $\epsilon>0$. Apply Lemma~\ref{lem:reductiontocountablegraph} to construct a countable subgraph $F\subseteq E$ such that the inclusion induces an injective $^*$-homomorphism $C^*(F) \into C^*(E)$ whose ranges contains $\cF$.  We will identify $C^*(F)$ as a subalgebra of $C^*(E)$.  Further, if $E$ has a cycle or has no path from an isolated cycle to a cycle, we can arrange $F$ to have the same property.
		
		By Lemma~\ref{lem:countable-graphs}, $\ndim{C^*(F)} \leq k$.  Then there are finite-dimensional $C^*$-algebras $F^0,\ldots, F^k$ and maps $\phi\colon C^*(F)\to \bigoplus_{i=0}^k F^i$ and $\psi\colon \bigoplus_{i=0}^k F^i\to C^*(F)$ such that $\phi$ is cpc, $\restrict{\psi}{F^i}$ is cpc$_\perp$ for $i=0,\dots,k$, and $\norm{\psi(\phi(a))-a}<\epsilon$ for all $a \in \cF$. Using Arveson's extension theorem (e.g.\ Theorem~1.6.1 combined with Proposition~2.2.1 in \cite{BrownOzawa}), let $\tilde \phi \colon C^*(E)\to \bigoplus_{i=0}^k F^i$ be a cpc extension of $\phi \colon C^*(F)\to \bigoplus_{i=0}^k F^i$, and view $\psi$ as a map into $C^*(E)$.  Then $\tilde\phi$ and $\psi$ show $\ndim{C^*(E)} \leq k$.
	\end{proof}

\end{document}